\documentclass[11pt]{amsart} 
\usepackage{amscd}           
\usepackage{amssymb}         
\usepackage[curve,matrix,arrow]{xy}

\newtheorem{thm}{Theorem}[section]
\newtheorem{lem}[thm]{Lemma}

\newtheorem{cor}[thm]{Corollary}

\newtheorem{prop}[thm]{Proposition}

\theoremstyle{definition}

\renewcommand{\thecase}{}

\newtheorem{defn}[thm]{Definition}

\newtheorem{rmk}[thm]{Remark} 
\renewcommand{\thestep}{}

\theoremstyle{remark}

\makeatletter
\def\alphenumi{
  \def\theenumi{\alph{enumi}}
  \def\p@enumi{\theenumi}
  \def\labelenumi{(\@alph\c@enumi)}}
\makeatother

\makeatletter
\def\thecase{\@arabic\c@case}
\makeatother

\numberwithin{equation}{section}

\makeatletter
\def\thestep{\@arabic\c@step}
\makeatother

\newcommand\embed{\hookrightarrow}
\newcommand\too{\longrightarrow}

\newcommand\barbe{{\bar\beta}}

\newcommand\barM{{\bar{M}}}

\newcommand\barsO{{\bar{\mathcal{O}}}}

\newcommand\barZ{{\bar{Z}}}

\newcommand\ubarRR{{\underline{\mathbb{R}}}}

\newcommand\CC{\mathbb{C}}

\newcommand\LL{\mathbb{L}}

\newcommand\NN{\mathbb{N}}

\newcommand\PP{\mathbb{P}}

\newcommand\RR{\mathbb{R}}

\newcommand\ZZ{\mathbb{Z}}

\newcommand\bdelta{{\boldsymbol{\delta}}}
\newcommand\bDelta{{\boldsymbol{\Delta}}}

\newcommand\bga{{\boldsymbol{\gamma}}}
\newcommand\bgamma{{\boldsymbol{\gamma}}}

\newcommand\bl{{\boldsymbol{\ell}}}

\newcommand\bvarphi{{\boldsymbol{\varphi}}}

\newcommand\bsD{{\boldsymbol{\mathcal{D}}}}
\newcommand\btau{{\boldsymbol{\tau}}}

\newcommand\bd{{\mathbf{d}}}
\newcommand\bD{{\mathbf{D}}}

\newcommand\bE{{\mathbf{E}}}
\newcommand\boldf{{\mathbf{f}}}
\newcommand\bF{{\mathbf{F}}}

\newcommand\bg{{\mathbf{g}}}

\newcommand\bh{{\mathbf{h}}}
\newcommand\bH{{\mathbf{H}}}

\newcommand\bL{{\mathbf{L}}}

\newcommand\br{{\mathbf{r}}}

\newcommand\bS{{\mathbf{S}}}

\newcommand\bV{{\mathbf{V}}}

\newcommand\bW{{\mathbf{W}}}
\newcommand\bx{{\mathbf{x}}}

\newcommand{\cov}{\nabla}

\newcommand{\rd}{\partial}

\newcommand\thalf{{\textstyle{\frac{1}{2}}}}
\newcommand\tquarter{{\textstyle{\frac{1}{4}}}}
\newcommand\tthreehalf{{\textstyle{\frac{3}{2}}}}

\newcommand\half{{{\frac{1}{2}}}}

\newcommand\sixtyfourth{{{\frac{1}{{64}}}}}

\newcommand\fg{{\mathfrak{g}}}

\newcommand\fH{{\mathfrak{H}}}

\newcommand\fM{{\mathfrak{M}}}
\newcommand\fN{{\mathfrak{N}}}

\newcommand\fs{{\mathfrak{s}}}
\newcommand\fS{{\mathfrak{S}}}
\newcommand\ft{{\mathfrak{t}}}
\newcommand\fT{{\mathfrak{T}}}

\newcommand\fV{{\mathfrak{V}}}
\newcommand\fW{{\mathfrak{W}}}

\newcommand\al{\alpha}
\newcommand\be{\beta}
\newcommand\De{\Delta}
\newcommand\de{\delta}
\newcommand\eps{\varepsilon}

\newcommand\La{\Lambda}
\newcommand\ka{\kappa}
\newcommand\om{\omega}
\newcommand\Om{\Omega}
\newcommand\si{\sigma}

\newcommand\hatA{{\hat A}}

\newcommand\hatsC{{\hat{\mathcal{C}}}}

\newcommand\hatsG{{\hat{\mathcal{G}}}}

\newcommand\hatu{{\hat u}}

\newcommand\gl{{\mathfrak{g}\mathfrak{l}}}
\newcommand\fsl{{\mathfrak{s}\mathfrak{l}}}
\newcommand\so{{\mathfrak{s}\mathfrak{o}}}

\newcommand\su{{\mathfrak{s}\mathfrak{u}}}
\newcommand\fu{{\mathfrak{u}}}

\newcommand\GL{\operatorname{GL}}
\newcommand\Or{\operatorname{O}}

\newcommand\PU{\operatorname{PU}}

\newcommand\SO{\operatorname{SO}}

\newcommand\U{\operatorname{U}}

\newcommand\less{\setminus}

\newcommand{\8}{\infty}

\newcommand\ad{{\operatorname{ad}}}

\newcommand\Aut{\operatorname{Aut}}

\newcommand\ch{\operatorname{ch}}
\newcommand\Cl{\operatorname{C\ell}}
\newcommand\CCl{\operatorname{{\mathbb{C}\ell}}}

\newcommand\Coker{\operatorname{Coker}}

\newcommand\End{\operatorname{End}}

\newcommand\Fr{\operatorname{Fr}}

\newcommand\Hom{\operatorname{Hom}}

\newcommand\ind{\operatorname{Index}}

\newcommand\Imag{\operatorname{Im}}

\newcommand\Jac{\operatorname{Jac}}
\newcommand\Ker{\operatorname{Ker}}

\newcommand\loc{\operatorname{loc}}
\newcommand\Map{\operatorname{Map}}

\newcommand\PD{\operatorname{PD}}

\newcommand\Real{\operatorname{Re}}

\newcommand\Ran{\operatorname{Ran}}
\newcommand\rank{\operatorname{rank}}

\newcommand\Stab{\operatorname{Stab}}

\newcommand\Sym{\operatorname{Sym}}
\newcommand\Tor{\operatorname{Tor}}
\newcommand\tr{\operatorname{tr}}
\newcommand\Tr{\operatorname{Tr}}

\newcommand\vol{\operatorname{vol}}

\newcommand\Ad{{\mathrm{Ad}\,}}

\newcommand\id{{\mathrm{id}}}

\newcommand\spin{\text{spin}}
\newcommand\spinc{\text{$\text{spin}^c$ }}
\newcommand\spinu{\text{$\text{spin}^u$ }}
\newcommand\Spinc{\text{$\text{Spin}^c$}}
\newcommand\Spinu{\text{$\text{Spin}^u$}}

\newcommand\sA{{\mathcal{A}}}
\newcommand\sB{{\mathcal{B}}}
\newcommand\sC{{\mathcal{C}}}
\newcommand\sD{{\mathcal{D}}}
\newcommand\sE{{\mathcal{E}}}

\newcommand\sG{{\mathcal{G}}}

\newcommand\sM{{\mathcal{M}}}

\newcommand\sO{{\mathcal{O}}}
\newcommand\sP{{\mathcal{P}}}

\newcommand\sU{{\mathcal{U}}}

\newcommand\tsC{{\tilde\sC}}

\newcommand\tM{{\tilde M}}
\newcommand\tN{{\tilde N}}

\newcommand\tU{{\tilde U}}

\begin{document}
\title[PU(2) Monopoles and Links]
{PU(2) Monopoles and Links of Top-level Seiberg-Witten Moduli Spaces} 
\author[Paul M. N. Feehan]{Paul M. N. Feehan}
\address{Department of Mathematics\\
Ohio State University\\
Columbus, OH 43210}
\email{feehan@math.ohio-state.edu}
\urladdr{http://www.math.ohio-state.edu/$\sim$feehan/} 
\author[Thomas G. Leness]{Thomas G. Leness}
\address{Department of Mathematics\\
Florida International University\\
Miami, FL 33199}
\email{lenesst@fiu.edu}
\urladdr{http://www.fiu.edu/$\sim$lenesst/} 
\date{July 31, 2000, \texttt{math.DG/0007190}.
First version: December 8, 1997, \texttt{dg-ga/9712005 (v1)}, \S 1-3.}
\thanks{The first author was supported in part by an NSF Mathematical 
Sciences Postdoctoral Fellowship under grant DMS 9306061 and by NSF grant
DMS 9704174} 

\begin{abstract}
  This is the first of two articles in which we give a proof---for a broad
  class of four-manifolds---of Witten's conjecture that the Donaldson and
  Seiberg-Witten series coincide, at least through terms of degree less
  than or equal to $c-2$, where $c=-\frac{1}{4}(7\chi+11\sigma)$ and
  $\chi$ and $\sigma$ are the Euler characteristic and signature of the
  four-manifold. In the present article, we construct virtual normal
  bundles for the Seiberg-Witten strata of the moduli space of PU(2)
  monopoles and compute their Chern classes.
\end{abstract}

\maketitle


\section{Introduction}
\label{sec:Intro2}
\subsection{Main results}
\label{subsec:Statement}
The purpose of the present article and its sequel \cite{FL2b}, is to prove
that Witten's conjecture \cite{Witten} relating the Donaldson and
Seiberg-Witten invariants holds in ``low degrees'' for a broad class of
four-manifolds, using the $\PU(2)$-monopole cobordism \cite{PTLocal}.  We
shall assume throughout that $X$ is a closed, connected, smooth
four-manifold with an orientation for which $b_2^+(X)>0$. We state the
simplest version of our result here; more general results are given in
\cite[\S 1]{FL2b}. The Seiberg-Witten (SW) invariants (see
\cite[\S 4.1]{FL2b}) comprise a function, $SW_X:\Spinc(X)\to\ZZ$, where
$\Spinc(X)$ is the set of isomorphism classes of \spinc structures on
$X$. For $w\in H^2(X;\ZZ)$, one defines
\begin{equation}
\label{eq:SWSeries}
\bS\bW_X^{w}(h)
=
\sum_{\fs \in \Spinc(X)}(-1)^{\half(w^{2}+c_{1}(\fs)\cdot w)}
SW_X(\fs)e^{\langle c_{1}(\fs),h\rangle},
\quad
h\in H_2(X;\RR),
\end{equation}
by analogy with the structure of the Donaldson series $\bD_X^w(h)$
\cite[Theorem 1.7]{KMStructure}. There is a map $c_1:\Spinc(X)\to
H^2(X;\ZZ)$ and the image of the support of $SW_X$ is the set $B$ of
SW-basic classes \cite{Witten}. A four-manifold $X$ has SW-simple type when
$b_1(X)=0$ if $c_1(\fs)^2=2\chi+3\sigma$ for all $c_1(\fs)\in B$,
where $\chi$ and $\sigma$ are the Euler characteristic and signature of
$X$.  We let $B^\perp\subset H^2(X;\ZZ)$ denote the orthogonal complement
of $B$ with respect to the intersection form $Q_X$ on $H^2(X;\ZZ)$. Let
$c(X)=-\frac{1}{4}(7\chi+11\sigma)$.

\begin{thm}
\label{thm:DSWSeriesRelation}
Let $X$ be four-manifold with $b_1(X)=0$ and odd $b_2^+(X)\geq 3$. Assume
$X$ is abundant, SW-simple type, and effective. For any $\La\in B^\perp$
and $w\in H^2(X;\ZZ)$ for which $\La^2=2-(\chi+\si)$ and $w-\Lambda\equiv
w_2(X)\pmod{2}$, and any $h\in H_2(X;\RR)$, one has
\begin{equation}
\label{eq:LowDegreeEquality}
\begin{aligned}
\bD^{w}_X(h) 
&\equiv 0\equiv \bS\bW^{w}_X(h) \pmod{h^{c(X)-2}},
\\
\bD^{w}_X(h)
&\equiv 
2^{2-c(X)}e^{\frac{1}{2}Q_X(h,h)}\bS\bW^{w}_X(h) 
\pmod{h^{c(X)}}.
\end{aligned}
\end{equation}
\end{thm}

The vanishing assertion in low degrees for the series $\bD^{w}_X(h)$ and
$\bS\bW^{w}_X(h)$ in equation \eqref{eq:LowDegreeEquality} was proved in
\cite{FKLM}.

We shall explain below the terminology and notation in the statement of
Theorem \ref{thm:DSWSeriesRelation}.  Witten's conjecture \cite{Witten}
asserts that a four-manifold $X$ with $b_1(X)=0$ and odd $b_2^+(X)\geq 3$
has SW-simple type if and only if it has ``KM-simple type'', that is,
simple type in the sense of Kronheimer and Mrowka (see Definition 1.4 in
\cite{KMStructure}), and that the SW-basic and ``KM-basic'' classes (see
Theorem 1.7 in \cite{KMStructure}) coincide; if $X$ has simple type, then
\begin{equation}
\label{eq:WConjecture}
\bD^{w}_X(h) 
=
2^{2-c(X)}e^{\frac{1}{2}Q_X(h,h)}\bS\bW^{w}_X(h),
\quad
h\in H_2(X;\RR).
\end{equation}
The quantum field theory argument giving equation
\eqref{eq:WConjecture} when $b_2^+(X)\ge 3$ has been extended by Moore and
Witten \cite{MooreWitten} to allow $b_2^+(X) \ge 1$, $b_1(X)\ge 0$, and
four-manifolds $X$ of non-simple type. Recall that $b_2^+(X)$ is the
dimension of a maximal positive-definite linear subspace $H^{2,+}(X;\RR)$
for the intersection pairing $Q_X$ on $H^2(X;\RR)$.

With our stated hypotheses, Theorem \ref{thm:DSWSeriesRelation} therefore
proves that equation \eqref{eq:WConjecture} holds mod $h^{c(X)}$ where
$\delta=c(X)$. To prove the equivalence mod $h^{\delta}$ for $\delta>c(X)$,
more work is required. For example, in \cite{FLLevelOne}, we use the gluing
theory of \cite{FL3}, \cite{FL4}---allowing ``one bubble''--- to prove that
equation \eqref{eq:WConjecture} holds mod $h^{c(X)+2}$ under hypotheses
similar to those of Theorem \ref{thm:DSWSeriesRelation}. If one desires a
mod $h^\delta$ relation such as \eqref{eq:LowDegreeEquality} for larger
values of $\delta$ relative to $c(X)$, one must allow ``more bubbles'' and
the difficulty of the calculations rapidly increases: see \cite[\S 1]{FL2b}
for a more detailed discussion.

The concept of ``abundance'' was first introduced in \cite{FKLM}:

\begin{defn}
  \cite[p. 169]{FKLM} We say that a closed, oriented four-manifold $X$ is
  {\em abundant\/} if the restriction of the intersection form to
  $B^{\bot}$ contains a hyperbolic sublattice.
\end{defn}

The abundance condition is merely a convenient way of formulating the
weaker, but more technical condition that one can find (for example)
classes $\Lambda_j\in B^\perp$ such that $\Lambda_j^2=2j-(\chi+\sigma)$,
for $j=1,2,3$: this is the only property of $Q_X|B^\perp$ which we use to
prove Theorem \ref{thm:DSWSeriesRelation}, while slightly different choices
of classes $\Lambda$ with even squares are used to prove the main results
(Theorems 1.1 and 1.3) of \cite{FKLM}. While all compact, complex
algebraic, simply-connected surfaces with $b_2^+\geq 3$ are abundant (see
Theorem \ref{thm:Abundance} and its proof in Appendix \ref{sec:Abundance}),
some examples of manifolds which are not abundant but for which one can
still find classes in $B^\perp$ with the desired even squares are described
in \cite[\S 2]{FKLM}. It remains an interesting problem to determine whether
all smooth four-manifolds have this property, whether or not they are
abundant.

To see that such classes $\Lambda_j$ exist when $X$ is abundant, note that
because $Q_X|B^\perp$ contains a hyperbolic factor there are classes
$e_0,e_1\in H^2(X;\ZZ)$ which are orthogonal to the SW-basic classes and
which satisfy $e_0^2=0=e_1^2$ and $e_0\cdot e_1=1$.  Now set
$t=\frac{1}{2}(\chi+\sigma)$ and define $\La_j\in H^2(X;\ZZ)$ by $\La=e_0
+(j-t)e_1$, and observe that $\Lambda_j^2$ has the desired values for
$j=1,2,3$.

In the present article and its companion \cite{FL2b} we prove Theorem
\ref{thm:DSWSeriesRelation} using the moduli space $\sM_{\ft}$ of $\PU(2)$
monopoles \cite{PTLocal} to provide a cobordism between the link of the
moduli space $M_\kappa^w$ of anti-self-dual connections and links of moduli
spaces of Seiberg-Witten moduli spaces, $M_{\fs}$, these moduli spaces
being (topologically) embedded in $\sM_{\ft}$. Let $\bar\sM_{\ft}$ denote
the Uhlenbeck compactification (see Theorem \ref{thm:Compactness}) of
$\sM_{\ft}$ in the space of ideal $\PU(2)$ monopoles,
$\cup_{\ell=0}^\8(\sM_{\ft_\ell}\times\Sym^\ell(X))$.

\begin{defn}
\label{defn:Effective}
  We say that a closed, oriented, smooth four-manifold $X$ with $b_1(X)\geq
  0$ and $b_2^+(X)\geq 1$ is {\em effective\/} if $X$ satisfies Conjecture
  3.1 in \cite{FKLM}.  This conjecture asserts that for a Seiberg-Witten
  moduli space $M_{\fs}$ appearing in level $\ell\geq 0$ of $\bar\sM_{\ft}$,
  the pairings of products of
  Donaldson-type cohomology classes on the top stratum of $\sM_{\ft}/S^1$
  with a link of $(M_{\fs}\times\Sym^\ell(X))\cap\bar\sM_{\ft}$ in
  $\bar\sM_{\ft}$ are multiples of
  the Seiberg-Witten invariants for $M_{\fs}$. In particular, these
  pairings are zero when the Seiberg-Witten invariants for $M_{\fs}$ are
  trivial. 
\end{defn}

The motivation for Conjecture 3.1 of \cite{FKLM} and a more detailed
explanation appears in \cite[\S 3.1]{FKLM}; see also
\cite{FLGeorgia}. It is almost certainly true that
this conjecture holds for all four-manifolds, based on our work in
\cite{FL3}, \cite{FL4}, \cite{FL5}.  We verify Conjecture 3.1 in
\cite{FKLM} by direct calculation in the present pair of articles
for Donaldson invariants defined by $M_\kappa^w\embed \sM_{\ft}$
and Seiberg-Witten moduli spaces $M_{\fs}\embed \sM_{\ft}$, while
in \cite{FLLevelOne} we verify the conjecture for Seiberg-Witten
moduli spaces $M_{\fs}$ contained in the ``first level'' of the
Uhlenbeck compactification $\bar\sM_{\ft}$. However, we strongly
expect the conjecture to hold for Seiberg-Witten moduli spaces
$M_{\fs}$ contained in any level of the compactification
$\bar\sM_{\ft}$.

Theorem \ref{thm:DSWSeriesRelation} gives a relation, mod $h^{\delta}$,
between the Donaldson and Seiberg-Witten series when $\delta=c(X)$. Now if
$\tilde X=X\#\overline{\CC\PP}^2$ is the blow-up of $X$, then $c(\tilde
X)=c(X)+1$.  However, this does not imply that one can compute a term of
higher degree in $h$ in the Donaldson series $\bD_X^w(h)$ from Theorem
\ref{thm:DSWSeriesRelation} merely by blowing up.
To understand why, it suffices to examine the case where $X$ also
has KM-simple type.  Let $\PD[e]\in H^2(\tilde X;\ZZ)$ denote
the Poincar\'e dual of the exceptional class $e\in H_2(\tilde X;\ZZ)$. Then
Proposition 1.9 in \cite{KMStructure} implies that for all  
$h\in H_2(\tilde X;\ZZ)\cong \RR[e]\oplus H_2(X;\ZZ)$,
$$
\bD^{w+\PD[e]}_{\tilde X}(h) 
= 
-\bD^{w}_X(h)\exp(-\textstyle{\frac{1}{2}}e\cdot e)\sinh(e\cdot h).
$$
Thus, $\bD^{w+\PD[e]}_{\tilde X}$ has a zero at the origin of order one
greater than that of $\bD^{w}_X$.

On the other hand, If $\Lambda\in B^\perp(X)$, then we also have
$\Lambda\in B^\perp(\tilde X)$ and $w+\PD[e]-\Lambda \equiv
w_2(X)+\PD[e]\equiv w_2(\tilde X)\pmod{2}$, while $\chi(\tilde X) +
\sigma(\tilde X) = \chi(X)+\sigma(X)$, so $\Lambda^2$ still has the desired
even square in the case of $\tilde X$. If $X$ is SW-simple type and
abundant then the same holds for $\tilde X$ (and as mentioned, it is almost
certainly true that all four-manifolds are effective). Theorem
\ref{thm:DSWSeriesRelation} then yields, for all $h\in H_2(\tilde X;\ZZ)$,
\begin{align*}
\bD^{w+\PD[e]}_{\tilde X}(h) 
&\equiv 0\equiv \bS\bW^{w+\PD[e]}_{\tilde X}(h) \pmod{h^{c(X)-1}},
\\
\bD^{w+\PD[e]}_{\tilde X}(h)
&\equiv 
2^{2-c(X)-1}\exp\left(\textstyle{\frac{1}{2}}Q_{\tilde X}(h,h)\right)
\bS\bW^{w+\PD[e]}_{\tilde X}(h) 
\pmod{h^{c(X)+1}}.
\end{align*}
One sees that the formula does not give new higher-degree terms in the
Donaldson series.

As we mentioned just prior to its statement, Theorem
\ref{thm:DSWSeriesRelation} is the simplest formulation of the main results
of the present article and its companion \cite{FL2b}. More general results
are described in \cite[\S 1]{FL2b}, where we allow four-manifolds with
$b_2^+(X)\geq 1$ and drop the assumptions that $X$ has SW-simple type or is
abundant and relax the constraints on the existence of classes $\Lambda\in
B^\perp$ with prescribed squares; we also consider some limited cases where
$b_1(X)>0$. A full account of the case $b_1(X)>0$ is of potential interest
but is beyond the scope of the present pair of articles: we shall describe
this case elsewhere.

\subsection{An outline of the proof of Theorem \ref{thm:DSWSeriesRelation}}
The proof of Theorem \ref{thm:DSWSeriesRelation} splits into two steps.
Step (i), which we carry out in this article, is to construct links of
Seiberg-Witten moduli spaces in the top level of the Uhlenbeck-compactified
moduli space of $\PU(2)$ monopoles, as boundaries of tubular neighborhoods
in certain ``thickened'' or ``virtual'' moduli spaces of $\PU(2)$ monopoles
(see Theorems \ref{thm:DefnOfStabilizeBundle} and
\ref{thm:ThickenedModuliSpace}) and then compute the Chern character and
Chern classes of the vector bundles defining these tubular neighborhoods
(see Theorem \ref{thm:ChernCharacterOfNormal} and Corollary
\ref{cor:SimpleNormalChernClass}).

Step (ii), which we take up in the companion article \cite{FL2b}, is to
compute the pairings of products of cohomology classes on the moduli space
of $\PU(2)$ monopoles with the links of the anti-self-dual moduli space of
$\SO(3)$ connections and with links of the moduli spaces of Seiberg-Witten
monopoles. These computations rely on our calculation of the Chern
characters of the normal bundles of the strata of Seiberg-Witten monopoles,
and a comparison of the orientations of the moduli spaces of anti-self-dual
connections and Seiberg-Witten monopoles, and their links in the moduli
space of $\PU(2)$ monopoles. Applying the $\PU(2)$-monopole cobordism then
yields an expression for the Donaldson invariants in terms of
Seiberg-Witten invariants and hence completes the proof of Theorem
\ref{thm:DSWSeriesRelation}.

\subsection{A guide to the article}
An index of notation appears just before \S \ref{sec:Prelim}. The
present article is a revision of sections 1--3 of the preprint
\cite{FL2}, while the companion article \cite{FL2b} is a revision of
sections 4--7 of \cite{FL2}, which was distributed in December 1997.

Section \ref{sec:Prelim} of this paper gathers together the principal gauge
theory results developed in \cite{FeehanGenericMetric}, \cite{FL1} that we
shall need here for the moduli space of $\PU(2)$ monopoles $\sM_{\ft}$. In
\S \ref{subsec:NonAbelianMonopoles}, we review the construction of the moduli
spaces of $\PU(2)$ monopoles $\sM_{\ft}$ from \cite{FL1}, but now phrased
in the convenient and more compact framework of ``\spinu structures''
$\ft$. In \S
\ref{subsec:CompactnessTransversality} we recall our Uhlenbeck compactness
and transversality results for the moduli spaces of $\PU(2)$ monopoles from
\cite{FL1} and \cite{FeehanGenericMetric}.  In \S \ref{subsec:PertSW}, we
recall the construction of the moduli space of Seiberg-Witten monopoles
$M_{\fs}$ as in \cite{KMThom}, \cite{KMContact}, \cite{MorganSWNotes}, but
with non-standard perturbations so we can directly identify these moduli
spaces with strata of reducible $\PU(2)$ monopoles.  Our transversality
result
\cite{FeehanGenericMetric} ensures that the natural stratification of the
Uhlenbeck-compactified moduli space of $\PU(2)$ monopoles is smooth.  In \S
\ref{subsec:AbelianCohom} we compute the cohomology ring of the
configuration space of \spinc pairs and describe the cohomology classes
arising in the definition of Seiberg-Witten invariants and in our later
calculation (see \S \ref{subsec:TopOfNormal}) of the Chern character of
certain universal families of vector bundles over $X$ parameterized by
$M_{\fs}$.

The construction of the links in $\sM_\ft$
of the stratum $M^w_{\ka}$ of anti-self-dual
connections and of the strata $M_{\fs}$ of {\em reducible\/} (or {\em
  Seiberg-Witten\/}) $\PU(2)$ monopoles occupies \S
\ref{sec:Singularities}. In \S \ref{subsec:TypesOfSing} we classify the
possible singularities of $\sM_{\ft}$. A  link $\bL_{\ft,\kappa}^{w}$
of the stratum $M^w_{\ka}$ of anti-self-dual connections in $\sM_{\ft}$ is
constructed in \S \ref{subsec:ASDLink} using the $L^2$ norm of the spinor
components of $\PU(2)$ monopoles to define the distance to the stratum
$M^w_{\ka}$.  In \S \ref{subsec:ReducibleSingularities} we show that the
subspaces of reducible $\PU(2)$ monopoles in $\sM_{\ft}$ can be identified
with Seiberg-Witten moduli spaces $M_{\fs}$, as defined in \S
\ref{subsec:PertSW}.  As we explain in \S
\ref{subsec:LinkOfReduciblesLocal}, the elliptic deformation complex for
the $\PU(2)$ monopole equations \eqref{eq:PT} at a reducible pair splits
into a {\em tangential deformation complex\/} --- which can be identified
with the elliptic deformation complex for the Seiberg-Witten monopole
equations --- and a {\em normal deformation complex\/}. Via the Kuranishi
model, these deformation complexes describe the local structure of the
moduli space $\sM_{\ft}$ of $\PU(2)$ monopoles near a reducible solution.
The description of a neighborhood of $M_{\fs}$ in $\sM_{\ft}$ is
complicated by the fact that we may have both a tangential deformation
complex with positive index, so $\dim M_{\fs} > 0$, and a normal
deformation complex with negative index. (This contrasts with the simpler
situation considered in \cite{DonApplic}, \cite[pp. 65--69]{FU}, where
abstract perturbations are used in conjunction with the local Kuranishi
model to describe the local structure of the moduli space $M^w_{\ka}$ near
an isolated reducible connection when the four-manifold has $b_2^+(X) = 0$:
in this case the normal deformation complex has positive index.)  The
stratum $M_{\fs}$ will not in general be a smooth submanifold of
$\sM_{\ft}$ as reducible $\PU(2)$ monopoles cannot be shown to be regular
points of the zero locus of the $\PU(2)$ monopole equations \eqref{eq:PT}.
Therefore, in \S \ref{subsec:LinkOfReduciblesGlobal}, we construct an
ambient finite-dimensional, smooth manifold $\sM_{\ft}(\Xi,\fs)$ containing
an open neighborhood in $\sM_{\ft}$ of the stratum $M_{\fs}$ and containing
$M_{\fs}$ as a smooth submanifold. We can then define a link
$\bL_{\ft,\fs}$ of $M_{\fs}$ in $\sM_{\ft}/S^1$ as the $S^1$ quotient of
the intersection with $\sM_{\ft}/S^1$ of the boundary of a tubular
neighborhood of $M_{\fs}$ in the ambient manifold $\sM_{\ft}(\Xi,\fs)$. In
order to compute the pairings of cohomology classes on
$\sM_{\ft}^{*,0}/S^1$ (the smooth locus or top stratum of $\sM_{\ft}/S^1$)
with the link $\bL_{\ft,\fs}$ we shall need the Chern
classes of the normal bundle $N_{\ft}(\Xi,\fs)$ of the stratum $M_{\fs}\embed
\sM_{\ft}(\Xi,\fs)$: we accomplish this in \S \ref{subsec:TopOfNormal},
using the Atiyah-Singer index theorem for families, by computing the Chern
character of this normal bundle (Theorem \ref{thm:ChernCharacterOfNormal})
and then, after imposing a constraint on $H^1(X;\ZZ)$ to simplify our
calculations, its Chern classes (Corollary
\ref{cor:SimpleNormalChernClass}). 

In Appendix \ref{sec:Abundance}, we include a proof that all compact, complex
algebraic, simply connected surfaces with $b_2^+\geq 3$ are abundant
(Theorem \ref{thm:Abundance}). For minimal surfaces of general type, this
fact was asserted in \cite[p. 175]{FKLM}.

\subsubsection*{Acknowledgments}
The authors thank Adebisi Agboola, Ron Fintushel, Tom Mrowka, Peter
Ozsv\'ath, Andr\'as Stipsicz, and Zolt\'an Szab\'o for helpful
conversations during the course of our work on this article and its
companion \cite{FL2b}. We especially thank Tom Mrowka for his many helpful
comments during the course of this work, for bringing a correction to some
examples in \cite{FL2} to our attention, as well pointing out that our
Theorem 1.4 in \cite{FL2} could be more usefully rephrased and specialized to
give the version stated as Theorem 1.4 in \cite{FL2b} (and Theorem 2.1 in
\cite{FKLM}). We are also very 
grateful to Andr\'as Stipsicz for his considerable help with the proof of
Theorem \ref{thm:Abundance} (the proofs of the difficult cases are all due
to him), as well as to Adebisi Agboola for his assistance with the number
theoretical aspects of that proof. We thank the Columbia University
Mathematics Department, the Institute for Advanced Study, Princeton, and
the Max Planck Institute f\"ur Mathematik, Bonn for their generous support
and hospitality during a series of visits while this article and its
companion \cite{FL2b} were being prepared. Finally, we thank the anonymous
referee and Simon Donaldson for their editorial suggestions and comments on
the previous versions of this article.


\newpage
\twocolumn
\centerline{{\sc Index of Notation}}
\bigskip

$(A,\Phi)$, $[A,\Phi]$ \hfill \S \ref{subsubsec:PU(2)PairConfigSpaceAndEqns}

$[A,\Phi,\bx]$ \hfill \S \ref{subsec:CompactnessTransversality}

$A^{\ad}$,  $\hat A$ \hfill \S \ref{subsubsection:SpinuConnections}

$A^{\det}$ \hfill Equation \eqref{eq:FixedDetConnection}

$(B,\Psi)$, $[B,\Psi]$  \hfill Equation \eqref{eq:SWPreconfigSpace}

$A_L$ \hfill Lemma \ref{lem:ReducibleSpinu}

$A_\La$ \hfill Equation \eqref{eq:FixedDetConnection}

$F_A$, $F_A^+$ \hfill Equation \eqref{eq:PT}

$F_B$, $F_B^+$, $\Tr(F_B)$ 
\hfill Equations \eqref{eq:CurvDiracSpincConnAction} \&
\eqref{eq:U1Monopole}

$\sA_{\fs}$ \hfill \S \ref{subsubsec:SpincStructures}

$\sA_{\ft}$ \hfill Lemma \ref{lem:AdjointConnection}

$\sA_\kappa^w$, $\sB_\kappa^w$ \hfill \S \ref{subsubsec:SpaceSO(3)Connections}

$\sB_{\ft}$ \hfill \S \ref{subsec:ASDLink}

$\sC_{\fs}$, $\sC^0_{\fs}$, $\tsC_{\fs}$ \hfill Equations
\eqref{eq:SWConfigSpace} \& \eqref{eq:SWPreconfigSpace} 

$\hatsC_{\fs}$, $\hatsC_{\fs}^0$ \hfill Equation \eqref{eq:DefineS1SemiSlice}

$\sC^{p,0}_{\fs}$ \hfill Lemma \ref{lem:FramedSWConfigurationSpace}

$\sC_{\ft}$, $\tsC_{\ft}$ \hfill Equations
\eqref{eq:DefineSpinuConfigurationSpace} \& \eqref{eq:PreConfiguration}
 
$\sC^0_{\ft}$, $\sC^*_{\ft}$, $\sC^{*,0}_{\ft}$ \hfill Equation
\eqref{eq:DefineSmoothConfiguration} 

$\Cl(T^*X)$, $\CCl(T^*X)$ \hfill \S \ref{subsubsec:HermitianCliffordModules}

$D_A$, $D_B$ \hfill \S \ref{subsubsec:SpincStructures}

$\sD_{A,\Phi}$ \hfill Equation \eqref{eq:ConstDetDefOperator}

$\sD_{A,\Phi}^n$, $\sD_{A,\Phi}^t$
\hfill Equation \eqref{eq:DefOperatorReducSplit}

$\sD_{B,\Psi}$ \hfill Equation \eqref{eq:SWDeformationOperator}

$\bsD$  \hfill \S \ref{subsec:LinkOfReduciblesGlobal}

$\bsD^n$ \hfill Equation \eqref{eq:NormalBundle}

$\bD'$, $\bD''$ \hfill Equations \eqref{eq:DiracEndWL2L1Family} \&
\eqref{eq:LambdaWSplitUnivBundleOps}

$\bE$, $\bF$ \hfill Equation \eqref{eq:UniversalBundle}

$\bE'$, $\bF'$, $\bE''$, $\bF''$
\hfill Equation \eqref{eq:LambdaWSplittingUnivBundle}

$E_\ell$ \hfill Equation \eqref{eq:LowerLevelE}

$F_j$, $F_j^t$, $F_j^n$ \hfill Equations \eqref{eq:DefSequenceBundles} \&
\eqref{eq:HilbertSpaceRedCplxSplitting} 

$\sG_\fs$ \hfill \S \ref{subsec:PertSW}

$\hatsG_{\fs}$ \hfill Equation \eqref{eq:DefineHarmonicGT}

$\sG_{\fs}^p$ \hfill Equation \eqref{eq:DefineBasedSWGaugeTransformations}

$\sG_{\ft}$ \hfill Definition \ref{defn:SpinuGaugeTransformation}

$\sG_\kappa^w$ \hfill Lemma \ref{lem:DecompSpinuGaugeTransformation}

$H_{A,\Phi}^i$, $H_{B,\Psi}^i$  \hfill Lemma
\ref{lem:IdentityOfReducibleCohomology} 

$\Jac(X)$ \hfill Equation \eqref{eq:DefineJacobian}

$\bL^{w}_{\ft,\kappa}$ \hfill Definition \ref{defn:ASDLink}

$\bL_{\ft,\fs}$ \hfill Equation \eqref{eq:DefineReducibleLink}

$\LL_{\fs}$ \hfill Equation \eqref{eq:SWUnivLineBundle}

$M_{\fs}$, $M^{0}_{\fs}$ \hfill Equation \eqref{eq:SWModuli}

$\sM_{\ft}$ \hfill Equation \eqref{eq:DefinePU2MonopoleModuliSpace}

$\sM^*_{\ft}$, $\sM^0_{\ft}$,
$\sM^{*,0}_{\ft}$ \hfill Equation \eqref{eq:DefineSmoothPU2Monopole}

$\bar\sM_{\ft}$ \hfill Equation \eqref{eq:UhlCompactPUMonModSpace}

$M_\kappa^w$, $M^{w,*}_{\ka}$ \hfill Equation \eqref{eq:ASDModuliDefn}

$\sM_{\ft}(\Xi,\fs)$ \hfill Definition \ref{defn:ThickenedModuliSpace}

\newpage
\phantom{{\sc Index of Notation}}
\medskip

$N_{\ft}(\Xi,\fs)$
\hfill Equation \eqref{eq:NormalBundle}

$N_{\ft}^\eps(\Xi,\fs)$, $\PP N_{\ft}(\Xi,\fs)$
\hfill Definition \ref{defn:LinkOfReducible}

$\tN_{\ft}(\Xi,\fs)$
\hfill \S \ref{subsubsec:GroupActionsNormalBundleLifts}

$\fS$  \hfill Equations \eqref{eq:PT} \& \eqref{eq:U1Monopole}

$\Spinc(X)$ \hfill \S \ref{subsec:Statement}

$V$, $V^\pm$ \hfill \hfill \S \ref{subsubsec:HermitianCliffordModules}

$\fV$, $\fV^t$, $\fV^n$ 
\hfill Equations \eqref{eq:S1EquivInfRankObstBundle} \& \eqref{eq:DefineV}

$W$, $W^\pm$ \hfill Definition \ref{defn:SpincStructure}

$c(X)$ \hfill \S \ref{subsec:Statement}

$c_1(\fs)$ \hfill  Equation \eqref{eq:C1Spinc}

$c_1(\ft)$, $p_1(\ft)$, $w_2(\ft)$ 
\hfill Equation \eqref{eq:SpinuCharacteristics}

$d_\RR$, $d_\CC$ \hfill Lemma \ref{lem:ReducibleSpinu}

$d_a(\ft)$ \hfill Equation \eqref{eq:Defndana}

$d_s(\fs)$ \hfill Equation \eqref{eq:DimSWModuliSpace}

$d_{\hat A}^i$ \hfill Equation \eqref{eq:ASDDeformationComplex}

$d^{i}_{A,\Phi}$ \hfill Equations \eqref{eq:d1PT} \& \eqref{eq:d0PT}

$d^{i,t}_{A,\Phi}$
\hfill Equations \eqref{eq:ReducibleTangentDerivative} \&
\eqref{eq:ReducibleTangentGauge}

$d^{i,n}_{A,\Phi}$ \hfill Equations \eqref{eq:ReducibleNormalDerivative} \&
\eqref{eq:ReducibleNormalGauge}

$d^{i}_{B,\Psi}$ \hfill Equations \eqref{eq:d0SW} \& \eqref{eq:d1SW}

$\det^{\frac{1}{2}}(V^+)$
\hfill Equation \eqref{eq:SquareRootDeterminantLineBundle}

$\fg_{\ft}$ \hfill Equation \eqref{eq:CanonicalSO(3)Bundle}

$\bh_p$ \hfill \S \ref{subsubsec:AbelConfig}

$n_a(\ft)$ \hfill Equation \eqref{eq:Defndana}

$n_s'(\ft,\fs)$, $n_s''(\ft,\fs)$
\hfill Equation \eqref{eq:DefOfNormalIndices}

$\br$ \hfill Lemma \ref{lem:Retraction}

$r_\Xi$ \hfill Theorem \ref{thm:ChernCharacterOfNormal}

$\fs$ \hfill Definition \ref{defn:SpincStructure}

$\fs\otimes L$ \hfill Equation \eqref{eq:TensorSpincStructure}

$\ft$ \hfill Definition \ref{defn:DefineSpinu}

$\ft_\ell$ \hfill Equation \eqref{eq:DefnEellVell}

$w$ \hfill \S \ref{subsubsec:SpaceSO(3)Connections}

$\bDelta$ \hfill Equation \eqref{eq:JacobianBundle}

$\Xi$, $\Pi_{\Xi}$ \hfill Lemma \ref{thm:DefnOfStabilizeBundle}

$\bgamma$, $\tilde\bgamma$  \hfill Equations
\eqref{eq:NormalBundleNbhdEmbedding} \& \eqref{eq:DefnLiftedNormalEmbedding}

$\gamma_i$, $\gamma_i^*$, $\gamma^{J,*}_i$
\hfill Definition \ref{defn:RelatedBasis}

$\bdelta_L$ \hfill Equation \eqref{eq:LambdaWSplitUnivBundleOps}

$\eta$ \hfill Equation \eqref{eq:U1Monopole}

$\vartheta$ \hfill Equation \eqref{eq:PT}

$\iota$ \hfill Equations \eqref{eq:EmbedSO(3)ConnsIntoPairs} \&
\eqref{eq:DefnOfIota}

$\mu_{\fs}$ \hfill Equation \eqref{eq:SWMuMap}

$\rho$ \hfill Equation \eqref{eq:CliffordMap}

$\varrho$ \hfill Equation \eqref{defn:GInclusion}

$\varrho_L$ \hfill Proof of Lemma \ref{lem:TopU1Embedding}

$\sigma$, $\chi$ \hfill \S \ref{subsec:Statement}

$\tau$ \hfill Equation \eqref{eq:PT}

$\bvarphi$ \hfill Equation \eqref{eq:ObstructionSection}

\onecolumn


\section{Preliminaries}
\label{sec:Prelim}
In this section we recall the framework for gauge theory for $\PU(2)$
monopoles established in \cite{FeehanGenericMetric}, \cite{FL1},
\cite{FLGeorgia}, though we
introduce some notational and other simplifications here.  In \S
\ref{subsec:NonAbelianMonopoles} we
describe the $\PU(2)$ monopole equations while in \S
\ref{subsec:CompactnessTransversality} we recall our Uhlenbeck compactness
and transversality results from
\cite{FeehanGenericMetric}, \cite{FL1}. We review the usual construction of the
moduli space of Seiberg-Witten monopoles in
\S \ref{subsec:PertSW}, although
we use the perturbation parameter $\tau\in\Om^0(\GL(\Lambda^+))$ to achieve
transversality for the moduli space of Seiberg-Witten monopoles rather than
employ the customary perturbation parameter $\eta\in\Om^+(X,i\RR)$.

One slightly awkward issue here then concerns the possible presence of
zero-section Seiberg-Witten monopoles. We recall from \cite[\S
6.3]{MorganSWNotes} that these may be avoided when $b_2^+(X)>0$ by choosing a
generic perturbation $\eta\in\Om^+(X,i\RR)$. For trivial reasons, variation
of the parameter $\tau$ has no effect on the presence or absence of
zero-section solutions and so, even for a generic Riemannian metric on $X$,
they cannot be avoided in all cases. Instead, we circumvent this problem by
a employing a restriction on the second Stiefel-Whitney class of an
$\SO(3)$ bundle comprising part of the definition of the ``\spinu
structure'' introduced below, following \cite{FS1}, which precludes the
existence of flat, reducible connections.  With the aid of the blow-up
trick of \cite{MorganMrowkaPoly}, there is no resulting loss of generality when
computing Donaldson or Seiberg-Witten invariants
\cite{FSBlowUp}, \cite{FSTurkish}.

Finally, in \S \ref{subsec:AbelianCohom} we compute the cohomology ring of
the configuration space of Seiberg-Witten pairs and express the
Seiberg-Witten $\mu$-classes in terms of these generators.

Although we use the generic metrics and Clifford maps of
\cite{FeehanGenericMetric} to achieve transversality for the moduli space
of $\PU(2)$ monopoles, we note that all of the results in this article can be
obtained using the holonomy perturbations of \cite{FL1}.

\subsection{PU(2) monopoles}
\label{subsec:NonAbelianMonopoles}
Throughout this article, $X$ denotes a closed, connected, oriented, smooth
four-manifold. We begin our discussion with a brief review of the
definition of \spinc structures and introduce the concept of a ``\spinu
structure'', which will then allow us to give a convenient definition of
the $\PU(2)$-monopole equations.

\subsubsection{Hermitian Clifford modules, module derivations, and spin
  connections} 
\label{subsubsec:HermitianCliffordModules}
Let $V$ be a complex vector bundle over a Riemannian manifold $(X,g)$. A
real-linear map $\rho:T^*X\to \End(V)$ is a {\em Clifford map}, which
is {\em compatible with $g$}, if it satisfies
\begin{equation}
\rho(\alpha)^2 = -g(\alpha,\alpha)\id_{V},
\qquad \alpha\in C^\8(T^*X).
\label{eq:CliffordMap}
\end{equation}
Equation \eqref{eq:CliffordMap} and the universal property of Clifford
algebras \cite[Proposition I.1.1]{LM} imply that $\rho$ extends to a real
algebra homomorphism $\rho:\Cl(T^*X)\to \End(V)$, where $\Cl(T^*X)$
is the real Clifford algebra, and a complex algebra homomorphism,
$\rho:\CCl(T^*X)\to\End(V)$, where
$\CCl(T^*X)=\Cl(T^*X)\otimes_\RR\CC$ is the complex Clifford algebra.  In
particular, the pair $(\rho,V)$ defines a complex Clifford module.
Recall that $\Cl(T^*X)$ is canonically isomorphic to $\La^\bullet(T^*X)$ as
an orthogonal vector bundle (see \cite[Proposition II.3.5]{LM}
or \cite[Proposition 3.5]{BerlineGetzlerVergne}).  The fibers of
$V$ are irreducible modules for the fibers of $\CCl(T^*X)$,
\cite[Theorem I.5.8]{LM} if and only if $V$ has complex rank $2^n$, when $X$
has dimension $2n$.  If $X$ is even-dimensional, as we shall assume from
now on, the bundle $V$ admits a
splitting $V=V^+\oplus V^-$, where the subbundles $V^\pm$ are the 
$\mp 1$ eigenspaces
of $\rho(\vol)$ and are irreducible $\Cl^0(T^*X)$ modules
\cite[p. 98--99]{LM}, where $\Cl^0(T^*X)\cong \La^{\text{even}}(T^*X)$.
If $V$ is a Hermitian vector bundle, then we require that
$\rho(\alpha)\in\End(V)$ be skew-Hermitian for all $\alpha\in
C^\8(T^*X)$ and call $(\rho,V)$ a {\em Hermitian Clifford module\/}.

A unitary connection $A$ on $V$, where $(\rho,V)$ is a Hermitian Clifford
module, defines a {\em Clifford module derivation\/} $\cov_A$ on $C^\8(V)$ if
\begin{equation}
\label{eq:CovDerivCliffordMult}
\cov^A_\eta(\rho(\alpha)\Phi) 
= 
\rho(\cov_\eta\alpha)\Phi + \rho(\alpha)\cov^A_\eta\Phi,
\end{equation}
where $\alpha\in C^\8(T^*X)$, $\eta\in C^\8(TX)$, and $\Phi \in C^\8(V)$
and $\cov$ is an orthogonal connection on $T^*X$.  (We follow the
convention of \cite[Definition 3.39]{BerlineGetzlerVergne}.)  If $\cov_A$
is a module derivation, it preserves the splitting $V=V^+\oplus V^-$
\cite[Corollary II.4.12]{LM}. We shall also let $\cov$ denote the
canonically induced orthogonal connections on
$\Lambda^\bullet=\Lambda^\bullet(T^*X)$ and the real Clifford algebra
$\Cl(T^*X)$ \cite[Proposition II.4.8]{LM}.

Conversely, if we require that $\cov_A$ be a Clifford module derivation,
then the relation \eqref{eq:CovDerivCliffordMult} defines an orthogonal
connection $\cov$ on $T^*X$
\cite[\S 6.1]{SalamonSWBook}. The connection $A$ is called {\em \spin\/} if
$\cov$ is the Levi-Civita connection on $T^*X$.

We assume for the remainder of the article that $(X,g)$ is an oriented,
Riemannian four-manifold.  The subbundles $V^\pm$ may then be characterized by
requiring that $\rho(\om)V^-=0$ for all $\om\in \Omega^+(X,\RR)$ and
similarly for $V^+$.

\subsubsection{\Spinc\  structures}
\label{subsubsec:SpincStructures} Here we
specialize to the case where the complex Clifford module has rank four.

\begin{defn}
\label{defn:SpincStructure}
We call $\fs=(\rho,W)$ a {\em \spinc structure\/} over an oriented
Riemannian four-manifold $(X,g)$ if $(\rho,W)$ is a Hermitian Clifford
module and $W$ has complex rank four.
\end{defn}

Given a \spinc structure $\fs=(\rho,W)$, one defines
\begin{equation}
\label{eq:C1Spinc}
c_1(\fs) = c_1(W^+).
\end{equation}
If $L$ is a complex line bundle over $X$,  we obtain a new \spinc structure
on $(X,g)$,
\begin{equation}
\label{eq:TensorSpincStructure}
\fs\otimes L= (\rho,W\otimes L),
\end{equation}
with $c_1(\fs\otimes L)= c_1(\fs)+2c_1(L)$. If $\fs$, $\fs'$ are any two
\spinc structures on $(X,g)$ then one has $\fs' = \fs\otimes L$, for a
complex line bundle $L$ on $X$ uniquely determined by $\fs$ and $\fs'$
\cite{SalamonSWBook}.

There is a canonical isomorphism of orthogonal vector bundles,
$$
\fu(W) \cong i\underline{\RR}\oplus\su(W).
$$
where we abbreviate the trivial real line subbundle $i\RR\,\id_{W}\subset
\fu(W)$ by $i\underline{\RR}=X\times i\RR$. The Clifford multiplication
$\rho$ defines canonical isomorphisms $\Lambda^\pm\cong \su(W^\pm)$, where
$\Lambda^\pm=\Lambda^\pm(T^*X)$ are the bundles of self-dual and
anti-self-dual two-forms, with respect to the Riemannian metric $g$ on
$T^*X$.  A unitary connection $\cov_A$ on $W$ determines a unitary
connection on $\det(W)$; conversely, a choice of orthogonal connection on
$\su(W)$ and a unitary connection on $\det(W)$ uniquely determine a unitary
connection on $W$ \cite[\S 2.1.1]{FL1}. In particular, any two unitary
connections on $W$, which are both spinorial with respect to $\cov$,
differ by an element of $\Omega^1(X,i\RR)$, since the induced connection on
$\su(W)\cong\Lambda^2$ is constrained by the choice of $\cov$ on $T^*X$.
The Dirac operator $D=\rho\circ\cov$ on $C^\8(W)$ is not self-adjoint
unless the connection $\cov$ on $T^*X$ is torsion-free---that is,
$\cov$ is the Levi-Civita connection---as one can see from examples.

For $k\ge 2$, we let $\sA_{\fs}$ denote the space of $L^2_k$ spin
connections on $W$. From the preceding remarks (see
\cite[Lemma 6.1]{SalamonSWBook} for details), $\sA_{\fs}$ is an affine space for the
Hilbert space $L^2_k(X,i\Lambda^1)$ of imaginary one-forms:
following the convention of \cite[\S 2(i)]{KMThom}, the action is
given by
\begin{equation}
\label{eq:AffineSpincConnAction}
(B,b)\mapsto B + b\,\id_{W}.
\end{equation}
Thus, denoting the trace on two-by-two complex matrices by $\Tr$,
one has 
\begin{equation}
\label{eq:CurvDiracSpincConnAction}
\Tr(F_{B+b}) = \Tr(F_B) + 2db
\quad\text{and}\quad
D_{B+b} = D_B + \rho(b),
\end{equation}
as we later use when describing the deformation complex for the
Seiberg-Witten equations.

\subsubsection{\Spinu\  structures}
\label{subsubsec:SpinuStructures}
An elegant reformulation of the $\PU(2)$ monopole equations, as discussed
in \cite{PTLocal} and \cite{FL1}, has been described by Mrowka
\cite{MrowkaPrincetonMorseTalk}, based on his joint work with Kronheimer
and, motivated by his comments, we shall give a more invariant definition
of $\PU(2)$ monopoles than the one we presented in \cite{FL1}.

\begin{defn}
\label{defn:DefineSpinu}
We call $\ft=(\rho,V)$ a {\em \spinu structure\/} over an oriented
Riemannian four-manifold $(X,g)$ if $(\rho,V)$ is a Hermitian Clifford
module and $V$ has complex rank eight.
\end{defn}

In the case of even-dimensional manifolds $X$ with spin structures and
complex modules (of arbitrary dimension) for the real Clifford algebra,
$\Cl(T^*X)$, the following result appears as Proposition 3.35 in
\cite{BerlineGetzlerVergne}. 

\begin{lem}
\label{lem:VTensorProdDecomp}
Let $W$ and $V$ be complex Clifford modules over a Riemannian four-manifold
$X$, of rank four and eight respectively. Then there is a rank-two complex
vector bundle $E$ over $X$, unique up to isomorphism, and an isomorphism of
complex Clifford modules,
$$
V \cong W\otimes E.
$$
If $W$ and $V$ are Hermitian Clifford modules, then the bundle $E$ can be
assumed Hermitian and the isomorphism $V \cong W\otimes E$ can be taken
to be an isomorphism of Hermitian Clifford modules.
\end{lem}

\begin{proof}
  Let $\Delta$ be an irreducible $\CCl(\RR^4)$ module, recalling that any
  such module is unique up to equivalence \cite[Theorem I.5.7]{LM}, with
  $\CCl(\RR^4)\cong M_4(\CC)$ as a complex algebra acting on
  $\Delta\cong\CC^4$ by the standard representation. (We use $M_d(\CC)$ to
  denote the algebra of complex $d\times d$ matrices.) Hence, the only
  $\CCl(\RR^4)$-module endomorphisms of $\Delta$ are given by complex
  scalar multiplication and we have an isomorphism of complex vector
  spaces, $\CC\cong\End_{\CCl(\RR^4)}(\Delta)$, given by $z\mapsto
  z\,\id_\Delta$. Therefore, we obtain an isomorphism of complex vector
  spaces
\begin{equation}
  \label{eq:DefineC2fromDeltaDelta}
\CC^2\cong\Hom_{\CCl(\RR^4)}(\Delta,\Delta\oplus\Delta),
\end{equation}
where $(z_1,z_2)\in\CC^2$ is identified with the $\CCl(\RR^4)$-module
homomorphism $\Delta\to \Delta\oplus \Delta$ given by $v\mapsto (z_1v,z_2v)$.
Indeed, the map \eqref{eq:DefineC2fromDeltaDelta} is surjective because we
can compose any homomorphism $\Delta\to \Delta\oplus \Delta$ with
projection onto each factor and then use the fact that
$\End_{\CCl(\RR^4)}(\Delta) \cong \CC$.  Moreover, the map
\eqref{eq:DefineC2fromDeltaDelta} is injective by construction.

Given a complex Clifford module $W$ of rank four, define
\begin{equation}
\label{eq:DefineEfromV}
E=\Hom_{\CCl(T^*X)}(W,V).
\end{equation}
For every $x\in X$ and isomorphism $\CCl(T^*X)|_x\cong
\CCl(\RR^4)$ of complex Clifford algebras, there are
isomorphisms of complex Clifford modules, $W|_x\cong \Delta$ and
$V|_x\cong \Delta\oplus\Delta$.  The isomorphism
\eqref{eq:DefineC2fromDeltaDelta} then implies that $E$ is a
rank-two complex vector bundle over $X$. The map $W\otimes E\to
V$ given by $\Phi\otimes M \to M(\Phi)$,
where $\Phi\in W|_x$ and $M\in E|_x$ for some
$x\in X$, is a $\CCl(\RR^4)$-module
isomorphism since it is fiberwise injective, the ranks of the
two bundles agree, and it is a complex Clifford module
homomorphism by construction. To see that the map $W\otimes E\to
V$ is injective, observe that if $M(\Phi)=0$, then either
$\Phi=0$ or $\Phi\neq 0$ and $M$ has a non-trivial kernel, in
which case one can see from the explicit form of the map
\eqref{eq:DefineC2fromDeltaDelta} that $M|_x$ must then be zero.

Next we show that $E$ is unique up to isomorphism. {}From equation
\eqref{eq:SquareRootDeterminantLineBundle} below, we have 
$$
c_1(E) = \textstyle{\frac{1}{2}}c_1(V^+) - c_1(W^+).
$$
{}From equation \eqref{eq:OrthogDecompsuVForms} below we see that the
bundle $\su(V)$ determines the subbundle $\su(E)$ up to isomorphism,
independently of $W$, and hence $\su(V)$ uniquely determines
$p_1(\su(E))$. Then 
$$
c_2(E) = -\textstyle{\frac{1}{4}}(p_1(\su(E))-c_1(E)^2),
$$
and $E$ is determined up to isomorphism by $c_1(E)$ and $c_2(E)$.
Finally, if $W$ and $V$ are Hermitian then $E$ is Hermitian and the
isomorphism $V\cong W\otimes E$ can be taken to be Hermitian.
\end{proof}

The presentation of $V$ as $W\otimes E$ is not unique.  If $W$ is
replaced by $W\otimes L$, where $L$ is a complex line bundle over $X$,
then $E$ may be replaced by $L^{-1}\otimes E$ and so one also has
$V\cong (W\otimes L)\otimes (L^{-1}\otimes E)$.  By Lemma
\ref{lem:VTensorProdDecomp}, we may always write
$$
(\rho,V) = (\rho,W\otimes E),
$$
for some \spinc structure $\fs=(\rho,W)$ and Hermitian bundle $E$. Although
$W$ and $E$ are not determined by $(\rho,V)$, the \spinu
structure $(\rho,V)$ does determine a complex line bundle,
which is independent of tensor-product decomposition $V\cong W\otimes E$,
\begin{equation}
\label{eq:SquareRootDeterminantLineBundle}
{\det}^{\frac{1}{2}}(V^+) = \det(W^+)\otimes\det(E),
\end{equation}
noting that $\det(V^+)\cong \det(V^-)$, so $\det(V) \cong \det(V^+)^{\otimes
  2}$ since $V=V^+\oplus V^-$.

The \spinu structure $(\rho,V)$ also defines an $\SO(3)$
bundle over $X$. To see this, recall that there are isometries
of orthogonal vector bundles,
\begin{equation}
\label{eq:SkewHermDecomp}
\fu(E)\cong i\underline{\RR}\oplus\su(E),
\quad
M \mapsto (\Tr M)\id_E \oplus \left(M-(\Tr M)\id_E\right),
\end{equation}
where $\underline{\RR}=\RR\,\id_E$ here, and
\begin{equation}
\label{eq:SymmSkewHermComp}
\gl(E) \cong i\fu(E) \oplus \fu(E),
\quad
M\mapsto \thalf(M+M^\dagger) \oplus \thalf(M-M^\dagger),
\end{equation}
with similar isomorphisms for $W$. Since $\gl(V) \cong
\gl(W)\otimes\gl(E)$, the decomposition \eqref{eq:SymmSkewHermComp}
yields an isometry,
\begin{equation}
\label{eq:SkewHermTensorProd}
\fu(V) \cong i\fu(W)\otimes_{\RR}\fu(E).
\end{equation}
Combining the identifications \eqref{eq:SkewHermDecomp} and
\eqref{eq:SkewHermTensorProd},
we obtain isometries of orthogonal vector bundles,
$\fu(V)\cong i\underline{\RR}\oplus\su(V)$ where
$\underline{\RR}=\RR\,\id_{V}$ here, and an orthogonal
decomposition,
\begin{equation}
\label{eq:OrthogDecompsuV}
\su(V) 
\cong 
\su(W)
\oplus
i\su(W)\otimes\su(E)
\oplus
\su(E),
\end{equation}
where we have identified $\RR\,\id_{W}\otimes\su(E)$ with $\su(E)$
and $\su(W)\otimes\RR\,\id_E$ with $\su(W)$.  

Recall that the Clifford map $\rho:T^*X\to\End(W)$ defines
isometries of orthogonal vector bundles $\rho:\Lambda^{\pm}\to\su(W^\pm)$ and
$\rho:\Lambda^2\to\su(W)$ \cite[\S 4.8]{SalamonSWBook}. More generally, one
has:

\begin{lem}
\label{lem:su(W)isomorphism}
  Let $(\rho,W)$ be a \spinc structure over an oriented, Riemannian
  four-manifold $X$. Then the map
  $\rho:\Lambda^\bullet\otimes_\RR\CC\to\End(W)$ defined by composing the
  isomorphism $\Lambda^\bullet\otimes_\RR\CC\cong\CCl(T^*X)$ and the
  Clifford algebra representation $\CCl(T^*X)\to\End(W)$ yields an isometric
  isomorphism of orthogonal vector bundles,
\begin{equation}
\label{eq:su(W)isomorphism}
\Lambda^1\oplus\Lambda^2\oplus i(\Lambda^3\oplus\Lambda^4)
\cong
\su(W).
\end{equation}
\end{lem}

\begin{proof}
The isomorphism $\Lambda^\bullet\cong\Cl(T^*X)$ is given explicitly by
\cite[Proposition 3.5]{BerlineGetzlerVergne}
$$
e^{i_1}\wedge \cdots \wedge e^{i_p} \mapsto \rho(e^{i_1})\cdots\rho(e^{i_p}),
$$
if $\{e^1,\dots,e^4\}$ is a local oriented, orthonormal frame for
$T^*X$. Using this isomorphism and the fact that
$\rho(\alpha)^\dagger=-\rho(\alpha)$ when $\alpha\in\Omega^1(X,\RR)$, it is
easy to see that the left-hand side of \eqref{eq:su(W)isomorphism} is
mapped into $\fu(W)$, the skew-Hermitian endomorphisms of $W$. The elements
of $\rho(\Lambda^1)$ and $\rho(i\Lambda^3)$ give isomorphisms in
$\Hom(W^\pm, W^\mp)$ and their block-matrix representations in
$\End(W^+\oplus W^-)$ have zero components in $\End(W^\pm)$, so they have
zero trace. As remarked earlier, the elements of $\rho(\Lambda^2)$ are
known to have zero trace. Next, $\rho(\vol)|_{W^\pm} = \pm \id_W$ by
definition of $W^\pm$ and thus $\rho(i\Lambda^4)$ has zero trace. Hence the
left-hand side of \eqref{eq:su(W)isomorphism} is also mapped into $\fsl(W)$,
the traceless endomorphisms of $W$.

The map $\CCl(T^*X)\to\End(W)$ is an isomorphism \cite[Proposition
3.19]{BerlineGetzlerVergne} of Clifford modules and in particular is
injective. Since the sum of the ranks of the bundles on the left-hand side
of \eqref{eq:su(W)isomorphism} is $15$, which is equal to the rank of
$\su(W)$, we see that the map \eqref{eq:su(W)isomorphism} must be an
isomorphism, as claimed.
\end{proof}

Note that $\rho:\Lambda^\bullet\otimes_\RR\CC$ assigns $\rho(1) = \id_V$.
For convenience, we define
\begin{equation}
  \label{eq:FormsToSu(W)}
\Lambda^{\sharp}
=  
\Lambda^1\oplus\Lambda^2\oplus i(\Lambda^3\oplus\Lambda^4).
\end{equation}
From the decomposition \eqref{eq:OrthogDecompsuV}, we see that the Clifford
map $\rho:\Lambda^{\sharp}\to\su(V)$ embeds $\Lambda^{\sharp}$
as a subbundle of $\su(V)$ and $\Lambda^\pm\cong\su(W^\pm)$ as subbundles
of $\su(V^\pm)$ via $\omega\mapsto\rho(\omega)$. In particular, the
orthogonal decomposition \eqref{eq:OrthogDecompsuV} is equivalent
to
\begin{equation}
\label{eq:OrthogDecompsuVForms}
\begin{aligned}
\su(V) 
\cong 
\rho(\Lambda^\sharp)
\oplus
i\rho(\Lambda^\sharp)\otimes\su(E)
\oplus
\su(E).
\end{aligned}
\end{equation}
and, upon restriction to $\su(V^\pm)$, 
\begin{equation}
\label{eq:OrthogDecompsuVTwoForms}
\begin{aligned}
\su(V^\pm) 
\cong 
\rho(\Lambda^\pm)
\oplus
i\rho(\Lambda^\pm)\otimes\su(E)
\oplus
\su(E).
\end{aligned}
\end{equation}
Plainly, the decompositions \eqref{eq:OrthogDecompsuVForms} and
\eqref{eq:OrthogDecompsuVTwoForms} are independent of the
tensor-product decomposition $V\cong W\otimes E$, and so the
\spinu structure $\ft$ determines an $\SO(3)$ bundle,
\begin{equation}
\label{eq:CanonicalSO(3)Bundle}
\fg_{\ft} = \su(E).
\end{equation}
Indeed, the subbundle $\fg_{\ft} \subset \su(V)$ can be characterized
invariantly as span of the sections $\xi\in C^\8(\su(V))$ for which
$[\rho(\omega),\xi] = 0$ for all $\omega\in C^\8(\Lambda^{\bullet})$.

The \spinu structure $\ft$ thus defines the characteristic
classes,
\begin{equation}
\label{eq:SpinuCharacteristics}
c_1(\ft) = c_1({\det}^{\frac{1}{2}}(V^+)), \quad
p_1(\ft) = p_1(\fg_{\ft}), \quad\text{and}\quad
w_2(\ft) = w_2(\fg_{\ft}).
\end{equation}
These classes obey the constraints that $w_2(\ft)\in H^2(X;\ZZ/2\ZZ)$ has
an integral lift and
\begin{equation}
\label{eq:SpinuCharConstraints}
\begin{aligned}
c_1(\ft) - w_2(\ft) &\equiv w_2(X)\pmod{2}, 
\\
p_1(\ft) &\equiv w_2(\ft)^2\pmod{4},
\end{aligned}
\end{equation}
recalling that $w_2(\su(E)) \equiv c_1(E)\pmod{2}$ and $c_1(W^+)\equiv
w_2(X)\pmod{2}$. Note that $w_2(\ft)$ is {\em not\/} necessarily equal
to $c_1(\ft)\pmod{2}$.  

Conversely, given a triple $(\Lambda,w,p)$, where $\Lambda,w\in H^2(X;\ZZ)$
obey $\Lambda-w\equiv w_2(X)\pmod{2}$ and $p\in\ZZ$ obeys $p\equiv
w^2\pmod{4}$, there is a \spinu structure $\ft=(\rho,V)$ with
$c_1(\ft)=\Lambda$, $w_2(\ft)\equiv w\pmod{2}$, and $p(\ft)=p$.
Indeed, it suffices to choose a \spinc structure $(\rho,W)$ with
$c_1(W^+)=\Lambda-w$ and a $\U(2)$ bundle $E$ with $c_1(E)=w$ and
$c_2(E)=-\frac{1}{4}(p-w^2)$; one then sets $V=W\otimes E$.

\subsubsection{\Spinu\ structures and spin connections}
\label{subsubsection:SpinuConnections}
Suppose we are given a \spinu structure $(\rho,V)$.
If $A$ is a unitary connection on $V$, it induces an orthogonal connection
on $\fu(V)\cong i\underline{\RR}\oplus\su(V)$ and $\su(V)$, denoted $A^{\ad}$,
by setting
\begin{equation}
\label{eq:CovDerivOnsuV}
\cov^A_\eta(M\Phi) = (\cov^{A^{\ad}}_\eta M)\Phi + M\cov^A_\eta\Phi,
\end{equation}
where $M\in C^\8(\su(V))$, $\eta\in C^\8(TX)$, and $\Phi \in C^\8(V)$; more
succinctly, we have
\begin{equation}
\label{eq:CovDerivOnsuVIsBracket}
\cov^{A^{\ad}}_\eta M = [\cov^A_\eta,M].
\end{equation} 
Suppose $A$ defines a Clifford module derivation $\cov_A$.
The bundle $\su(V)$ is a real Clifford module and since $\cov_A$ is a
Clifford module derivation on $C^\8(V)$, it induces a Clifford module
derivation on $C^\8(\su(V))$, denoted $\cov_{A^{\ad}}$, via equations
\eqref{eq:CovDerivCliffordMult} and \eqref{eq:CovDerivOnsuV}, so that
\begin{equation}
\label{eq:CovDerivCliffordMultsuV}
\cov^{A^{\ad}}_\eta(\rho(\alpha)M) 
= 
\rho(\cov_\eta\alpha)M + \rho(\alpha)\cov^{A^{\ad}}_\eta M,
\end{equation}
where $\alpha\in C^\8(T^*X)$, $\eta\in C^\8(TX)$, and $M \in C^\8(\su(V))$.
Hence, if $A$ is a spin connection on $V$ then $A^{\ad}$ is a spin
connection on $\su(V)$, inducing the Levi-Civita connection $\cov$ on
$T^*X$ and hence on $\Lambda^2$.

\begin{lem}
\label{eq:DecompOfsu(V)Connection}
Let $(\rho,V)$ be a \spinu structure over a Riemannian four-manifold
$(X,g)$ and let $A$ be a spin connection on $V$. Then the induced
orthogonal connection $A^{\ad}$ on $\su(V)$ is spin and preserves the three
subbundles in the orthogonal decomposition \eqref{eq:OrthogDecompsuVForms}
for $\su(V)$, inducing the Levi-Civita connection $\cov$ on the subbundle
$\rho(\Lambda^\sharp)$, an $\SO(3)$ connection $\hat A$ on the subbundle
$\su(E)$, and the tensor-product connection $\cov\otimes\hat A$ on the
subbundle $i\rho(\Lambda^\sharp)\otimes\su(E)$.
\end{lem}

\begin{proof}
The connection $A^{\ad}$ is spin by the remarks preceding the statement of
the lemma. {}From equations \eqref{eq:CovDerivCliffordMultsuV} and
\eqref{eq:CovDerivOnsuVIsBracket} one can see that
$$
\cov^{A^{\ad}}_\eta(\rho(\omega)\id_V)
=
\rho(\cov_\eta\omega)\id_V,
$$
for all $\eta\in C^\8(T^*X)$ and $\omega\in C^\8(\Lambda^\sharp)$,
so $A^{\ad}$ preserves the subspace $\rho(\Lambda^\sharp)\subset\su(V)$,
inducing the Levi-Civita connection $\cov$.

Recall that sections $\xi$ of $\su(E)\subset\su(V)$ can be
characterized as sections of $\su(V)$ having zero commutator with all
$\rho(\omega)$, for $\omega\in C^\8(\Lambda^\sharp)$. Thus, for any such $\xi$,
$\omega$, and $\eta\in C^\8(T^*X)$, we have $[\xi,\rho(\omega)]=0$ and
\begin{align*}
\cov^{A^{\ad}}_\eta[\xi,\rho(\omega)]
&=
[\cov^{A^{\ad}}_\eta\xi,\rho(\omega)] + [\xi,\cov^{A^{\ad}}_\eta\rho(\omega)]
\\
&=
[\cov^{A^{\ad}}_\eta\xi,\rho(\omega)] + [\xi,\rho(\cov_\eta\omega)]
=
[\cov^{A^{\ad}}_\eta\xi,\rho(\omega)].
\end{align*}
Hence, $[\cov^{A^{\ad}}_\eta\xi,\rho(\omega)] = 0$ for all $\omega$ and
therefore $A^{\ad}$ preserves the subspace $\su(E)\subset\su(V)$, on which
it induces an orthogonal connection $\hat A$.

Sections of the subbundle $i\rho(\Lambda^\sharp)\otimes\su(E)$
are linear combinations of sections of the form $i\rho(\omega)\xi$, 
where $\omega\in C^\8(\Lambda^\sharp)$ and $\xi\in C^\8(\su(E))$. Since
$\cov_{A^{\ad}}$ is a Clifford module derivation and induces the
connections $\cov$ and $\hat A$ on the subbundles $\rho(\Lambda^\sharp)$ and
$\su(E)$, respectively, we have
$$
\cov^{A^{\ad}}_\eta(\rho(\omega)\xi)
=
\rho(\cov_\eta\omega)\xi + \rho(\omega)\cov^{\hat A}_\eta\xi,
$$
so the orthogonal connection induced by $A$ on
$i\rho(\Lambda^\sharp)\otimes\su(E)$ is given by $\cov\otimes\hat A$.
\end{proof}

We shall fix, once and for all, a smooth, unitary connection $A_{\La}$ on
the square-root determinant line bundle, ${\det}^{\frac{1}{2}}(V^+)$, and
henceforth require that our unitary connections $A$ on $V=V^+\oplus V^-$
induce the resulting unitary connection on $\det(V^+)$,
\begin{equation}
\label{eq:FixedDetConnection}
A^{\det} = 2A_{\Lambda} \text{ on }\det(V^+),
\end{equation}
where we write $A^{\det}$ for the connection on $\det(V^+)$ induced by
$A|_{V^+}$. If a unitary connection $A$ on $V$ induces a connection $A^{\det} =
2A_\Lambda$ on $\det (V^+)$, then it necessarily induces the connection
$A_\Lambda$ on ${\det}^{\frac{1}{2}}(V^+)$.

\subsubsection{$\PU(2)$ monopoles}
\label{subsubsec:PU(2)PairConfigSpaceAndEqns}
For $k\ge 2$, we let $\sA_{\ft}$ denote the space of $L^2_k$ spin
connections on $V$. From the preceding subsection,
$\sA_{\ft}$ is an affine space for the Hilbert space
$L^2_k(\Lambda^1\otimes\fg_{\ft})$,
\begin{equation}
\label{eq:AffineSpinuConnAction}
(A,a) \mapsto A + a,
\end{equation}
via the inclusion \eqref{eq:OrthogDecompsuV} given by
$\su(E)\subset\su(V)$, $a\mapsto \id\otimes a$.
This descends to an action on the affine space of $\SO(3)$ connections on
$\fg_{\ft}$, 
\begin{equation}
\label{eq:AffineSO(3)ConnAction}
(\hatA,a) \mapsto \hatA + \ad(a),
\end{equation}
with $\ad(a) \in L^2_k(\Lambda^1\otimes\so(\fg_{\ft}))$. We have
\begin{equation}
\label{eq:CurvDiracSpinuConnAction}
\ad^{-1}(F_{\widehat{A+a}}) 
=
\ad^{-1}(F_{\hat A}) + d_{\hat A}a + a\wedge a
\quad\text{and}\quad
D_{A+a} = D_A + \rho(a),
\end{equation}
as we later use when describing the deformation complex for the
$\PU(2)$-monopole equations; note that the map
$\ad:\fg_{\ft}\to\so(\fg_{\ft})$ is an isomorphism.

Recall from the proof of Lemma \ref{lem:VTensorProdDecomp} that the fibers
of $V$ are isomorphic to $\Delta\oplus\Delta$ as complex Clifford algebra
modules, where $\Delta$ is the unique (up to equivalence) irreducible
$\CCl(\RR^4)$ module. {}From equation
\eqref{eq:DefineC2fromDeltaDelta} that there is an isomorphism $\CC^2
\cong\Hom_{\CCl(\RR^4)}(\Delta,\Delta\oplus\Delta)$, identifying
$(z_1,z_2)$ with the homomorphism $\Delta\to\Delta\oplus\Delta$, $v\mapsto
(z_1v,z_2v)$. Hence, there is an isomorphism of complex algebras,
\begin{equation}
\label{eq:MatrixAlgebraIsom}
m:M_2(\CC)\cong\End_{\CCl(\RR^4)}(\Delta\oplus\Delta), 
\quad M\mapsto \id_\Delta\otimes M,
\end{equation}
where $M_2(\CC)$ is the space of complex $2\times 2$ matrices.  If
$u\in\End_{\CCl(\RR^4)}(\Delta\oplus\Delta)$, we define the {\em Clifford
  determinant\/} of $u$ by
\begin{equation}
  \label{eq:CliffordDeterminant}
{\det}_{\CCl}(u) = \det(m^{-1}(u)),
\end{equation}
where $\det(m^{-1}(u))$ is the usual complex determinant of $m^{-1}(u)$ as
a $2\times 2$ complex matrix. Since the Clifford determinant
\eqref{eq:CliffordDeterminant} is invariant under conjugation by
Clifford automorphisms of $\Delta\oplus\Delta$ and
complex automorphisms of $\CC^2$, we can therefore define the Clifford
determinant, $\det_{\CCl}(u)$, of a complex Clifford algebra endomorphism
$u$ of the bundle $V$ by taking the corresponding definition on the fibers
of $V$.

\begin{defn}
\label{defn:SpinuGaugeTransformation}
We say that $u$ is a {\em \spinu automorphism\/} of $V$ if it is an
$L^2_{k+1}$ unitary, complex Clifford algebra automorphism of $V$ with
Clifford-determinant one. We let $\sG_{\ft}$ denote the Hilbert Lie group
of \spinu automorphisms of $V$.
\end{defn}

Let $\sG_\kappa^w$ be the group of $L^2_{k+1}$ unitary, determinant-one
automorphisms of $E$, with Lie algebra $L^2_{k+1}(\su(E))$, where
$w=c_1(E)$ and $\kappa=c_2(E)-\frac{1}{4}c_1(E)^2$. While at first
glance it might seem more natural to relax the requirement that
$u\in\sG_\ft$ have Clifford-determinant one to that of complex-determinant
one, our refinement gives us a useful identification of $\sG_\ft$ with
$\sG_\kappa^w$:

\begin{lem}
\label{lem:DecompSpinuGaugeTransformation}
Suppose $(\rho,V)$ is a \spinu structure with $V=W\otimes E$, for some \spinc
structure $(\rho,W)$ and rank-two Hermitian bundle $E$. Then the following
map is an isomorphism of Hilbert Lie groups:
\begin{equation}
  \label{eq:DecompSpinuGaugeTransformation}
m:\sG_\kappa^w \cong \sG_{\ft}, \quad u\mapsto \id_W\otimes u.  
\end{equation}
\end{lem}

\begin{proof}
  Plainly, the map is an injective homomorphism, so it remains to show that
  it is surjective. Suppose $v\in\sG_\ft$. For any $x\in X$, the
  remarks preceding Definition \ref{defn:SpinuGaugeTransformation} imply
  that we may write $v|_x = \id_{W_x}\otimes u_x$, where $u_x\in
  \End_\CC(E_x)$ and so we have $v = \id_W\otimes u$ for some
  $u\in\End_\CC(E)$. But $\det(u) = \det_{\CCl}(v) = 1$ and as
  $\id_W\otimes u^\dagger u = v^\dagger v= \id_V$, we must have $u^\dagger
  u = \id_E$ and thus $u\in\sG_\kappa^w$, as desired.
\end{proof}

Consequently, the Hilbert Lie group $\sG_{\ft}$ has Lie algebra
$L^2_{k+1}(\fg_{\ft})\subset L^2_{k+1}(\su(V))$. Note that if $u$ is a
\spinu automorphism of $V$, then
$\Ad(u)$ preserves three factors in the orthogonal decomposition
\eqref{eq:OrthogDecompsuVTwoForms} of $\su(V)$; it acts as the identity on
the subbundle $\rho(\Lambda^2)$, as an orthogonal gauge transformation
$\hatu$ on $\su(E)$, and as $\id\otimes \hat u$ on the subbundle
$i\rho(\Lambda^2)\otimes\su(E)$.

For an $L^2_k$ section $\Phi$ of $V^+$, we let $\Phi^*$ denote its
pointwise Hermitian dual and let $(\Phi\otimes\Phi^*)_{00}$ be the
component of $\Phi\otimes\Phi^* \in i\fu(V^+)$
which lies in the factor
$\su(W^+)\otimes\su(E)$ of the decomposition \eqref{eq:OrthogDecompsuV} of
$i\fu(V^+)\cong\underline{\RR}\oplus i\su(V^+)$ (with $V^+$ in place of $V$).
The Clifford multiplication
$\rho$ defines an isomorphism $\rho:\La^+\to\su(W^+)$ and thus an
isomorphism $\rho=\rho\otimes\id_{\su(E)}$ of $\La^+\otimes\su(E)$ with
$\su(W^+)\otimes\su(E)$.

The pre-configuration space of pairs on $V$ is given by
\begin{equation}
\label{eq:PreConfiguration}
\tsC_{\ft} = \sA_{\ft} \times L^2_k(V^+),
\end{equation}
with tangent spaces $L^2_k(\fg_{\ft})\oplus L^2_k(V^+)$. We call a
pair $(A,\Phi) \in \tsC_{\ft}$ a {\em $\PU(2)$ monopole\/} if
\begin{equation}
\label{eq:PT}
\fS(A,\Phi)
=
\begin{pmatrix}
\ad^{-1}(F_{\hat A}^+) - \tau\rho^{-1}(\Phi\otimes\Phi^*)_{00}
\\
D_A\Phi + \rho(\vartheta)\Phi
\end{pmatrix}
=
0,
\end{equation}
where $F_{\hat A}^+ \in L^2_{k-1}(\Lambda^+\otimes\so(\fg_{\ft}))$ is the
self-dual component of the curvature $F_{\hat A}$ of $\hat A$ while
$D_A=\rho\circ\cov_A:L^2_k(V^+)\to L^2_{k-1}(V^-)$
is the Dirac operator, $\tau \in L^2_{k+1}(X,\GL(\Lambda^+))$ and
$\vartheta \in L^2_{k+1}(\Lambda^1\otimes\CC)$ are perturbation
parameters. We let
\begin{equation}
\label{eq:DefinePU2MonopoleModuliSpace}
\sM_{\ft} = \{[A,\Phi] \in \sC_{\ft}:\text{ $(A,\Phi)$ satisfies
\eqref{eq:PT}}\},
\end{equation}
be the moduli space of solutions to \eqref{eq:PT} cut out of the
configuration space,
\begin{equation}
\label{eq:DefineSpinuConfigurationSpace}
\sC_{\ft} = \tsC_{\ft}/\sG_{\ft},
\end{equation}
where $u\in\sG_{\ft}$ acts by $u(A,\Phi) =
(u_*A,u\Phi)$. The linearization of the map $\sG_{\ft}\to \tsC_{\ft}$,
$u\mapsto u(A,\Phi)$, at $\id_V\in\sG_{\ft}$, is given by
\begin{equation}
\label{eq:LinearizationGaugeGroupAction}
\zeta
\mapsto
-d^0_{A,\Phi}\zeta = (-d_{\hat A}\zeta,\zeta\Phi),
\end{equation}
with $L^2$-adjoint $d^{0,*}_{A,\Phi}$, where $\zeta\in L^2_{k+1}(\fg_{\ft})$.

The circle $S^1$ defines a family of unitary gauge transformations on $V$
acting by scalar multiplication, so that
\begin{equation}
\label{eq:S1ZAction}
S^1\times\tsC_{\ft}\to \tsC_{\ft},
\quad
(e^{i\theta}, (A,\Phi)) \mapsto (A,e^{i\theta}\Phi).
\end{equation}
Because scalar multiplication by $S^1$ commutes with $\sG_{\ft}$, the action
\eqref{eq:S1ZAction} descends to an action on $\sC_{\ft}$ and on
$\sM_{\ft}$.

\begin{rmk}
\label{rmk:GaugeGroupConvention}
Note that we break here with our former convention
\cite{FL1}, \cite{FLGeorgia} of considering
$\sM_{\ft}$ and $\sC_{\ft}$ as quotients by
$S^1\times_{\{\pm 1\}}\sG_{\ft}$, rather than $\sG_{\ft}$ as defined here.
\end{rmk}

We call a spin connection $A$ on $V$ {\em reducible\/} if it splits as a
direct sum of connections on $V=W\oplus W'$, where $\fs = (\rho,W)$ and
$\fs' = (\rho,W')$ are \spinc structures, and call $A$ {\em irreducible\/}
otherwise. We write $\ft=\fs\oplus \fs'$ and $A=B\oplus B'$, for the induced
spin connections $B$ on $W$ and $B'$ on $W'$.

We call a connection $\hat A$ on an $\SO(3)$ bundle $F$ {\em reducible\/}
if it splits as a direct sum $d_\RR\oplus A_L$ on $F =
i\underline{\RR}\oplus L$, where $A_L$ is a unitary connection on a complex
line bundle $L$ over $X$ and $d_\RR$ is the product connection on
$i\underline{\RR}$, and call $\hat A$ {\em irreducible\/} otherwise.  The
following lemma summarizes the relationship between these notions of
reducibility:

\begin{lem}
\label{lem:ReducibleSpinu}
Let $(\rho,V)$ be a \spinu structure on $X$.  Then a spin connection $A$
on $V$ is reducible with respect to a splitting $V=W\oplus W'$, where
$(\rho,W)$ and $(\rho,W')$ are \spinc structures, if and only if the induced
$\SO(3)$ connection $\hat A$ on $\fg_{\ft}$ is reducible with respect to a
splitting $\fg_{\ft}=i\underline{\RR}\oplus L$, where $L$ is a complex
line bundle over $X$ such that $W'=W\otimes L$.  If $A=B\oplus B'$
is reducible, then there is a unitary connection $A_L$ on $L$ such that
\begin{equation}
\label{eq:ReducibleSpinuConn}
B' = B\otimes A_L
\quad\text{and}\quad
\hat A= d_{\RR}\oplus  A_L.
\end{equation}
If $A^{\det} = 2A_\Lambda$ on $\det(V^+)$, then 
\begin{equation}
\label{eq:InducedConnOnLn}
A_L = A_\Lambda\otimes (B^{\det})^*.
\end{equation}
\end{lem}

\begin{rmk}
We write $B^{\det}$ for the connection on $\det(W^+)$ induced by
$B|_{W^+}$, where $B$ is a spin connection on $W=W^+\oplus W^-$.
\end{rmk}

\begin{proof}
Suppose that $A=B\oplus B'$ is a reducible spin connection on $V$ with
respect to the splitting $V=W\oplus W'$. There is a unique complex line
bundle $L$ over $X$ such that $W'=W\otimes L$. If $A_L$ is a
unitary connection on $L$, then both $B'$ and $B\otimes A_L$ are
spin connections on $W'$. But any two spin connections on $W'$ differ by an
element of $\Omega^1(X,i\RR)$, via the action
\eqref{eq:AffineSpincConnAction}, and thus we can write $B'=B\otimes
A_L$ for some unitary connection $A_L$ on $L$. 

The unitary connection $\tilde A = d_\CC\oplus A_L$ on
$E=\underline{\CC}\oplus L$ induces the $\SO(3)$ connection $d_\RR\oplus
A_L$ on $\su(E)\cong i\underline{\RR}\oplus L$.  To see this, we pass to
a local trivialization of $E$ and view $\tilde A$ as a connection matrix
one-form,
$$
\tilde A = \begin{pmatrix}0 & 0 \\ 0 & A_L\end{pmatrix}
\in\Omega^1(\fu(E)).
$$
The matrix $\tilde A$ acts on a section $\xi$ of $\su(E)$ via the
adjoint representation,
$$
\hat A\cdot\xi
= 
\ad(\tilde A)\xi
=
[\tilde A,\xi]
=
\begin{pmatrix}0 & -\bar A_L\bar z \\ A_Lz & 0\end{pmatrix},
\text{ where }
\xi = \begin{pmatrix}\nu & -\bar z \\ z & -\nu\end{pmatrix}
\in C^\8(\su(E)).
$$
Thus, $\hat A\cdot(\nu,z) = (0,A_Lz)$, with respect to the isomorphism
$\su(E)\cong i\underline{\RR}\oplus L$ of Lemma \ref{lem:su(E)Split} and so,
as a connection, we see that $\hat A = d_\RR \oplus A_L$.

Therefore, the spin connection $A=B\otimes (d_\CC\oplus A_L)$ on
$W\otimes E$ induces the $\SO(3)$ connection $d_\RR\oplus A_L$ on
$\su(E)\cong i\underline{\RR}\oplus L$ via the decomposition
\eqref{eq:OrthogDecompsuV},
$$
\su(V) \cong \su(E)\oplus \su(W)\oplus i\su(W)\otimes\su(E). 
$$
The connections $B$ on $W$ and $d_\CC\oplus A_L$ on $E$ induce the
connection $B^{\det}\otimes A_L$ on $\det(W^+)\otimes\det(E) =
{\det}^{\frac{1}{2}}(V^+)$. By convention, our spin connections $A$ on $V$
induce the connection $2A_\Lambda$ on $\det(V^+)$, and thus $A_\Lambda =
B^{\det}\otimes A_L$.

Conversely, if $\hat A=d_\RR\oplus A_L$ is an $\SO(3)$ connection which is
reducible with respect to the splitting $\su(E) = i\underline{\RR}\oplus
L$, then $\hat A$ lifts to a $\U(2)$ connection $d_\CC\oplus A_L$ on $E$
which is reducible with respect to the splitting $E=\underline{\CC}\oplus
L$, and lifts to a spin connection $A=B\otimes(d_\CC\oplus
A_L)=B\oplus B'$ on $V$ which is reducible with respect to the splitting
$W\otimes(\underline{\CC}\oplus L) = W\oplus W\otimes L$.
\end{proof}

We say that a pair $(A,\Phi)\in\tsC_{\ft}$ is irreducible
(respectively, reducible) if the connection $A$ is irreducible
(respectively, reducible). We let $\sC^*_{\ft}\subset \sC_{\ft}$ be the open
subspace of gauge-equivalence classes of irreducible pairs. If $\Phi\equiv
0$ on $X$, we call $(A,\Phi)$ a {\em zero-section\/} pair.  We let
$\sC^0_{\ft}\subset \sC_{\ft}$ be the open subspace of gauge-equivalence
classes of non-zero-section pairs and recall that
\begin{equation}
\label{eq:DefineSmoothConfiguration}
\sC^{*,0}_{\ft} =
\sC^*_{\ft}\cap\sC^0_{\ft}
\end{equation}
is a Hausdorff, Hilbert manifold
\cite[Proposition 2.8]{FL1} represented by pairs with trivial stabilizer
in $\sG_{\ft}$. Let
\begin{equation}
\label{eq:DefineSmoothPU2Monopole}
\sM_{\ft}^{*,0} = \sM_{\ft}\cap\sC^{*,0}_{\ft}
\end{equation}
be the open subspace of the moduli space $\sM_{\ft}$ represented by
irreducible, non-zero-section $\PU(2)$ monopoles; the subspaces
$\sM_{\ft}^*$ and $\sM_{\ft}^0$ are defined analogously.

\subsubsection{Spaces of $\SO(3)$ connections}
\label{subsubsec:SpaceSO(3)Connections}
Given a class in $H^2(X;\ZZ/2\ZZ)$ with an integral lift $w\in H^2(X;\ZZ)$
and an integer $p$ obeying $p^2\equiv w^2\pmod{4}$, set $\kappa=-\frac{1}{4}p$,
and choose any Hermitian rank-two bundle $E$ with $c_1(E)=w$, so 
$w_2(\su(E))=w\pmod{2}$, and $p_1(\su(E))=-4\kappa$. Let $\sA_\kappa^w(X)$
denote the affine space of $L^2_k$ orthogonal connections on $\su(E)$ and
let $\sB_\kappa^w(X) = \sA_\kappa^w(X)/\sG_\kappa^w$ be the quotient, where
$\sG_\kappa^w$ acts on $\su(E)$ via the adjoint action.

Let $\sA_{\kappa}^{w,*}(X)$ and $\sB^{w,*}_{\ka}(X)$ be the subspace of
irreducible $\SO(3)$ connections and its quotient.  The moduli space of
anti-self-dual connections on $\su(E)$ is then defined by
\begin{equation}
\label{eq:ASDModuliDefn}
M_{\kappa}^{w}(X)  = \{[\hat A]\in\sB_{\kappa}^{w}(X): F_{\hat A}^+ = 0\},
\end{equation}
with $M_\kappa^{w,*} = M_\kappa^w\cap\sB_\kappa^{w,*}$. We
follow the convention of \cite{KMStructure} in taking the quotient by
the group of determinant-one, unitary automorphisms of $E$ rather than the
group of determinant-one, orthogonal automorphisms of $\su(E)$.

Lemma \ref{eq:DecompOfsu(V)Connection} implies that a spin connection $A$
on $V=W\otimes E$ determines a unique orthogonal connection $\hat A$ on
$\su(E)$ and, conversely, that every orthogonal connection on $\su(E)$ is
in the image of the map $A\mapsto \hat A$.  Hence, one can easily translate
between the conventions employed in the present article and those of
\cite{FL1}, \cite{FeehanGenericMetric}, \cite{FLGeorgia}:

\begin{lem}
\label{lem:AdjointConnection}
Let $A_\La$ be the fixed unitary connection on the complex line bundle
${\det}^{\frac{1}{2}}(V^+)$ and let $\sA_{\ft}$ denote the corresponding
space of spin connections on $V$, as described in \S
\ref{subsubsection:SpinuConnections}. For each $\SO(3)$ connection $\hat A$
on $\su(E)$, let $A$ denote the corresponding spin connection on $V$. With
respect to the action of $\sG_\kappa^w$ on $\sA_{\kappa}^w$ and of
$\sG_{\ft}$ on $\sA_{\ft}$, and the identification
$\sG_\kappa^w\cong\sG_{\ft}$, the following map is a gauge-equivariant
bijection:
\begin{equation}
\label{eq:EConnToSpinuConn}
\sA_{\kappa}^w\cong \sA_{\ft}, \quad \hat A\mapsto A.
\end{equation}
\end{lem}

\subsection{Uhlenbeck compactness and transversality for PU(2) monopoles}
\label{subsec:CompactnessTransversality}
We briefly recall our Uhlenbeck compactness and transversality results
\cite{FeehanGenericMetric}, \cite{FL1} for the moduli space of $\PU(2)$
monopoles with the perturbations discussed in the preceding section. The
moduli space of $\PU(2)$ monopoles is non-compact but has an Uhlenbeck
compactification analogous to the compactification $\bar M_\kappa^w$ of the
moduli space of anti-self-dual connections on an $\SO(3)$ bundle \cite[\S
4.4]{DK}.

Given a non-negative integer $\ell$ and a \spinu structure $\ft=(\rho,V)$,
Lemma \ref{lem:VTensorProdDecomp} allows us to write $V=W\otimes E$ for
some choice of \spinc structure $(\rho,W)$ and corresponding Hermitian,
rank-two bundle $E$. Let $E_\ell\to X$ be the Hermitian,
rank-two bundle with 
\begin{equation}
\label{eq:LowerLevelE}
c_1(E_\ell)=c_1(E) \quad\text{and}\quad c_2(E_\ell)=c_2(E)-\ell.
  \end{equation}
We then define a \spinu structure $\ft_\ell=(\rho,V_\ell)$ on $X$ by setting 
\begin{equation}
  \label{eq:DefnEellVell}
V_\ell = W\otimes E_\ell,
  \end{equation}
and observe that
\begin{equation}
\label{eq:LowerLevelSpinu}
c_1(\ft_\ell) = c_1(\ft),
\quad
p_1(\ft_\ell) = p_1(\ft)+4\ell,
\quad\text{and}\quad
w_2(\ft_\ell) = w_2(\ft). 
\end{equation}
We say that a sequence of points $[A_\alpha,\Phi_\alpha]$ in $\sC_{\ft}$
{\em converges\/} to a point $[A,\Phi,\bx]$ in
$\sC_{\ft_\ell}\times\Sym^\ell(X)$ if the following hold:
\begin{itemize}
\item There is a sequence of $L^2_{k+1,\loc}$ \spinu
bundle isomorphisms $u_\alpha:V|_{X\less\bx}\to
V_{\ell}|_{X\less\bx}$ such that the sequence of monopoles
$u_\alpha(A_\alpha,\Phi_\alpha)$ converges to
$(A,\Phi)$ in $L^2_{k,\loc}$ over $X\less\bx$, and
\item The sequence of measures
$|F_{A_\alpha}|^2$ converges
in the weak-* topology on measures to $|F_A|^2 +
8\pi^2\sum_{x\in\bx}\delta(x)$.
\end{itemize}
There is a natural extension of this definition of convergence of points in
$\sC_{\ft}$ to one for sequences in the space of ``ideal pairs'',
$\sqcup_{\ell=0}^{\8}(\sC_{\ft_\ell}\times\Sym^\ell(X))$, and this serves
to define the ``Uhlenbeck topology'' on the space of ideal pairs. We define
the topological space
\begin{equation}
  \label{eq:UhlCompactPUMonModSpace}
\bar\sM_{\ft}
=
\text{Closure}(\sM_{\ft})\subset
\bigsqcup_{\ell=0}^{\8}(\sM_{\ft_\ell}\times\Sym^\ell(X)),
\end{equation}
where the closure is taken with respect to the Uhlenbeck topology on the
(second countable, Hausdorff) space of {\em ideal $\PU(2)$ monopoles},
$\sqcup_{\ell=0}^{\8}(\sM_{\ft_\ell}\times\Sym^\ell(X))$.  We call the
intersection of $\bar\sM_{\ft}$ with $\sM_{\ft_\ell}\times\Sym^\ell(X)$ a {\em
  lower-level\/} of the compactification $\bar\sM_{\ft}$ if $\ell>0$ and call
$\sM_{\ft}$ the {\em top\/} or {\em highest level\/}.

\begin{thm}
\label{thm:Compactness}
\cite{FL1} Let $X$ be a closed, oriented, Riemannian, smooth four-manifold
with \spinu structure $\ft$. Then there is a positive integer $N$,
depending at most on the curvature of the fixed unitary connection on
$\det(V^+)$ together with $p_1(\ft)$, such that the topological space
$\bar\sM_{\ft}$ is second countable, Hausdorff, compact, and given by the
closure of $\sM_{\ft}$ in the space of ideal $\PU(2)$ monopoles,
$\sqcup_{\ell=0}^{N}(\sM_{\ft_\ell}\times\Sym^\ell(X))$, with respect to
the Uhlenbeck topology.
\end{thm}

Theorem \ref{thm:Compactness} is a special case of the more general result
proved in \cite{FL1} for the moduli space of solutions to the $\PU(2)$
monopole equations in the presence of holonomy perturbations.
The existence of an Uhlenbeck compactification for the moduli space of
solutions to the unperturbed $\PU(2)$ monopole equations
\eqref{eq:PT} was announced by Pidstrigatch
\cite{PTLectures} and an argument was outlined in \cite{PTLocal}.
A similar argument for the equations
\eqref{eq:PT} (without perturbations)
was outlined by Okonek and Teleman in
\cite{OTQuaternion}. An independent proof of Uhlenbeck compactness for
\eqref{eq:PT} and other perturbations of these equations is
also given in \cite{TelemanGenericMetric}.

We recall from \cite[Equation (2.37)]{FL1} that the elliptic deformation
complex for the moduli space $\sM_{\ft}$ is given by
\begin{equation}
\begin{CD}
L^2_{k+1}(\fg_{\ft})
@>{d_{A,\Phi}^0}>>
\begin{matrix}
L^2_k(\Lambda^1\otimes\fg_{\ft}) \\
\oplus \\
L^2_k(V^+)
\end{matrix}
@>{d_{A,\Phi}^1}>>
\begin{matrix}
L^2_{k-1}(\Lambda^+\otimes\fg_{\ft}) \\
\oplus \\
L^2_{k-1}(V^-)
\end{matrix}
\end{CD}
\label{eq:ConstDetDefComplex}
\end{equation}
with elliptic deformation operator
\begin{equation}
\sD_{A,\Phi} = d_{A,\Phi}^{0,*} + d^1_{A,\Phi}
\label{eq:ConstDetDefOperator}
\end{equation}
and cohomology $H_{A,\Phi}^\bullet$. Here, $d^1_{A,\Phi}$ is the
linearization at the pair $(A,\Phi)$ of the gauge-equivariant map $\fS$
defined by the equations \eqref{eq:PT}, so
\begin{equation}
  \label{eq:d1PT}
d^1_{A,\Phi}(a,\phi)
=
(D\fS)_{A,\Phi}(a,\phi)
=
\begin{pmatrix}
d_{\hat A}^+a-\tau\rho^{-1}(\phi\otimes\Phi^*+\Phi\otimes\phi^*)_{00} \\
(D_A+\rho(\vartheta))\phi+\rho(a)\Phi
\end{pmatrix},
\end{equation}
while $-d^0_{A,\Phi}$ is the differential of the map
$\sG_{\ft}\to\tsC_{\ft}$, $u\mapsto u(A,\Phi)=(A-(d_Au)u^{-1},u\Phi)$, so
\begin{equation}
  \label{eq:d0PT}
d^0_{A,\Phi}\zeta
=
(d_A\zeta,\zeta\Phi).
\end{equation}
The space $H_{A,\Phi}^0 = \Ker d_{A,\Phi}^0$ is the Lie algebra of the
stabilizer in $\sG_{\ft}$ of a pair $(A,\Phi)$ and $H^1_{A,\Phi}$ is the
Zariski tangent space to $\sM_{\ft}$ at a point $[A,\Phi]$.  If
$H_{A,\Phi}^2=0$, then $[A,\Phi]$ is a regular point of the zero locus of
the $\PU(2)$ monopole equations \eqref{eq:PT} on $\sC_{\ft}$.

We now turn to the question of transversality. Let $D_{A,\cov}$ and $D_A$
be the Dirac operators on $V^+$ defined, respectively, by unitary
connections $A$ on $V$ which are spin with respect to an $\SO(4)$
connection $\cov$ and spin, in the usual sense, with respect to the
Levi-Civita connection on $T^*X$: the two operators differ by an element
$\rho(\vartheta') \in \Hom(V^+,V^-)$
\cite[Lemma 3.1]{FeehanGenericMetric}, where $\vartheta' \in
\Omega^1(X,\CC)$. Even though a unitary connection on $V$ will not
necessarily induce a torsion-free connection on $T^*X$ for generic
pairs $(g,\rho)$ of Riemannian metrics and compatible Clifford maps, we can
assume that the Dirac operator $D_A$ in
\eqref{eq:PT} is defined using the Levi-Civita connection for the metric
$g$ by absorbing the difference term $\rho(\vartheta')$ into the
perturbation term $\rho(\vartheta)$. Given any fixed pair $(g_0,\rho_0)$
satisfying \eqref{eq:CliffordMap} and automorphism $f\in
C^\8(\GL(T^*X))$, then $(g,\rho) = (f^*g_0,f^*\rho_0)$ is again a compatible
pair; the pair $(g,\rho)$ is {\em generic\/} if $f$ is generic.

\begin{thm}
\label{thm:Transversality}
\cite{FeehanGenericMetric}
Let $X$ be a closed, oriented, smooth four-manifold with \spinu structure
$\ft$. Then for a generic, $C^\8$ pair $(g,\rho)$ satisfying
\eqref{eq:CliffordMap} and generic, $C^\8$ parameters
$(\tau,\vartheta)$, the moduli space
$\sM_{\ft}^{*,0}=\sM_{\ft}^{*,0}(g,\rho,\tau,\vartheta)$ of $\PU(2)$ monopoles is
a smooth manifold of the expected dimension,
$$
\dim \sM_{\ft}^{*,0}
=
d_a + 2n_a,
$$
where
\begin{equation}
\label{eq:Defndana}
\begin{aligned}
d_a(\ft)
&=
-2p_1(\ft)-\tthreehalf(\chi+\sigma),
\\
n_a(\ft)
&=
\tquarter( p_1(\ft)+c_1(\ft)^2-\sigma).
\end{aligned}
\end{equation}
\end{thm}

In equation \eqref{eq:Defndana}, the quantity $d_a(\ft)$ is the
expected dimension of the moduli space of anti-self-dual connections on
$\fg_{\ft}$, while $n_a(\ft)$ is the complex index of the Dirac
operator on $C^\8(V^+)$.
Theorem \ref{thm:Transversality} was proved independently, using a somewhat
different method, by A. Teleman \cite{TelemanGenericMetric}.

\subsection{Moduli spaces of Seiberg-Witten monopoles}
\label{subsec:PertSW}
We recall the definition of the space of Seiberg-Witten monopoles.
Let $\fs=(\rho,W)$ be a \spinc structure on $X$ and let
\begin{equation}
\label{eq:SWPreconfigSpace}
\tsC_{\fs} = \sA_{\fs} \times L^2_k(W^+)
\end{equation}
be the pre-configuration space of pairs $(B,\Psi)$, with tangent space
$L^2_k(X,i\RR)\oplus L^2_k(W^+)$ due to the affine structure
\eqref{eq:AffineSpincConnAction}. We call a unitary automorphism of $W$
a {\em \spinc automorphism\/} if it is a complex Clifford module
endomorphism of $W$. Thus, \spinc automorphisms induce the identity on the
factor $\su(W)\cong\Lambda^2(T^*X)$ of the orthogonal decomposition
$\fu(W)=i\underline{\RR}\oplus\su(W)$. We let $\sG_{\fs}$ denote the group
of $L^2_{k+1}$ \spinc automorphisms of $W$ and observe that there is an
isomorphism \cite[\S\S4.3--4.5]{MorganSWNotes}
$$
L^2_{k+1}(X,S^1)\to \sG_{\fs}, \quad s\mapsto s\cdot\id_{W}.
$$
Hence, the group $\sG_{\fs}$ has Lie algebra $L^2_{k+1}(X,i\RR)$ via the
identification 
$$
L^2_{k+1}(X,i\RR) \cong T_{\id}\sG_{\fs},
\quad
f\mapsto f\cdot\id_W,
$$
and acts smoothly by pushforward on $\tsC_{\fs}$,
\begin{equation}
\label{eq:SWGaugeAction}
(s,(B,\Psi))
\mapsto
(s_*B,s\Psi) 
= 
(B-s^{-1}ds,s\Psi).
\end{equation}
We thus obtain a configuration space 
\begin{equation}
\label{eq:SWConfigSpace}
\sC_{\fs} = \tsC_{\fs}/\sG_{\fs}
\end{equation}
and let $\tsC_{\fs}^0\subset\tsC_{\fs}$ denote the open subspace of pairs
$(B,\Psi)$ with $\Psi\not\equiv 0$, where $\sG_{\fs}$ acts freely on
$\tsC_{\fs}^0$ with quotient $\sC_{\fs}^0$, a smooth Hilbert manifold.

We call a pair $(B,\Psi) \in \tsC_{\fs}$ a {\em Seiberg-Witten monopole\/} if
\begin{equation}
\label{eq:U1Monopole}
\fS(B,\Psi)
=
\begin{pmatrix}
\Tr(F^+_B) - \tau\rho^{-1}(\Psi\otimes\Psi^*)_0-\eta \\
D_B\Psi + \rho(\vartheta)\Psi,
\end{pmatrix}
= 0,
\end{equation}
where $F_B^+ \in L^2_{k-1}(\Lambda^+\otimes\fu(W^+))$ is the self-dual
component of the curvature $F_B$ of $B$ and $\Tr(F_B^+) \in
L^2_{k-1}(\Lambda^+\otimes\su(W^+))$ is the trace-free part,
$D_B = \rho\circ\cov_B:C^\8(W^+)\to C^\8(W^-)$ is the Dirac operator
defined by the spin connection $B$, the perturbation terms $\tau$ and
$\vartheta$ are as defined in our version of the $\PU(2)$ monopole
equations \eqref{eq:PT}, and where $\eta\in C^\8(i\La^+)$ is an additional
perturbation (see Remark \ref{rmk:NonGenericSWTwoFormPert}). The quadratic
term $\Psi\otimes\Psi^*$ lies in $C^\8(i\fu(W^+))$ and
$(\Psi\otimes\Psi^*)_0$ denotes the traceless component lying in
$C^\8(i\su(W^+))$, so $\rho^{-1}(\Psi\otimes\Psi^*)_0 \in C^\8(i\La^+)$.

\begin{rmk}
\label{rmk:NonGenericSWTwoFormPert}
We note that in the usual presentation of the Seiberg-Witten equations
\cite{KMThom}, \cite{MrowkaOzsvathYu}, one takes $\tau=\id_{\Lambda^+}$ and
$\vartheta=0$, while $\eta$ is a generic perturbation. However, we shall
see in Lemma \ref{lem:RestrictionOfPU2MonopoleEquation} that in order to
identify solutions to the Seiberg-Witten equations \eqref{eq:U1Monopole}
with reducible solutions to the $\PU(2)$ monopole equations \eqref{eq:PT},
we need to employ the perturbations given in equation
\eqref{eq:U1Monopole} and choose
\begin{equation}
  \label{eq:EtaSWchoice}
  \eta = F_{A_\Lambda}^+,
\end{equation}
where $A_\Lambda$ is the fixed unitary connection on the line bundle
$\det^{\frac{1}{2}}(V^+)$ with Chern class denoted by $c_1(\ft)=\Lambda\in
H^2(X;\ZZ)$ and represented by the real two-form $(1/2\pi
i)F_{A_\Lambda}$. In particular, this choice of $\eta$ is {\em not
  generic\/} in the sense of \cite[Proposition 6.3.1 \& Corollary
6.3.2]{MorganSWNotes} if $c_1(\fs)-\Lambda\in H^2(X;\ZZ)$ is a torsion
class. If $c_1(\fs)-\Lambda$ is not torsion, then this class is represented by
$(1/ 2\pi i)(\Tr(F_B)-F_{A_\Lambda})$ and, if $b_2^+(X)>0$ and the metric $g$ is
generic, then there are no zero-section solutions to the Seiberg-Witten
equations \eqref{eq:U1Monopole} by \cite[Proposition 6.3.1]{MorganSWNotes}.
\end{rmk}

The moduli space of Seiberg-Witten monopoles is defined by
\begin{equation}
\label{eq:SWModuli}
\begin{aligned}
M_{\fs} 
&= 
\{(B,\Psi)\in \tsC_{\fs}: \fS(B,\Psi)=0\}/\sG_{\fs}.
\\
M^{0}_{\fs} &= M_{\fs}\cap \sC^0_{\fs}.
\end{aligned}
\end{equation}
The usual Seiberg-Witten moduli space---obtained from our definition with
$\tau=\id_{\Lambda^+}$ and $\vartheta=0$---is compact
\cite[Corollary 3]{KMThom},
\cite[Proposition 6.4.1]{MorganSWNotes} and, for generic $\eta$,
the open subspace $M^{0}_{\fs}$ is a smooth
manifold of the expected dimension \cite[Lemma 5]{KMThom},
\cite[Proposition 6.2.2]{MorganSWNotes}.
In the subsections below we describe transversality and compactness
properties for the space $M_{\fs}$, as defined here.

\subsubsection{The deformation complex}
Given $(B,\Psi)\in\tsC_{\fs}$, the smooth map $\sG_{\fs}\to \tsC_{\fs}$
defined by \eqref{eq:SWGaugeAction} has differential at $\id_{W}$,
\begin{equation}
f\mapsto -d^{0}_{B,\Psi}f = (-df,f\Psi),
\qquad f \in L^2_{k+1}(X,i\RR).
\label{eq:d0SW}
\end{equation}
For all $(b,\psi)\in L^2_k(i\Lambda^1)\oplus L^2_k(W^+)$ and
$f\in L^2_{k+1}(X,i\RR)$ we have
\begin{align*}
\left((df,-f\Psi),(b,\psi)\right)_{L^2}
&= (df,b)_{L^2} + (-f\Psi,\psi)_{L^2} \\
&= (f,d^*b)_{L^2} - \left(f,\overline{\langle\Psi,\psi\rangle}\right)_{L^2} \\
&= (f,d^*b)_{L^2} - \left(f,\langle\psi,\Psi\rangle\right)_{L^2} \\
&= (f,d^*b)_{L^2} - \left(f,i\Imag\langle\psi,\Psi\rangle\right)_{L^2},
\end{align*}
noting that $\Omega^0(X,i\RR)$ has a real $L^2$ inner product, and thus
(compare  \cite[Lemma 4.5.5]{MorganSWNotes})
\begin{equation}
d^{0,*}_{B,\Psi}(b,\psi) = d^*b - i\Imag\langle\psi,\Psi\rangle.
\label{eq:d0SWL2adjoint}
\end{equation}
The Seiberg-Witten equations \eqref{eq:U1Monopole} define a
$\sG_{\fs}$-equivariant map $\fS:\tsC_{\fs}\to L^2_{k-1}(i\Lambda^+)\times
L^2_{k-1}(W^-)$ with differential at the point $(B,\Psi)$ given by
\begin{equation}
d^{1}_{B,\Psi}(b,\psi) = (D\fS)_{B,\Psi}(b,\psi)
=
\begin{pmatrix}
2d^+b-\tau\rho^{-1}(\psi\otimes\Psi^*+\Psi\otimes\psi^*)_0 \\
D_B\psi+\rho(\vartheta)\psi+\rho(b)\Psi
\end{pmatrix}.
\label{eq:d1SW}
\end{equation}
Since $\fS(s(B,\Psi)) = (\fS_1(B,\Psi),s\fS_2(B,\Psi)) =
s\cdot\fS(B,\Psi)$, the differential of the composition $\sG_{\fs}\to
L^2_k(i\Lambda^+)\times L^2_k(W^-)$ is given by
$$
f\mapsto d^{1}_{B,\Psi}\circ d^{0}_{B,\Psi}f
= (0,f\fS_2(B,\Psi)), \qquad f\in L^2_{k+1}(X,i\RR),
$$
and so $d^{1}_{B,\Psi}\circ d^{0}_{B,\Psi} = 0$ if and only if
$\fS_2(B,\Psi)=0$. In particular, if $\fS(B,\Psi)=0$ we have an
elliptic deformation complex
\begin{equation}
\label{eq:SWDeformationComplex}
\begin{CD}
L^2_{k+1}(X,i\RR)
@> d^{0}_{B,\Psi} >>
\begin{matrix}
L^2_k(i\Lambda^1) \\ \oplus \\ L^2_k(W^+)
\end{matrix}
@> d^{1}_{B,\Psi} >>
\begin{matrix}
L^2_{k-1}(i\Lambda^+) \\ \oplus \\
L^2_{k-1}(W^-)
\end{matrix}
\end{CD}
\end{equation}
with cohomology $H^{\bullet}_{B,\Psi}$. We write
\begin{equation}
\label{eq:SWDeformationOperator}
\sD_{B,\Psi} = d^{0,*}_{B,\Psi}+d^{1}_{B,\Psi}
\end{equation}
for the rolled-up deformation operator of the complex
\eqref{eq:SWDeformationComplex}.

The deformation complex \eqref{eq:SWDeformationComplex} is the same, up to
slight differences in the zeroth order terms of the differential
$d^{1}_{B,\Psi}$, as the deformation complex \cite[\S
4.6]{MorganSWNotes} for the usual Seiberg-Witten equations
and so has the same index. Therefore, provided the
zero locus $M^{0}_{\fs}$ of the section $\fS$ is regular we have
\begin{equation}
\label{eq:DimSWModuliSpace}
d_s(\fs)
=
\dim M^{0}_{\fs}
= 
\textstyle{\frac{1}{4}}\left(c_1(\fs)^2 - 2\chi - 3\sigma\right).
\end{equation}
With the deformation complex in place, we now turn to the issues of
compactness and transversality.

\subsubsection{Compactness and transversality}
We first observe that the standard arguments
\cite{KMThom}, \cite{MorganSWNotes}, \cite{Witten} establishing that the
Seiberg-Witten 
equations as stated in those references define a compact moduli space carry
over to show that the moduli space $M_{\fs}$ defined by
\eqref{eq:U1Monopole} is compact; the only slight change is the requirement
that $\tau$ be $C^0$-close to the identity on $\Lambda^+$:

\begin{prop}
\label{prop:U1MonopoleCompact}
Let $\fs=(\rho,W)$ be a \spinc structure on a closed, oriented,
Riemannian four-manifold, $(X,g)$, parameters $\eta\in\Om^0(i\Lambda^+)$
and $\vartheta \in \Omega^1(X,\CC)$, and perturbation
$\tau\in\Om^0(\GL(\Lambda^+))$ such that $\|\tau-\id_{\Lambda^+}\|_{C^0} <
\sixtyfourth$.  Then the moduli space
$M_{\fs}(g,\eta,\tau,\vartheta)$ is compact.
\end{prop}

Our next task is to establish an transversality result for
$M_{\fs}^0(g,\eta,\tau,\vartheta)$ with generic parameter $\tau$ analogous
to the one in \cite{KMThom} where $\eta$ is generic.

\begin{prop}
\label{prop:SmoothFamilyOfReducibles}
Let $\fs=(\rho,W)$
be a \spinc structure on a closed, oriented, Riemannian four-manifold,
$(X,g)$.
Then for any fixed
parameters $\eta\in\Om^0(i\Lambda^+)$ and $\vartheta \in \Omega^1(X,\CC)$
there is a first-category subset of $\Om^0(\GL(\Lambda^+))$ such that
for all $C^\8$ parameters $\tau$ in the complement of this subset, the
moduli space $M^{0}_{\fs}(g,\tau,\eta,\vartheta)$
is a smooth manifold of the expected dimension.
\end{prop}

\begin{proof}
We closely follow the method described in \cite[\S 5.1 \& \S 5.2]{FL1}.
Our Banach space of parameters is given by
$$
\sP = C^r(\GL(\Lambda^+)).
$$
We define an extended $\sG_{\fs}$-equivariant map,
$$
\underline{\fS} = (\underline{\fS}_1,\underline{\fS}_2):
\sP\times \tsC_{\fs}
\to L^2_{k-1}(i\Lambda^+)\times
L^2_{k-1}(W^-),
$$
by setting
$$
\underline{\fS}(\tau,B,\Psi)
=
\begin{pmatrix}
\Tr(F^+_B)-\tau\rho^{-1}(\Psi\otimes\Psi^*)_0-\eta \\
D_B\Psi+\rho(\vartheta)\Psi
\end{pmatrix}.
$$
The parametrized moduli space $\fM_{\fs}$ is then
$(\underline{\fS})^{-1}(0)/\sG_{\fs}\subset\sC_{\fs}$ and $\fM_{\fs}^0 =
\fM_{\fs}\cap(\sP\times\sC_{\fs}^0)$. The map $\underline{\fS}$ has
differential at the point $(\tau,[B,\Psi])\in \fM^0_{\fs}$ given by
$$
D\underline{\fS}(\de\tau,b,\psi)
=
\begin{pmatrix}
2d^+b-\tau\rho^{-1}(\psi\otimes\Psi^*+\psi\otimes\Psi^*)_0
+ (\de\tau)\tau\rho^{-1}(\Psi\otimes\Psi^*)_0\\
(D_B+\rho(\vartheta))\psi+\rho(b)\psi
\end{pmatrix},
$$
where $(b,\psi)\in L^2_k(i\Lambda^1)\oplus L^2_k(W^+)$.  Suppose
$$
(c,\varphi)\in L^2_{k-1}(i\Lambda^+)\oplus L^2_{k-1}(W^-)
$$
is $L^2$ orthogonal to the image of $D\underline{\fS}$. We may
assume without loss of generality that $(\tau,B,\Psi)$ is a $C^r$
representative of the point $(\tau,[B,\Psi])$ (see, for example,
\cite[Proposition 3.7]{FL1}) and so, as
$(D\underline{\fS})(D\underline{\fS})^*(c,\varphi)=0$, elliptic
regularity implies that $(c,\varphi)$ is $C^r$. Then
\begin{equation}
0
=
(D\underline{\fS}(\delta\tau,0,0),(c,\varphi))_{L^2}
=
((\delta\tau)\tau(\Psi\otimes\Psi^*)_0,c)_{L^2},
\label{eq:L2Orthog}
\end{equation}
for all $\delta\tau\in C^r(\gl(\Lambda^+))$, which yields the pointwise
identity
\begin{equation}
\langle(\delta\tau)_x\tau_x(\Psi_x\otimes\Psi_x^*)_0,c_x\rangle_x = 0,
\qquad x\in X.
\label{eq:PointwiseOrthog}
\end{equation}
If $c_x\ne 0$ for some $x\in X$ (and thus $c\ne 0$ on a non-empty open
neighborhood in $X$), then \eqref{eq:PointwiseOrthog} implies that $\Psi =
0$ on the non-empty open subset $\{c\ne 0\}\subset X$. Aronszajn's
theorem then implies that $\Psi\equiv 0$ on $X$, contradicting our assumption
that $(\tau,B,\Psi)$ is a point in $\fM_{\fs}^0$. Hence, we must have
$c\equiv 0$ on $X$, so \eqref{eq:L2Orthog} now yields
$$
0
=
(D\underline{\fS}(0,b,0),(0,\varphi))_{L^2}
=
(\rho(b)\Psi,\varphi)_{L^2},
$$
for all $b \in L^2_k(i\Lambda^1)$.  Since $\Psi\not\equiv 0$, this
implies that $\varphi\equiv 0$, just as in the proof of \cite[Lemma
5]{KMThom}.  Hence, $\Coker (D\underline{\fS})_{\tau,[B,\Psi]} = 0$
at each point $(\tau,[B,\Psi])$ in $(\underline{\fS})^{-1}(0)$ and so
the parametrized moduli space $(\underline{\fS})^{-1}(0)=\fM_{\fs}^0$
is a smooth Banach submanifold of $\sP\times\sC_{\fs}^0$. Therefore, the
Sard-Smale theorem, in the form of \cite[Proposition 4.3.11]{DK} (see the
proof of Corollary 5.3 in \cite{FL1}), implies that there is a
first-category subset of $C^r(\GL(\Lambda^+))$ such that for all $C^r$
parameters $\tau$ in the complement of this subset, the zero locus
$M_{\fs}^0(\tau)$ is a regular submanifold of $\sC_{\fs}^0$.  Lastly, we
can constrain the parameter $\tau$ to be $C^\8$ (and not just $C^r$) by the
argument used in \cite[\S 5.1.2]{FL1}---see \cite[\S 8.4]{SalamonSWBook}.)
\end{proof}

The following lemma relates the Seiberg-Witten moduli spaces defined with
the perturbations of \cite{KMThom}, \cite{MorganSWNotes},
\cite{Witten} to those used in this article.

\begin{lem}
\label{lem:CompareParameters}
Let $(X,g)$ be a closed, oriented, Riemannian four-manifold with
$b_2^+(X)>0$ and \spinc structure $\fs$.  
\begin{enumerate}
\item Suppose that the moduli spaces
  $M^{0}_{\fs}(g,\id_{\Lambda^+},\eta',0)$ and
  $M^{0}_{\fs}(g,\tau,\eta,\vartheta)$ are regular zero loci of the
  Seiberg-Witten maps $\underline{\fS}$ in \eqref{eq:U1Monopole} for these
  perturbation parameters.
  If $b^2_+(X)>1$, then these moduli spaces are related by an
  oriented, smooth cobordism.
\item Suppose that $M_{\fs}(g,\id_{\Lambda^+},\eta',0)$ contains no
  zero-section solutions. Then $M_{\fs}(g,\tau,\eta,\vartheta)$ contains no
  zero-section solutions for any triple $(\tau,\vartheta,\eta)$ which is
  $C^r$-close enough to $(\id_{\Lambda^+},0,\eta')$.
\end{enumerate}
\end{lem}

\begin{proof}
If $\sP$ denotes the Banach space of $C^r$ perturbation parameters
$(\tau,\eta,\vartheta)$, then a generic, smooth path in $\sP$ joining
$(\id_{\Lambda^+},\eta',0)$ to $(\tau,\eta,\vartheta)$ induces a smooth,
oriented cobordism in $\sP\times\sC_{\fs}^0$ joining
$M^{0}_{\fs}(g,\id_{\Lambda^+},\eta',0)$ to 
$M^{0}_{\fs}(g,\tau,\eta,\vartheta)$. 
This proves Assertion (1).

If the compact space $M_{\fs}(g,\id_{\Lambda^+},\eta',0)$ is contained in
the open subset $\sC_{\fs}^0\subset\sC_{\fs}$, then the same holds for
$M_{\fs}(g,\tau,\eta,\vartheta)$, that is, the latter subspace also
contains no zero-section pairs, for any triple $(\tau,\vartheta,\eta)$
which is $C^r$-close enough to $(\id_{\Lambda^+},0,\eta')$. This proves
Assertion (2).
\end{proof}

\subsection{Cohomology ring of the configuration space of \spinc pairs}
\label{subsec:AbelianCohom}
We next compute the cohomology ring of the configuration space
$\sC_{\fs}^0$ and describe the cohomology classes used to define
Seiberg-Witten invariants and arising in the calculation (see \S
\ref{subsec:TopOfNormal}) of the Chern classes of certain vector bundles
over $M_{\fs}$.  For this purpose we may make use of the fact since
$\U(1)$ is Abelian, we have a global slice theorem for $\sC_{\fs}$ modulo
the harmonic gauge transformations.
Computations of the cohomology ring of $\sC_{\fs}^0$ have
already appeared in, for example \cite[Chapter 7]{SalamonSWBook}
and \cite[Lemma 5]{OTWall}.  The version we present here
is convenient for the computation of the Chern class of the
universal line bundle in Lemma \ref{lem:MuMapDescription}.

\subsubsection{Harmonic gauge transformations}
\label{subsubsec:AbelConfig}
A gauge transformation $s\in\Map(X,S^1)$ is {\em harmonic\/} if
$s^{-1}ds\in L^2_k(i\Lambda)$ satisfies $d^*(s^{-1}ds)=0$ (we always have
$d(s^{-1}ds)=0$).  Every component of $\Map(X,S^1)$ contains a harmonic
representative which is unique up to multiplication by a constant element
of $S^1$. The harmonic elements of $\Map(X,S^1)$ form a subgroup,
\begin{equation}
\label{eq:DefineHarmonicGT}
\hatsG_{\fs}=\{s\in\Map(X,S^1): d^*(s^{-1}ds) = 0\}
 \cong S^1 \times H^1(X,i\ZZ),
\end{equation}
where, given a point $p\in X$, the isomorphism is defined by
$s\mapsto (s(p),s^{-1}ds)$. Conversely, we have a map
$$
\bh_p:H^1(X;i\ZZ)\to\Map(X,S^1)
$$
defined by setting $\bh_p(\beta)$ equal to the unique harmonic
gauge transformation satisfying $\bh_p(\beta)(p)=1$
(see, for example, \cite[\S 5.4]{SalamonSWBook}).
The group $S^1\times H^1(X;i\ZZ)$ acts on $\tsC_{\fs}$ by
\begin{equation}
\label{eq:HarmonicGaugeGrpAction}
\left((e^{i\theta},\beta),(B,\Psi)\right)
\mapsto
(B-\beta\,\id_W,e^{i\theta}\bh_p(\beta)\Psi),
\end{equation}
where $(e^{i\theta},\beta)\in S^1\times H^1(X;i\ZZ)$ and $(B,\Psi)\in
\tsC_{\fs}$.  We then have the following global ``slice result'' for
$\sC_{\fs}$:

\begin{lem}
\label{lem:S1SemiSlice}
Let $B_0$ be an $L^2_k$ spin connection on $W$. Then the inclusion
\begin{equation}
\label{eq:DefineS1SemiSlice}
\hatsC_{\fs}
=
(B_0 + (\Ker d^*)\id_W)
\times L^2_k(W^+)
\hookrightarrow \tsC_{\fs},
\end{equation}
where $\Ker d^*\subset L^2_k(i\Lambda^1)$,
is equivariant with respect to the action \eqref{eq:HarmonicGaugeGrpAction}
of $S^1\times H^1(X;i\ZZ)\subset\sG_{\fs}$ and descends to a homeomorphism
$$
(B_0 + (\Ker d^*)\id_W)
\times_{S^1\times H^1(X;i\ZZ)} L^2_k(W^+)
\to \sC_{\fs},
$$
and a diffeomorphism on $\hatsC_{\fs}^0$, the complement in
$\hatsC_{\fs}$ of the space of zero-section pairs.
\end{lem}

\begin{proof}
Any $L^2_k$ spin connection on $W$ can be written as $B=B_0+b\,\id_W$, where
$b\in L^2_k(i\La^1)$.  An element $s$ of the gauge group $\sG_{\fs}$ acts on
this representation of $B$ by sending $b$ to $b-s^{-1}ds$.
The argument in \cite[\S 7.2]{SalamonSWBook}
(also see \cite[pp. 54--55]{DK}) implies that there
is a solution of the equation $d^*(b- s^{-1}ds)=0$, which is
unique up to a harmonic gauge transformation, so the
map on quotients is surjective. If $(B_j,\Psi_j)$, $j=1,2$, are pairs in
$\tsC_{\fs}$ such that $d^*(B_j-B_0)=0$ and $s\in\sG_{\fs}$ satisfies
$s(B_1,\Psi_1) = (B_2,\Psi_2)$, that is
$$
(B_1-s^{-1}ds\,\id_W,s\Psi_1) = (B_2,\Psi_2),
$$
then $d^*(B_2-B_1)=0$ implies $d^*(s^{-1}ds)=0$, so $s$ is harmonic and
the map on quotients is injective.
\end{proof}

\subsubsection{Cohomology ring of the configuration space}
\label{subsubsec:AbelianCohom}
We now use the reduction of the gauge group from $\sG_{\fs}$to
$\hatsG_{\fs}$ to compute the integral cohomology ring of $\sC^0_{\fs}$.
We begin by defining the universal complex line bundle over
$\sC^0_{\fs}\times X$, 
\begin{equation}
\label{eq:SWUnivLineBundle}
\LL_{\fs}
=
\tsC^0_{\fs}\times_{\sG_{\fs}}\underline{\CC}
\cong
\hatsC^0_{\fs}\times_{\hatsG_{\fs}}\underline{\CC}.
\end{equation}
where $\underline{\CC}=X\times\CC$ and the action of $\sG_{\fs}$ is given
for $s\in\sG_{\fs}$, $x\in X$ and $z\in\CC$ by
\begin{equation}
\label{eq:DefineUnivSWBundleAction}
\left( s, (B,\Psi),(x,z)\right)
\to
\left( s_*(B,\Psi),(x,s(x)^{-1}z)\right).
\end{equation}
As in \cite[\S 7.4]{SalamonSWBook}, we define cohomology classes on
$\sC^0_{\fs}$ (the ``Seiberg-Witten $\mu$-classes'') by
\begin{equation}
\mu_{\fs}: H_\bullet(X;\RR)\to H^{2-\bullet}(\sC^0_{\fs};\RR),
\qquad
\mu_{\fs}(\alpha)
=
c_1(\LL_{\fs})/\alpha,
\label{eq:SWMuMap}
\end{equation}
where $\alpha$ is either the positive generator $x\in H_0(X;\ZZ)$ or a
class $\gamma\in H_1(X;\RR)$. Let $\sG^p_{\fs}$ be the subgroup of maps
$s\in\sG_{\fs}$ with $s(p)=1$ and note that
\begin{equation}
\label{eq:DefineBasedSWGaugeTransformations}
\sG_{\fs}=\sG^p_{\fs}\times S^1.
\end{equation}
The following alternative characterization of $\mu_{\fs}(x)$ in appears in
much of the literature on Seiberg-Witten invariants.

\begin{lem}
\label{lem:FramedSWConfigurationSpace}
Let $\sC^{p,0}_{\fs}=\tsC^0_{\fs}/\sG^p_{\fs}$ be the configuration
space of framed pairs.  Then
\begin{equation}
\label{eq:SWMuMapAndFramedConfig}
c_1(\sC^{p,0}_{\fs}\times_{S^1}\CC)=\mu_{\fs}(x).
\end{equation}
\end{lem}

\begin{proof}
The class $\mu_{\fs}(x)$ is the first Chern class of
$\LL_{\fs}$ restricted to $\sC^0_{\fs}\times\{p\}$.
Equation \eqref{eq:SWMuMapAndFramedConfig} then follows from the splitting
$\sG_{\fs}=\sG^p_{\fs}\times S^1$.
\end{proof}

Using the isomorphism of universal line bundles in
\eqref{eq:SWUnivLineBundle}, we shall now give another description of the
cohomology classes $\mu_{\fs}(x)$ and $\mu_{\fs}(\gamma)$, and show they
generate the cohomology ring $H^*(\sC^0_{\fs};\ZZ)$.  We define
\begin{equation}
\label{eq:DefineJacobian}
\Jac(X)
=
H^1(X;i\RR)/H^1(X;i\ZZ).
\end{equation}
We then have:

\begin{lem}
\label{lem:Retraction}
Let $\fs$ be a \spinc structure on a smooth, closed, oriented four-manifold
$X$.  Then, there is an $S^1$-equivariant retraction
$\br:\sC^{p,0}_{\fs}\to \Jac(X)\times
\left(L^2_k(W^+) - \{0\}\right)$ where
$S^1$ acts trivially on $\Jac(X)$ and by complex multiplication
on $L^2_k(W^+)$.
\end{lem}

\begin{proof}
{}From the Hodge decomposition,
$$
L^2_k(i\Lambda^1) = \bH^1(X;i\RR)\oplus \Ran d\oplus \Ran d^*,
$$
where $\bH^1(X;i\RR)$ denotes the harmonic, imaginary-valued one-forms on $X$,
we see that
$$
\Ker d^* = (\Ran d)^\perp = \bH^1(X;i\RR)\oplus \Ran d^*.
$$
Thus, we can write the pre-configuration space $\hatsC^0_{\fs}$
as
$$
\hatsC^0_{\fs}
=
\left(B_0+\bH^1(X;i\RR)\oplus\Ran d^*\right)
\times
\left(L^2_k(W^+) - \{0\}\right).
$$
As usual, the factor $S^1$ of the harmonic gauge group
$\hatsG_{\fs} = S^1\times H^1(X;i\ZZ)$ acts
trivially on $B_0+\bH^1(X;i\RR)\oplus\Ran d^*$ and
by complex multiplication on $L^2_k(W^+) - \{0\}$.
An element $\beta\in\bH^1(X;i\ZZ)\cong H^1(X;i\ZZ)$ acts on an element
$(B,\Psi) = (B_0+b'+b'',\Psi)$ of $\hatsC^0_{\fs}$
(where $b'\in\bH^1(X;i\RR)$ and $b''\in\Ran d^*$) by
$$
(\beta, (B_0+b'+b'',\Psi))
\mapsto
(B_0+(b'+\beta)+b'',\bh_p(\beta)\Psi).
$$
Thus, replacing $b''\in\Ran d^*$ with $tb'', 0\leq t\leq 1$,
defines a $\hatsG_{\fs}$-equivariant retraction of
$\hatsC^0_{\fs}$ onto
$$
\left(B_0+\bH^1(X;i\RR)\right)
\times
\left(L^2_k(W^+)-\{0\}\right).
$$
This yields a deformation retraction of $\hatsC^0_{\fs}$ onto the
projectivization of the complex vector bundle
\begin{equation}
\label{S1Config}
\left(B_0+\bH^1(X;i\RR)\right)\times_{H^1(X;i\ZZ)}L^2_k(W^+)
\to
\Jac(X).
\end{equation}
By Kuiper's theorem (see Corollary \ref{KuiperSmooth}) this vector bundle is
trivial, which completes the proof.
\end{proof}

\begin{defn}
\label{defn:RelatedBasis}
Let $\{\gamma_i\}$ be a basis for $H_1(X;\ZZ)/\Tor$ and let $\{\gamma_i^*\}$ be
the dual basis for $H^1(X;\ZZ)$, so $\langle
\gamma_i^*,\gamma_j\rangle=\delta_{ij}$.  Then
$\{{\sqrt {-1}}\gamma_j^*\}$ generates $\pi_1(\Jac (X))$.  We will
write these elements of $\pi_1(\Jac(X))$ as $\gamma^J_j$ to avoid
confusion.  Let $\gamma_j^{J,*}\in
H^1(\Jac(X);\ZZ)=\Hom(\pi_1(\Jac(X)),\ZZ)$ be defined by
$\langle\gamma_i^{J,*},\gamma_j^J\rangle=\delta_{ij}$.
Given such a basis $\{\gamma_i\}$ for $H_1(X;\ZZ)/\Tor$,
we call $\{\gamma_i^J\}$ and $\{\gamma_i^{J,*}\}$ the
{\em related generators\/} and {\em basis\/} for $\pi_1(\Jac(X))$ and
$H^1(\Jac(X);\ZZ)$ respectively.
\end{defn}

Note that $H^1(X;\ZZ)$ and $H^1(\Jac(X);\ZZ)$ are free $\ZZ$-modules by the
universal coefficient theorem \cite[Theorem 5.5.3]{Spanier}.

\begin{cor}
\label{cor:CohomOfS1Config}
Continue the notation and assumptions of Lemma \ref{lem:Retraction}.
Let $\br:\sC^0_{\fs}\to \Jac(X)\times\left(L^2_k(W^+) -
\{0\}\right)/S^1$ be the retraction given by Lemma \ref{lem:Retraction}.
Then the cohomology ring of $\sC^0_{\fs}$ is generated by
$\{\br^*h,\br^*\gamma_i^{J,*}\}$, where
$\{\gamma_1^{J,*},\dots,\gamma_{b_1}^{J,*}\}$ is a
basis for $H^1(\Jac(X);\ZZ)$ and $h$
the first Chern class of the $S^1$ bundle given by
$$
\begin{CD}
\left(B_0+\bH^1(X;i\RR)\right)\times_{H^1(X;i\ZZ)}
\left( L^2_k(W^+)- \{0\}\right)
\\
@VVV \\
\left(B_0+\bH^1(X;i\RR)\right)\times_{H^1(X;i\ZZ)}
\left( L^2_k(W^+)- \{0\}\right)/S^1.
\end{CD}
$$
\end{cor}

\subsubsection{Seiberg-Witten $\mu$-classes and generators of the
cohomology ring of the configuration space} Finally, we show that the
Seiberg-Witten $\mu$-classes can be expressed in terms of the classes
$\br^*h$ and $\br^*\gamma_i^*$ defined in Corollary
\ref{cor:CohomOfS1Config}.  For this purpose, we use the complex line
bundle,
\begin{equation}
\label{eq:JacobianBundle}
\bDelta = \bH^1(X;i\RR)\times X\times_{H^1(X;i\ZZ),\bh_p}\CC
\to
\Jac(X)\times X,
\end{equation}
whose first Chern class we now compute.

\begin{lem}
\label{lem:CanonicalLineBundle}
Continue the notation of Definition \ref{defn:RelatedBasis}.  Then
$c_1(\bDelta)=\sum_{i=1}^{b_1(X)}\gamma_i^{J,*}\times\gamma_i^*$.
\end{lem}

\begin{proof}
  Since $H^1(X;\ZZ)$ and $H^1(\Jac(X);\ZZ)$ are free, we
  can write $c_1(\bDelta)=c_{0,2}+c_{1,1}+c_{2,0}$ with respect to the
  decomposition given by the Kunneth formula:
\begin{align*}
H^2(\Jac(X)\times X;\ZZ)
& \cong
\bigoplus_{i=0}^2 H^i(\Jac(X);\ZZ)\otimes H^{2-i}(X;\ZZ),
\\
 c_{i,2-i} &\in H^i(\Jac(X);\ZZ)\otimes H^{2-i}(X;\ZZ).
\end{align*}
For $[B]\in\Jac(X)$ and $p\in X$,
the restrictions of $\bDelta$
to $\{[B]\}\times X$ and to $\Jac(X)\times
\{p\}$ are trivial so $c_1(\bDelta)=c_{1,1}$.
We now calculate the first Chern
class of $\bDelta$ restricted to the tori $\gamma_i^J\times\gamma_j$.
Let $\tilde\gamma_i^J(t)$ be the path in $\bH^1(X;i\RR)$,
starting at zero and covering $\gamma_i^J$.
Let $s:X\to S^1$ be an element of $\hatsG_{\fs}$
satisfying $s(p)=1$ and $s^{-1}ds=\tilde\gamma_i^J(1)$.
Thus, the degree of the map $s:\gamma_j\to S^1$
is $\delta_{ij}$.
The restriction of the line bundle $\bDelta$ to
$\gamma_i^J\times\gamma_j$ is given by taking the trivial
bundle over $I\times \gamma_j$ and making the identification
$$
(0,q,z) \sim (1,q,s^{-1}(q)(z)),\quad
\text{where}\quad q\in\gamma_j, z\in\CC.
$$
Since the degree of the map $s:\gamma_j\to S^1$ is
$\delta_{ij}$, we see
$$
c_1(\bDelta|_{\gamma^J_i\times\gamma_j})=\delta_{ij}
\gamma_i^{J,*}\times\gamma_j^*.
$$
This completes the proof.
\end{proof}

We now describe the relation between the cohomology classes
of Corollary \ref{cor:CohomOfS1Config} and the $\mu_{\fs}$-map.

\begin{lem}
\label{lem:MuMapDescription}
Let $\{\gamma_i\}$ be a basis for $H_1(X;\ZZ)/\Tor$, let $\{\gamma_i^*\}$
be the related basis for $H^1(\Jac(X);\ZZ)$, and let $x\in H_0(X;\ZZ)$ be a
generator.  Let $\br_J:\sC^0_{\fs}\to\Jac(X)$ be the composition of the
retraction $\br$ defined in Lemma \ref{lem:Retraction} with the projection
to $\Jac(X)$.  Then
\begin{equation}
\label{eq:ChernClassOfUniversalSW}
c_1(\LL_{\fs})= h\times 1 + (\br_J\times \id_X)^*c_1(\bDelta),
\end{equation}
and thus
$$
\mu_{\fs}(x)= \br^*h,\qquad
\mu_{\fs}(\gamma_i)=\br^*_J\gamma_i^{J,*}.
$$
\end{lem}

\begin{proof}
Because $\br$ is a retraction, it suffices to compute the restriction
of the universal line bundle $\LL_{\fs}$ to
the image of the retraction
where $\LL_{\fs}$ is given by
\begin{equation}
\label{eq:RestrictedUniversalSWBundle}
\bH^1(X;i\RR)\times \left(L^2_k(W^+)-\{0\}\right)
\times_{S^1\times H^1(X;i\ZZ)}\underline{\CC}.
\end{equation}
The bundle
\begin{equation}
\label{eq:FramedSWConfigSpace}
\begin{CD}
\bH^1(X;i\RR)\times_{H^1(X;i\ZZ)} \left(L^2_k(W^+)-\{0\}\right)
\times_{S^1}\CC
\\
@VVV \\
\bH^1(X;i\RR)\times_{H^1(X;i\ZZ)} \left(L^2_k(W^+)-\{0\}\right)/S^1
\end{CD}
\end{equation}
is the restriction of the framed configuration space,
$\sC^{p,0}_{\fs}\times_{S^1}\CC$ to the image of $\br$ and thus has first
Chern class $h$.  Then the restriction of $\LL_{\fs}$ in equation
\eqref{eq:RestrictedUniversalSWBundle} is isomorphic to the tensor product
of the bundle \eqref{eq:FramedSWConfigSpace} with the pullback of the
bundle $\bDelta$ by the map $\br_J\times\id_X$.  Thus
$c_1(\LL_{\fs})=h\times 1 +(\br_J\times\id_X)^*c_1(\bDelta)$.  The
computations of $\mu_{\fs}(\beta)$ then follow from the definition of the
map $\mu_{\fs}$ in equation \eqref{eq:SWMuMap} and the computation of
$c_1(\bDelta)$ in Lemma \ref{lem:CanonicalLineBundle}.
\end{proof}


\section{Singularities and their links}
\label{sec:Singularities}
As we shall see in this section, the moduli
space $\sM_{\ft}$ of solutions to the perturbed $\PU(2)$ monopole equations
\eqref{eq:PT} is a smoothly stratified space.
According to Theorem \ref{thm:Transversality}, the subspace
$\sM^{*,0}_{\ft}$ of solutions to
\eqref{eq:PT} which are neither zero-section nor reducible
is a smooth manifold of the expected dimension and this comprises the
top stratum of $\sM_{\ft}$. In \S \ref{subsec:TypesOfSing} we classify
the $\PU(2)$ monopoles where $S^1$ acts trivially, noting that these
are given by the zero-section or reducible solutions and comprise the
singular points of $\sM_{\ft}$. We exclude the possibility of $\PU(2)$
monopoles with flat associated $\SO(3)$ connections by a suitable
choice of $w_2(\ft)$.  We then analyze these lower strata of solutions
in $\sM_{\ft}$ which are either zero-section or reducible or both and
show that they are identified with the moduli space $M^w_{\ka}$ of
anti-self-dual connections on the associated $\SO(3)$ bundle
$\fg_{\ft}$ and moduli spaces of Seiberg-Witten monopoles,
$M_{\fs}$ (see \S \ref{subsec:ReducibleSingularities}). We
describe the link in $\sM^{*,0}_{\ft}$ of the stratum $M^w_\ka$ in \S
\ref{subsec:ASDLink} and the links in $\sM^{*,0}_{\ft}$ of the strata
$M_{\fs}$ in \S \ref{subsec:LinkOfReduciblesLocal} and \S
\ref{subsec:LinkOfReduciblesGlobal}. For $M^w_{\ka}$, it will
suffice to describe the normal cone (in the sense of
\cite[p. 41]{GorMacPh}) at a generic point in $M^w_\ka$, while in the
case of the link of a stratum $M_{\fs}$, a global description is
required. Finally, in \S \ref{subsec:TopOfNormal} we compute the Chern
character of the normal bundle of $M_{\fs}$ with respect to a
finite-dimensional, open, $S^1$-equivariant
smooth manifold containing a neighborhood of the image of
$M_{\fs}$ in $\sM_{\ft}$.

\subsection{Classification of fixed points under the circle action}
\label{subsec:TypesOfSing}
In this section we classify the fixed points in $\sM_{\ft}$ under the circle
action induced by scalar multiplication on $V$.

For pairs $(A,\Phi)$ in
$\tsC_{\ft}$, it is useful to distinguish between two kinds of circle action.
First, $S^1$ can act on $V$ by scalar multiplication: 
\begin{equation}
\label{eq:ScalarMultAction}
S^1\times V\to V,\quad (e^{i\theta},\Phi)\mapsto e^{i\theta}\Phi.
\end{equation}
Second, if $V=W\oplus W\otimes L$ is a direct sum of Clifford modules,
then $S^1$ can act by scalar multiplication on the factor $W\otimes L$ and
the identity on $W$:
\begin{equation}
\label{eq:VReducibleAction}
S^1\times V\to V,\quad (e^{i\theta},\Psi\oplus\Psi')
\mapsto \Psi\oplus e^{i\theta}\Psi',
\end{equation}
where $\Psi\in C^\8(W)$ and $\Psi'\in C^\8(W\otimes L)$.
With respect to the splitting $V=W\oplus W\otimes L$, the two actions are
related by
$$
\begin{pmatrix}1 & 0 \\ 0 & e^{i2\theta}\end{pmatrix}
=
e^{i\theta}u,
\text{ where }
u = \begin{pmatrix}e^{-i\theta} & 0 \\ 0 & e^{i\theta}\end{pmatrix}
\in\sG_{\ft},
$$
and so, when we pass to the induced circle actions on the quotient
$\sC_{\ft}=\tsC_{\ft}/\sG_{\ft}$, we obtain the same space $\sC_{\ft}/S^1$ for
both circle actions.

If we write $V=W\otimes E$, where $E=\underline{\CC}\oplus L$,
the circle action \eqref{eq:VReducibleAction} on $V$ is equivalent to the
one induced by scalar multiplication on the complex line bundle $L$ and
the trivial action on $\underline{\CC}$. Hence, the induced circle action
\eqref{eq:SO(3)ReducibleAction} on $\su(E)\subset\su(V)$ with respect to
the decomposition \eqref{eq:OrthogDecompsuV} takes the form
\begin{equation}
\label{eq:SO(3)ReducibleAction}
S^1\times \su(E)\to \su(E),\quad (e^{i\theta},(\nu,z))
\mapsto (\nu,e^{i\theta}z),
\end{equation}
with respect to the isomorphism $\su(E)\cong i\underline{\RR}\oplus L$ in
Lemma \ref{lem:su(E)Split}: the circle acts as the identity on
$i\underline{\RR}$ and by scalar multiplication on the complex line bundle
$L$. On the other hand, the circle action \eqref{eq:ScalarMultAction} on
$V$ induces the trivial action on both $\su(V)$ and $\su(E)$.

\begin{prop}
\label{prop:ClassificationOfStabilizers}
Let $\ft=(\rho,V)$ be a \spinu structure over a closed, oriented, smooth
four-manifold $X$, with $b_2^+(X)\geq 1$ and generic Riemannian metric. If
$\hat A$ is a non-flat connection on $\fg_{\ft}$, then $[A,\Phi]\in \sM_{\ft}$
is a fixed point with respect to the $S^1$ action on $\sM_{\ft}$ if and only if
one of the following hold:
\begin{enumerate}
\item
The connection $\hat A$ is anti-self-dual, irreducible, and $\Phi\equiv
0$. The pair $(A,\Phi)$ is a fixed point with respect to the circle action
\eqref{eq:ScalarMultAction}. 
\item
The pair $(A,\Phi)$ is reducible with respect to a splitting $V=W\oplus
W\otimes L$ but $\Phi\not\equiv 0$,
$$
(A,\Phi)
=
(B\oplus B\otimes A_L,\Psi\oplus 0).
$$
The pair $(A,\Phi)$ is a fixed point with respect to the circle action
\eqref{eq:VReducibleAction}.  The connection $\hat A$ on $\fg_{\ft}$ is
reducible with respect to the splitting $\fg_{\ft}\cong
i\underline{\RR}\oplus L$, with $\hat A=d_\RR\oplus A_L$.
\end{enumerate}
\end{prop}

\begin{proof}
Suppose $[A,\Phi]\in \sM_{\ft}$ is a fixed point of the $S^1$ action on
$\sM_{\ft}$. Consequently, there is an element $e^{i\theta}\neq \pm 1$ such that
$e^{i\theta}[A,\Phi]=[A,\Phi]$ and hence a gauge transformation
$u\in\sG_{\ft}$ such that
$$
(A,e^{i\theta}\Phi) = u(A,\Phi) = (u(A),u\Phi). 
$$
Thus, $u$ is in the stabilizer of $A$, and hence that of $\hat A$, and
$e^{i\theta}\Phi=u\Phi$.

(1) If $\hat A$ has the trivial stabilizer $\{\pm \id_V\}$ in $\sG_{\ft}$,
then $e^{i\theta}\Phi=\pm\Phi$ and so $\Phi\equiv 0$ because
$e^{i\theta}\neq \pm 1$. The connection $\hat A$ is irreducible
\cite[Theorem 10.8]{FU} since it is an $\SO(3)$ connection with minimal
stabilizer in $\sG_{\ft}$. The curvature equation in \eqref{eq:PT} implies
that $F^+_{\hat A}= 0$, so $\hat A$ is anti-self-dual. The pair $(A,\Phi)$
is fixed by the $S^1$ action \eqref{eq:ScalarMultAction}.

(2) If $\hat A$ has non-trivial stabilizer in $\sG_{\ft}$, then $\hat A$ is
reducible with respect to a splitting $\fg_{\ft}\cong i\underline{\RR}\oplus
L$ \cite[Theorem 10.8]{FU}, for some complex line bundle $L$, and takes
the form $\hat A = d_\RR\oplus A_L$. The connection $\hat A$ has stabilizer
$\SO(2)\cong S^1$ acting on $L$ by complex multiplication and trivially
on $i\underline{\RR}$. Lemma \ref{lem:ReducibleSpinu} implies
that $A$ is reducible with respect to the splitting $V=W\oplus W\otimes
L$, taking the form $A=B\oplus B\otimes A_L$.

If $u=\pm \id_V$, then $\Phi\equiv 0$ just as in case (1) and so $\hat A$ would
be a reducible, anti-self-dual $\SO(3)$ connection on $\fg_{\ft}$.  But
$b_2^+(X) > 0$ and the Riemannian metric $g$ on $X$ is generic, so this
possibility is excluded by Corollary 4.3.15 in \cite{DK} (the four-manifold
$X$ does not need to be simply connected).

Hence, we must have $u\neq\pm\id_V$ and $\Phi\not\equiv 0$.  Because
$\Ad(u)$ stabilizes $\hat A$, it induces an action on $L$
via complex multiplication by some $e^{-i2\mu}\neq 1$ and the trivial
action on $i\underline{\RR}$, as in
\eqref{eq:SO(3)ReducibleAction}. Hence, with respect to the splitting
$V=W\oplus W\otimes L$, the gauge transformation $u$ takes the form
$$
u=\begin{pmatrix}e^{i\mu} & 0 \\ 0 & e^{-i\mu}\end{pmatrix}.
$$
Because $e^{-i\theta}u\Phi = \Phi$ and writing $\Phi = \Psi\oplus\Psi'$
with respect to the splitting $V=W\oplus W\otimes L$, the section $\Phi$ is
fixed by
$$
e^{-i\theta}u
=
\begin{pmatrix}e^{i(\mu-\theta)} & 0 \\ 0 & e^{-i(\mu+\theta)}\end{pmatrix}.
$$
If $\Psi\not\equiv 0$, then $\mu=\theta\pmod{2\pi}$ and we must have
$\Psi' \equiv 0$; conversely, if $\Psi'\not\equiv 0$, then
$\mu=-\theta\pmod{2\pi}$ and we must have $\Psi \equiv 0$. The two cases
differ only by how the factors in the splitting of $V$ are labeled, so
we can assume that $\Psi'\equiv 0$ and $\Psi\in C^\8(W^+)$.  In
particular, the pair $(A,\Phi)$ is fixed by the $S^1$ action
\eqref{eq:VReducibleAction}.
\end{proof}

Our proof that $\Phi$ is a section of $W^+$ when $\hat A$ is a non-flat
reducible connection could be replaced by an appeal to Lemma 5.22 in
\cite{FL1}, but the argument here seems more direct.

Proposition \ref{prop:ClassificationOfStabilizers} only classifies the
fixed points $[A,\Phi]\in \sM_{\ft}$ of the circle action under the assumption
that $\hat A$ is not flat. However, the following lemma gives a simple
condition which guarantees that there will be no pairs in $\sM_{\ft}$ with
flat associated $\SO(3)$ connections.  Because it relies only on a choice
of integral class $w\pmod{2}$, the lemma applies simultaneously to all
\spinu structures $\ft=(\rho,V)$ with a fixed $w_2(\ft)\equiv
w\pmod{2}$ and all oriented, orthogonal 3-plane bundles $F$ with
$w_2(F)\equiv w\pmod{2}$, and thus simultaneously to all levels of the
Uhlenbeck compactifications $\bar M_\kappa^w$ and $\bar\sM_{\ft}$.

\begin{lem}
\cite[p. 226]{MorganMrowkaPoly}
\label{lem:MorganMrowka}
Let $X$ be a closed, oriented four-manifold. Then
\begin{itemize}
\item If $w\in H^2(X;\ZZ)$ and $e\in H_2(X;\ZZ)$ is a spherical class such
  that $\langle w,e\rangle\not\equiv 0\pmod{2}$, then no oriented,
  orthogonal 3-plane bundle $F$ over $X$ with $w_2(F)\equiv w\pmod{2}$
  admits a flat connection.
\item In particular, if $w\in H^2(X;\ZZ)$ and $\tilde X =
  X\#\overline{\CC\PP}^2$ is the blow-up, with exceptional class
  $e=[\overline{\CC\PP}^1]\in H_2(X;\ZZ)$ and Poincar\'e dual $e^*\in
  H^2(X;\ZZ)$, then no $\SO(3)$ bundle $F$ over $\tilde X$ with
  $w_2(F)=w+e^*\pmod{2}$admits a flat connection.
\end{itemize}
\end{lem}

Therefore, if $X$ has a spherical class $e\in H_2(X;\ZZ)$ such that
$\langle w,e\rangle\neq 0\pmod{2}$, then Lemma
\ref{lem:MorganMrowka} implies that there are no
pairs in $\sM_{\ft}$ with flat $\SO(3)$ connections when $w_2(\ft)\equiv
w\pmod{2}$. In this situation, Proposition
\ref{prop:ClassificationOfStabilizers} implies that the only fixed points
of the circle action on $\sM_{\ft}$ are given by either points $[A,\Phi]$ with
$A$ irreducible and $\Phi\equiv 0$ or $A$ reducible and $\Phi\not\equiv 0$.

Suppose $w\in H^2(X;\ZZ)$ is any integral class and that
$\ft=(\rho,V)$ is a \spinu structure over $\tilde X$ with
$w_2(\ft) \equiv w+e^*\pmod{2}$. Then the $\SO(3)$ bundle
$\fg_{\ft}$ over $\tilde X$ obeys the Morgan-Mrowka criterion of
Lemma \ref{lem:MorganMrowka}. Because the Donaldson and
Seiberg-Witten invariants of $\tilde X$ determine and are
determined by those of $X$, (see \cite{FSTurkish},
\cite{FSBlowUp}, \cite{MorganMrowkaPoly}), no information about
these invariants is lost by passing to the blow-up in this way.

In Lemma \ref{lem:RedAreSW} we show that the moduli space of
Seiberg-Witten monopoles, $M_{\fs}^{0}$, can be identified with
the subspace of $\PU(2)$ monopoles in $\sM_{\ft}^0$ which are
reducible with respect to the splitting $\ft=\fs\oplus\fs\otimes
L$, while zero-section points in $M_{\fs}$ are mapped to
zero-section reducibles in $\sM_{\ft}$ (though not necessarily
injectively). If $w_2(\ft)$ obeys the Morgan-Mrowka criterion
then Proposition \ref{prop:ClassificationOfStabilizers} implies
that reducible $\PU(2)$ monopoles cannot be zero-section pairs,
leading to the following

\begin{cor}
\label{cor:NoSWZeroSections}
If $M_{\fs}\hookrightarrow \sM_{\ft}$ is a Seiberg-Witten moduli subspace
and $w_2(\ft)$ obeys the Morgan-Mrowka criterion of Lemma
\ref{lem:MorganMrowka}, then $M_{\fs}$ contains no zero-section
solutions.
\end{cor}

\subsection{The link of the stratum of anti-self-dual PU(2) monopoles}
\label{subsec:ASDLink}
In this section we describe the Kuranishi model for a neighborhood of a
zero-section solution in the $\PU(2)$ monopole moduli space and define a
link of the stratum $M^w_\ka$ in $\sM_{\ft}$.  The differential $D\fS$ will
not be surjective at all the zero-section monopoles $(A,0)$, but we will
show that the cokernel of $(D\fS)_{A,0}$, namely $H^2_{A,0}$, can be
identified with the cokernel of the perturbed Dirac operator,
$D_{A}+\rho(\vartheta)$.

If we are given a decomposition, $V=W\otimes E$, Lemma
\ref{lem:DecompSpinuGaugeTransformation} provides a canonical
identification of automorphism groups, $\sG_\kappa^w\cong\sG_{\ft}$, and with
respect to this identification and choice of fixed connection $A_\Lambda$ on
$\det^{\frac{1}{2}}(V^+)$, Lemma
\ref{lem:AdjointConnection} provides a canonical gauge-equivariant isomorphism
$$
\sA_{\kappa}^w \cong \sA_{\ft}, \quad \hat A \mapsto A,
$$
where the space $\sA_{\kappa}^w$ of $\SO(3)$ connections on
$\su(E)\cong\fg_{\ft}$ was defined in \S
\ref{subsubsec:SpaceSO(3)Connections}, with $p_1(\fg_{\ft})=-4\kappa$ and
$w_2(\fg_{\ft})\equiv w\pmod{2}$.

Let $\sB_{\ft} = \sA_{\ft}/\sG_{\ft}$ be the quotient of the space of $L^2_k$
spin connections on $V$. By the preceding discussion there is a canonical
homeomorphism,
$$
\sB_{\kappa}^w \cong \sB_{\ft}, \quad [\hat A] \mapsto [A],
$$
restricting to diffeomorphisms on smooth strata.  There are a canonical
smooth embedding $\sA_{\ft}\hookrightarrow \tsC_{\ft}$ and a ``smoothly
stratified embedding'' $\sB_{\ft}\hookrightarrow \sC_{\ft}$---that is, a
topological embedding restricting to smooth embeddings on smooth
strata---given by $A\mapsto (A,0)$ and $[A]\mapsto [A,0]$. Combining these
identifications and embeddings, we obtain a gauge-equivariant
smooth embedding and a smoothly stratified embedding,
\begin{equation}
\label{eq:EmbedSO(3)ConnsIntoPairs}
\iota:\sA_{\kappa}^w \embed \tsC_{\ft}, \quad \hat A \mapsto (A,0)
\quad\text{and}\quad
\iota:\sB_{\kappa}^w \embed \sC_{\ft}, \quad [\hat A] \mapsto [A,0].
\end{equation}
Hence, we can identify the image of the induced smoothly stratified embedding
\begin{equation}
\label{eq:EmbedASDmodSpaceIntoPairs}
\iota:M_{\kappa}^w(X) \embed \sC_{\ft}, \quad [\hat A] \mapsto [A,0].
\end{equation}
with the subspace of $\sM_{\ft}$ given by the zero-section solutions to the
$\PU(2)$ monopole equations \eqref{eq:PT}.  Therefore, we shall refer to
pairs or connections representing points in $M^w_\kappa$ as {\em
zero-section $\PU(2)$ monopoles\/} or anti-self-dual connections, depending
on the context.

For a generic Riemannian metric on $X$, the
moduli space $M^{w,*}_\ka$ of irreducible anti-self-dual connections is
a smooth manifold \cite[Corollary 4.3.18]{DK}, \cite{FU}
of the expected dimension $d_a(\ft)$ given in equation
\eqref{eq:Defndana}.

At a zero-section $\PU(2)$ monopole $(A,0)$, the elliptic deformation
complex \eqref{eq:ConstDetDefComplex} splits into a direct sum of
complexes:
\begin{equation}
\label{eq:ASDDeformationComplex}
\begin{gathered}
\begin{CD}
C^\8(\fg_{\ft})
@> d_{\hat A}>>
C^\8(\Lambda^1\otimes \fg_{\ft})
@> d^+_{\hat A} >>
C^\8(\Lambda^+\otimes \fg_{\ft}),
\end{CD}
\\
\begin{CD}
C^\8(V^+)
@> D_{A}+\rho(\vartheta) >>
C^\8(V^-).
\end{CD}
\end{gathered}
\end{equation}
The first complex in \eqref{eq:ASDDeformationComplex} is simply the
elliptic deformation complex for the moduli space $M^w_{\ka}$ of
anti-self-dual connections, with cohomology $H_{\hat A}^\bullet$. The
following lemma is then a clear consequence of the preceding decomposition:

\begin{lem}
\label{lem:DecompOfZeroSectionDeformation}
If $(A,0)$ is a zero-section $\PU(2)$ monopole then there are
canonical isomorphisms,
\begin{align*}
H^0_{A,0} &\cong H^0_{\hat A}, \\
H^1_{A,0} &\cong  H^1_{\hat A}\oplus \Ker(D_{A}+\rho(\vartheta)), \\
H^2_{A,0} &\cong H^2_{\hat A}\oplus \Coker(D_{A}+\rho(\vartheta)).
\end{align*}
\end{lem}

If the manifold $X$ is not simply connected, there exist anti-self-dual
connections which, while not globally reducible, become reducible when
restricted to certain open sets in $X$: this complicates the construction
of the restriction maps used in the definition of geometric
representatives on moduli spaces of anti-self-dual connections in
\cite{KMStructure} and for $\PU(2)$ monopoles here.  Therefore, we recall
some of the technical points from \cite[pp. 586--588]{KMStructure} which we
shall need to address these additional complications when $X$ is not simply
connected.  A connection $\hat A$ on an $\SO(3)$ bundle $F$ over $X$
is called a {\em twisted
reducible\/} if it preserves a splitting $F=\lambda\oplus
N$, where $\lambda$ is a non-trivial real line bundle and $N$ is an
$\Or(2)$ bundle. The curvature $F_{\hat A}$ of a twisted reducible
connection $\hat A$
has rank one even though  $\hat A$ is not reducible. When restricted to open
sets in $X$ over which the real line bundle $\lambda$ is trivial, the
connection $\hat A$ becomes reducible. Let $H^i(X;\lambda)$ denote the
cohomology group of dimension $b^i(\lambda)$ with coefficients in the local
system given by $\lambda$. The bundle $\lambda\otimes\lambda$ is trivial
and so the cup product gives a pairing $H^2(X;\lambda)\otimes
H^2(X;\lambda)\rightarrow H^4(X;\RR)\cong\RR$.  We define $b^+(\lambda)$
and $b^-(\lambda)$ to be the maximum dimensions of subspaces of
$H^2(X;\lambda)$ on which this pairing is respectively positive or negative
definite.  If $F$ admits a reduction $\lambda\oplus N$,
then the corresponding subspace of twisted reducibles in $M_\kappa^w$
(where $p_1(F)=-4\kappa$ and $w_2(F)\equiv w\pmod{2}$)
defined by this splitting is a torus of dimension $b^1(\lambda)$.  Whether
or not such connections exist in 
$M^w_{\ka}$ for a generic metric
depends on $b^+(\lambda)$, as described in the following lemma.

\begin{lem}\cite[Lemma 2.4 \& Corollary 2.5]{KMStructure}
\label{lem:TwistedReducibles}
Let $X$ be a smooth four-manifold, $\lambda$ a non-trivial real line bundle
on $X$, and $F$ an oriented, orthogonal three-plane bundle for which $F$
admits a reduction $F=\lambda\oplus N$. Suppose $p_1(F)=-4\kappa\neq 0$ and
$w_2(F)\equiv w\pmod{2}$, where $w\in H^2(X;\ZZ)$.
\begin{enumerate}
\item
If $b^+(\lambda)=0$, then $b^1(\lambda)=-(b^+-b^1+1)$ and $M^w_\kappa$
contains twisted reducibles corresponding to $F=\lambda\oplus N$ for
all Riemannian metrics on $X$.  If, in addition, $b^1(\lambda)\ge 0$, then
$-2p_1(F)-3(1-b^1+b^+)> b^1(\lambda)$ and if the Riemannian metric
on $X$ is generic, each anti-self-dual connection $\hat A$,
which is a twisted reducible with respect to 
$F=\lambda\oplus N$, has $H^2_{\hat A}=0$.
\item
If $b^+(\lambda)\geq 1$, then for generic Riemannian metrics on $X$,
the space $M^w_\kappa$ contains no twisted reducibles for the splitting
$F=\lambda\oplus N$. 
\end{enumerate} 
\end{lem}

Therefore, when $b_2^+(X)>0$, the only twisted reducibles appearing
in $M^w_{\ka}$ have codimension at least
$-2p_1(\fg_{\ft})-2(b_2^+-b_1+1)$ and are smooth points.  The
Kuranishi lemma, \cite[Proposition 4.2.19]{DK} or \cite[Lemma 4.7]{FU},
gives the following description of a neighborhood of a
zero-section monopole in $\sM_{\ft}$.

\begin{cor}
\label{cor:ASDKuranishi}
Let $\ft$ be a \spinu structure on a closed, oriented, smooth four-manifold
$X$ with $b_2^+(X)>0$ and generic Riemannian metric.  Let $[A,0]$ be a point
in the image of $M^w_\kappa\embed \sM_{\ft}$, so $H_{\hat A}^0 = 0 = H^2_{\hat
A}$ for $[\hat A]\in M^w_\kappa$.  Then there are
\begin{itemize}
\item
An open, $S^1$-invariant neighborhood $\sO_A$ of the origin in
$T_{\hat A}M^w_{\ka}\oplus\Ker (D_A+\rho(\vartheta))$ together with a
smooth, $S^1$-equivariant embedding
$$
\bga_A:\sO_A\subset T_{\hat A}M^w_{\ka}\oplus\Ker (D_A+\rho(\vartheta))
\hookrightarrow
\tsC_{\ft},
$$
with $\bga_A(0,0) = (A,0)$ and $\sM_{\ft}\cap\bga_A(\sO_A)$ an open
neighborhood of $[A,0]$ in $\sM_{\ft}$, and
\item A smooth, $S^1$-equivariant map
$$
\bvarphi_A:\sO_A\subset T_{\hat A}M^w_{\ka}\oplus\Ker (D_A+\rho(\vartheta))
\to
\Coker (D_A+\rho(\vartheta))
$$
such that $\bga_A$ restricts to an $S^1$-equivariant,
smoothly-stratified diffeomorphism from
$\bvarphi_A^{-1}(0)\cap\sO_A$ onto $\sM_{\ft}\cap\bga_A(\sO_A)$.
\end{itemize}
\end{cor}

If $n_a\leq 0$, then at a generic point $[A,0]\in \sM_{\ft}$ where
$\Ker (D_A+\rho(\vartheta))=\{0\}$,
(assuming the map from $M^w_\kappa$ to the space of Fredholm operators
of index $n_a$ is transverse to the ``jumping lines strata'' as
described in \cite{Koschorke}) the Kuranishi model in Lemma
\ref{cor:ASDKuranishi} shows that a neighborhood of
$[A,0]\in \sM_{\ft}$ contains no elements of $M^0_{\ft}$.
Thus, such points in the anti-self-dual moduli space are
isolated from the subspace $\sM_{\ft}^0$ of non-zero-section points.
For this reason, we will restrict our attention to the
cases where $n_a>0$.

Although we can only describe a neighborhood of the anti-self-dual
connections locally, because of the problem of spectral flow, we can still
introduce a global, codimension-one subspace of the compactification
$\bar\sM_{\ft}$ which will serve as a link.  This space might not have a
fundamental class because it is not known to have locally finite topology
near the lower strata of the Uhlenbeck compactification $\bar\sM_{\ft}$.
However, we shall see that the local Kuranishi model
in Corollary \ref{cor:ASDKuranishi} will suffice to define
intersections under some additional assumptions.

\begin{defn}
\label{defn:ASDLink}
The {\em link of $\barM^w_{\kappa}$\/} in $\bar\sM_{\ft}$ is given by
$$
\bL^{w,\varepsilon}_{\ft,\kappa}
=
\{[A,\Phi,\bx]\in \bar\sM_{\ft}/S^1: \|\Phi\|_{L^2}^2=\varepsilon\}.
$$
We write $\bL^{w}_{\ft,\kappa}$ when the positive constant $\varepsilon$ is
understood.
\end{defn}

\begin{lem}
\label{lem:ContinuousL2Norm}
For generic $\varepsilon>0$, the link $\bL^{w,\varepsilon}_{\ft,\kappa}$ is
closed under the $S^1$ action, is a
smoothly stratified, closed subspace of $\bar\sM_{\ft}$, and has codimension
one in every stratum of $\bar\sM_{\ft}$ which it intersects.
\end{lem}

\begin{proof}
That $\bL^{w,\varepsilon}_{\ft,\kappa}$ is closed under the $S^1$ action
follows directly from the definitions.

It is clear that the function
$$
\bl:\bar\sM_{\ft} \to \RR, \qquad [A,\Phi,\bx]\mapsto \|\Phi\|_{L^2}^2,
$$
is smooth on each stratum of $\bar\sM_{\ft}$: we claim it is continuous
on $\bar\sM_{\ft}$.  Let $[A_\alpha,\Phi_\alpha]$ be a sequence of points in
$\bar\sM_{\ft}$ which converge to $[A_\infty,\Phi_\infty,\bx]$.  We may
assume, by an appropriate choice of a sequence of $L^2_{k+1}$ \spinu gauge
transformations of $V$, that the sequence $\{\Phi_\alpha\}$ converges to
$\Phi_\infty$ in $C^\infty$ on $X - B(\bx,r)$, where we define $B(\bx,r) =
\cup_{x\in\bx}B(x,r)$. Therefore,
\begin{align*}
\left| \|\Phi_\alpha\|_{L^2(X)}^2  - \|\Phi_\infty\|_{L^2(X)}^2\right|
&\leq
\|\Phi_\alpha - \Phi_\infty\|_{L^2(X-B(\bx,r))}^2
+ \|\Phi_\alpha - \Phi_\infty\|_{L^2(B(\bx,r))}^2
\\
&\leq
\|\Phi_\alpha - \Phi_\infty\|_{L^2(X-B(\bx,r))}^2 + Cr^4,
\end{align*}
where the second inequality follows from the universal {\em a priori\/}
$C^0$ bound for the sequence $\Phi_\alpha$ and $\Phi_\infty$ given by
\cite[Lemma 4.4]{FL1}. Thus,
$$
\limsup_{\alpha\to\8}
\left| \|\Phi_\alpha\|_{L^2(X)}^2  - \|\Phi_\infty\|_{L^2(X)}^2\right|
\le Cr^4, \quad\text{for all $r>0$},
$$
and so
$$
\lim_{\alpha\to\8}\|\Phi_\alpha\|_{L^2(X)}^2 = \|\Phi_\infty\|_{L^2(X)}^2,
$$
as desired. A generic $\varepsilon>0$ is a regular value for the function
$\bl$ on each smooth stratum of $\bar\sM_{\ft}$. For such an $\eps$, the
preimage $\bl^{-1}(\eps)$ is a smooth submanifold of each stratum and
because the function $\bl:\bar\sM_{\ft} \to \RR$ is continuous, these smooth
submanifolds fit together to form a smoothly stratified subspace of
$\bar\sM_{\ft}$ (see Remark 3.3 in \cite{FL2b}).
\end{proof}

\begin{rmk}
In defining a link, it might seem more natural to work with the image of
the $\varepsilon$-sphere in the normal bundle given by $\Ker
(D_A+\rho(\vartheta))$,  
at least on the image of $M^{w}_{\kappa}\embed \sM_{\ft}$ where the
cokernel of the Dirac operator vanishes.  This definition would have the
disadvantage of not being a global object because of the jumping-line
problem (that is, spectral flow).  However, it can be shown that the two
functions defined on the open set $\sO$ of
the Kuranishi model of $[A,0]$ in
Corollary \ref{cor:ASDKuranishi},
one given by the $L^2$ norm of the element of $\Ker (D_{A}+\rho(\vartheta))$,
the other defined by $\bl\circ\bga$
are $C^1$ close as $\varepsilon$ goes to zero, as they differ by the
difference of $\bga$ and the identity.
As we shall see in \S 3.4.2 in \cite{FL2b} it is sufficient that
these two links are cobordant.
\end{rmk}

\subsection{Reducible PU(2) monopoles and their identification with
Seiberg-Witten monopoles}
\label{subsec:ReducibleSingularities}
In this section we show that the subspaces of reducible $\PU(2)$ monopoles
in $\sM_{\ft}$ can be identified with the moduli spaces of Seiberg-Witten
monopoles defined in \S \ref{subsec:PertSW}.

\subsubsection{Decomposing $\PU(2)$ bundles}
We now describe the canonical isomorphisms of associated bundles induced
from a splitting $E=L_1\oplus L_2$ of a Hermitian two-plane bundle $E$; the
lemmas below are of course elementary and we just state them in order to
make our conventions clear.

\begin{lem}
\label{lem:su(E)Split}
If $E$ is a Hermitian two-plane bundle and $L_1$, $L_2$ are Hermitian line
bundles over a manifold $X$ such that $E=L_1\oplus L_2$, then there is a
canonical isometry of $\SO(3)$ bundles, $\su(E)\cong i\underline{\RR}\oplus
(L_2\otimes L_1^*)$, given by
\begin{equation}
\label{eq:su(E)Split}
\begin{pmatrix}\nu &-\bar z \\  z & -\nu\end{pmatrix}\longmapsto (\nu, z).
\end{equation}
\end{lem}

\begin{proof}
With respect to the decomposition $E=L_1\oplus L_2$, any element
$M\in\gl(E)$ takes the form
$$
M = \begin{pmatrix}M_{11} & M_{12} \\ M_{21} & M_{22}\end{pmatrix},
$$
where $M_{jk}\in L_j\otimes L_k^*$, for $1\le j,k\le 2$. Thus,
$M_{11},M_{22}\in \CC$, after identifying $L_1\otimes L_1^*\cong
L_2\otimes L_2^* \cong \CC$. If $M\in\su(E)$, then $M_{11}\in i\RR$ and
$M_{12}\in L_1\otimes L_2^*$ with $M_{22} = - M_{11}$ and $M_{21} = -\bar
M_{12}$, so any element of $\su(E)$ takes the shape
$$
M = \begin{pmatrix}\nu & -\bar z \\  z & -\nu\end{pmatrix},
\qquad \nu\in\CC, \quad  z\in L_2\otimes L_1^*,
$$
and the desired isomorphism $\su(E)\cong i\ubarRR\oplus
(L_2\otimes L_1^*)$ is given by $M\mapsto (\nu, z)$.

We recall from \cite[\S 2.4]{FL1} that the induced fiber inner
product on $\gl(E)$ is defined by $\langle M,M'\rangle =
\half\tr(M' M^\dagger)$. Thus, if $M',M\in C^\8(\su(E))$
correspond to $(\nu', z'),(\nu,
z)\in C^\8(i\Lambda^0))\oplus C^\8(L)$, respectively, we see
that
\begin{equation}
\label{eq:CompareInnerProduct}
\begin{aligned}
\langle M',M\rangle &= \thalf\tr(M' M^\dagger)
= \thalf\tr\begin{pmatrix}\nu' & -\bar z' \\  z' & -\nu'\end{pmatrix}
\begin{pmatrix}\bar\nu & -\bar z \\  z & -\bar\nu\end{pmatrix}
\\
&= \nu'\bar\nu + \thalf( z'\bar z + \bar z' z)
= \nu'\bar\nu + \Real\langle z',z\rangle_\CC, 
\end{aligned}
\end{equation}
and so the isomorphism of $\SO(3)$ bundles is an isometry.
\end{proof}

\subsubsection{The identification of reducible $\PU(2)$ monopoles with
Seiberg-Witten monopoles}
Next, we explain how reducible pairs in the moduli space of $\PU(2)$
monopoles may be identified with Seiberg-Witten monopoles.

Recall that $\sC_{\ft}$ is the quotient space of pairs whose
\spinu connections induce the fixed unitary connection $2A_\La$ on
the complex line bundle $\det(V^+)$.

\begin{lem}
\label{lem:TopU1Embedding}
Let $\fs=(\rho,W)$ be a \spinc structure over an oriented, Riemannian
four-manifold $X$ and let $\ft=(\rho,V)$ be a \spinu structure with
$V=W\oplus W\otimes L$. Then
\begin{equation}
\label{eq:DefnOfIota}
\iota:\tsC_{\fs}\hookrightarrow \tsC_{\ft},
\qquad
(B,\Psi)
\mapsto
(B\oplus B\otimes A_L,\Psi\oplus 0),
\end{equation}
is a smooth embedding, where $A_L=A_\La\otimes (B^{\det})^*$, and is
gauge-equivariant with respect to
\begin{equation}
\label{defn:GInclusion}
\varrho:\sG_{\fs}\hookrightarrow\sG_{\ft},
\qquad
s\mapsto
s\,\id_{W} \oplus s^{-1}\,\id_{W\otimes L},
\end{equation}
so that
$$
\iota(s(B,\Psi)) = \varrho(s)\iota(B,\Psi).
$$
The image of the map \eqref{eq:DefnOfIota} contains all pairs in $\tsC_{\ft}$
fixed by the action \eqref{eq:VReducibleAction} of $S^1$ on $V$. The
induced map,
\begin{equation}
\label{eq:DefnOfIotaOnQuotient}
\iota:\sC_{\fs}\hookrightarrow \sC_{\ft},
\qquad
[B,\Psi]
\mapsto
[B\oplus B\otimes A_L,\Psi\oplus 0],
\end{equation}
is continuous and, when restricted to $\sC_{\fs}^0$, a topological embedding.
If $w_2(\ft)\neq 0$, then \eqref{eq:DefnOfIotaOnQuotient} is a topological
embedding of $\sC_{\fs}$. If $w_2(\ft)=0$ then the map
\eqref{eq:DefnOfIotaOnQuotient} takes zero-section points in $\sC_{\fs}$ to
zero-section reducibles in $\sC_{\ft}$, although this identification of
zero-section points need not be injective if $b_1(X)> 0$.
\end{lem}

\begin{proof}
  The map $\iota:\tsC_{\fs}\rightarrow \tsC_{\ft}$ is clearly a $C^\infty$
  embedding. Furthermore,
$$
\varrho(s)\iota(B,\Psi) = (s_*B\oplus (s^{-1})_*(B\otimes A_L),s\Psi\oplus 0).
$$
Since $A_L = A_\Lambda\otimes (B^{\det})^*$, we see that $s\in \sG_{\fs}$
acts on $A_L$ as $(s^{-2})_*A_L$ (see also \S
\ref{subsubsec:DecompGroupActions}) and so $(s^{-1})_*(B\otimes A_L) =
s_*B\otimes (s^{-2})_*A_L$. Thus,
$\varrho(s)\iota(B,\Psi)=\iota(s_*B,s\Psi)$, as desired.

Next, we characterize the image of the map $\iota:\tsC_{\fs}\to\tsC_{\ft}$.
Suppose $(A,\Phi)\in \tsC_{\ft}$ is fixed by the $S^1$ action
\eqref{eq:VReducibleAction} on $V$: this action descends to the action
\eqref{eq:SO(3)ReducibleAction} on $\fg_{\ft} = 
i\underline{\RR}\oplus L$, which fixes the induced connection $\hat A$.
As in the proof of Assertion (2) of Proposition
\ref{prop:ClassificationOfStabilizers}, the connection $\hat A$ must then be
reducible with respect to this splitting, taking the form $\hat A =
d_\RR\oplus A_L$ on $\fg_{\ft}$ and $A=B\oplus B\otimes A_L$ on $V$.
Thus, $(A,\Phi)$ is in the image of $\iota$.

If $\sU\subset\sC_{\fs}$ is an open subset, so it is easy to see that
$\iota(\sU)$ is open in $\iota(\sC_{\fs})$ with respect to the subspace
topology induced by $\sC_{\ft}$.

We now show $\iota:\sC^{0}_{\fs}\to\sC_{\ft}$ is
injective and that
$\iota:\sC_{\fs}\to\sC_{\ft}$ is injective when $w_2(\ft)\not\equiv 0$.
Suppose $(B,\Psi)$ and $(B',\Psi')$ are pairs in
$\tsC^{0}_{\fs}$ and that $u\in\sG_{\ft}$ satisfies
$u\iota(B,\Psi)=\iota(B',\Psi')$, and thus
$[\iota(B,\Psi)]=[\iota(B',\Psi')]$ in $\sC_{\ft}$. With respect
to the decomposition $V=W\otimes E$ and splitting
$E=\underline{\CC}\oplus L$, we can consider $u$ to be a gauge
transformation of $E=\underline{\CC}\oplus L$ via the isomorphism
$\sG_\kappa^w\cong\sG_{\ft}$ in Lemma
\ref{lem:DecompSpinuGaugeTransformation} (and acting as the
identity on $W$) and write
$$
u
=
\begin{pmatrix} s & v \\ -\bar v & \bar s \end{pmatrix},
\quad \text{where $s\in\Map(X,\CC)$ and
$v\in C^\8(\Hom(L,\CC))$},
$$
noting that $u^{-1}=u^\dagger$ and $\det(u)=1$.  As in the proof of
Lemma \ref{lem:ReducibleSpinu}, we may write $\iota(B) = B\otimes A_E$
and $\iota(B') = B\otimes A_E'$, where $A_E$ and $A_E'=u(A_E)$ are
$\U(2)$ connections on $E$ which are reducible with respect to the splitting
$E=\underline{\CC}\oplus L$. Let $\varrho_L:S^1\to\End(E)$ denote the
action of $S^1$ on $E$ by the trivial action on $\underline{\CC}$ and
scalar multiplication on $L$. Then $\varrho_L(S^1)$ fixes $u(A_E)$
and thus for all real $\theta$
$$
u^{-1}\varrho_L(e^{i\theta})u \in \Stab_{\U(2)}(A_E) = S^1\times S^1,
$$
where $\Stab_{\U(2)}(A_E)$ is the stabilizer of $A_E$ in the space of
unitary automorphisms of $E$. Consequently, for every real
$\theta$ there are real constants $\mu$ and $\nu$ such that
$$
\begin{pmatrix} 1      & 0      \\ 0       & e^{i\theta} \end{pmatrix}
\begin{pmatrix} s      & v      \\ -\bar v & \bar s      \end{pmatrix}
=
\begin{pmatrix} s        & v      \\ -\bar v & \bar s      \end{pmatrix}
\begin{pmatrix} e^{i\mu} & 0      \\ 0       & e^{i\nu}    \end{pmatrix},
$$
and we can assume $e^{i\theta}\neq 1$. Simplifying, this becomes
$$
\begin{pmatrix} s      & v \\ -e^{i\theta}\bar v & e^{i\theta}\bar s      
\end{pmatrix}
=
\begin{pmatrix} e^{i\mu}s & e^{i\nu}v \\ -e^{i\mu}\bar v 
& e^{i\nu}\bar s      
\end{pmatrix}.
$$
If $s\not\equiv 0$, then we must have $e^{i\mu}=1$ and thus $v=0$
because $e^{i\theta}\bar v = \bar v$.  Since $\det(u)=|s|^2+|v|^2=1$, we
see that $s\in\Map(X,S^1)\cong\sG_{\fs}$. Hence, $u=\varrho(s)$ and
$[B,\Psi]=[B',\Psi'] \in \sC_{\fs}$ because
$s(B,\Psi)=(B',\Psi')$. 

It remains to consider the case $s\equiv 0$, for which we must
then have $v\bar v=1$. First suppose $\Psi\not\equiv 0$ and
observe that
$$
u(\Psi\oplus 0) = s\Psi\oplus -\bar v\Psi = \Psi'\oplus 0,
$$
so $\bar v\Psi \equiv 0$ on $X$ and hence $\Psi\equiv 0$,
contradicting our hypothesis that $[B,\Psi]\in\sC_{\fs}^0$.
Therefore the case $s\equiv 0$ cannot occur in this situation.
Thus, $\iota:\sC^0_{\fs}\to\sC_{\ft}$ is injective.

Otherwise, if $s\equiv 0$, suppose $w_2(\ft)\neq 0$ and observe that the
automorphism $u\in\Aut(\underline{\CC}\oplus L)$ induces an isomorphism
$L\cong\underline{\CC}$. But $c_1(E)=c_1(L)$ and $c_1(E)\equiv
w_2(\ft)\pmod{2}$: by the hypothesis in the final statement of the lemma,
$w_2(\ft)\neq 0$, so this contradicts $L\cong \underline{\CC}$ and
therefore the case $s\equiv 0$ cannot occur in this situation either.
Thus, if $w_2(\ft)\neq 0$, the map $\iota:\sC_{\fs}\to\sC_{\ft}$ is injective.
\end{proof}

The following lemma highlights the key property of the $\PU(2)$ monopole
equations: 

\begin{lem}
\label{lem:RestrictionOfPU2MonopoleEquation}
Let $(\rho,V)$ be a \spinu structure on $X$.
Suppose $(A,\Phi) = (B\oplus B\otimes A_L,\Psi\oplus 0)$ is a reducible
pair with respect to a splitting $V=W\oplus W\otimes L$, where $B$ is a
spin connection on $W$, $A_L$ is a unitary connection on a complex line
bundle $L$, and $\Psi$ is a section of $W^+$. Then $(A,\Phi)$ solves
the $\PU(2)$ monopole equations \eqref{eq:PT} if and only if $(B,\Psi)$
solves the Seiberg-Witten equations
\eqref{eq:U1Monopole} with perturbation $\eta=F_{A_\Lambda}^+$.
\end{lem}

\begin{proof}
Suppose that $(A,\Phi)$ solves the 
$\PU(2)$ monopole equations \eqref{eq:PT}.
With respect to the splitting $\fg_{\ft}=i\underline{\RR}\oplus L$, Lemma
\ref{lem:ReducibleSpinu} implies that $\hat A=d_{\RR}\oplus A_L$, where
$d_{\RR}$ is the product connection on $i\underline{\RR}$ and
$A_L=A_\Lambda\otimes (B^{\det})^*$ is a
unitary connection on the complex line bundle
$L=\det^{\frac{1}{2}}(V^+)\otimes\det(W^+)^*$. Let
$E=\underline{\CC}\oplus L$ and observe that $\tilde A = d_\CC\oplus A_L$
is a reducible unitary connection on $E$ which is a lift of $\hat A$ on
$\su(E)$. Then
$$
F_{\tilde A}^+ 
=
\begin{pmatrix}
0 & 0 \\ 0 & F_{A_L}^+
\end{pmatrix}
=
F_{A_L}^+
\otimes
\begin{pmatrix}
0 & 0 \\ 0 & 1
\end{pmatrix}
\in
C^\8(\Lambda^+\otimes \fu(E)),
$$
and therefore, since $F_{\hat A}^+ \in C^\8(\Lambda^+\otimes\so(\su(E)))$,
we have 
\begin{equation}
\label{eq:ReducibleCurvatureDecomposition}
\ad^{-1}(F_{\hat A}^+)
=
(F_{\tilde A}^+)_0
=
-F_{A_L}^+
\otimes
\begin{pmatrix}
\half  & 0 \\ 0 & -\half
\end{pmatrix}
\in
C^\8(\Lambda^+\otimes \su(E)).
\end{equation}
If $\Phi=\Psi\oplus 0$, where $\Psi\in\Omega^0(W^+)$ and writing
$\fu(V^+)=\fu(W^+)\otimes \fu(E)$, we have
$$
\Phi\otimes\Phi^*
=
\begin{pmatrix}
\Psi\otimes\Psi^* & 0 \\
0 & 0
\end{pmatrix}
=
(\Psi\otimes\Psi^*)
\otimes
\begin{pmatrix}
1 & 0 \\
0 & 0
\end{pmatrix}
\in
C^\8(\fu(W^+)\otimes \fu(E)).
$$
Hence, projecting to $\su(W^+)\otimes \su(E)$, we see that
\begin{equation}
\label{eq:QuadraticTermReducibleDecomposition}
(\Phi\otimes\Phi^*)_{00}
=
(\Psi\otimes\Psi^*)_0
\otimes
\begin{pmatrix}
\half & 0 \\
0     & -\half
\end{pmatrix}
\in
C^\8(\su(W^+)\otimes \su(E)).
\end{equation}
Furthermore, with respect to the splitting $V^+=W^+\oplus W^+\otimes L$,
we clearly have 
\begin{equation}
\label{eq:DiracSplit}
(D_{A}+\rho(\vartheta))\Phi = ((D_{B}+\rho(\vartheta))\Psi,0). 
\end{equation}
Then, combining equations \eqref{eq:PT},
\eqref{eq:ReducibleCurvatureDecomposition}, 
\eqref{eq:QuadraticTermReducibleDecomposition}, and \eqref{eq:DiracSplit},
and noting that
$$
F_{A_L} 
= 
F_{A_\Lambda} - F_{B^{\det}}
=
F_{A_\Lambda} - \Tr(F_{B}),
$$
shows that the pair $(B,\Psi)$ solves
\begin{equation}
\label{eq:PU2MonopoleComposition}
\begin{aligned}
\Tr(F_{B}^+) - \tau\rho^{-1}(\Psi\otimes\Psi^*)_0
- F_{A_\Lambda}^+
&=0,
\\
(D_{B}+\rho(\vartheta))\Psi
&=0.
\end{aligned}
\end{equation}
Comparing \eqref{eq:PU2MonopoleComposition} with the Seiberg-Witten
equations \eqref{eq:U1Monopole} concludes the proof.
\end{proof}

\begin{lem}
\label{lem:RedAreSW}
Let $X$ be a closed, oriented, Riemannian four-manifold with \spinu
structure $\ft$ having a splitting $\ft=\fs\oplus\fs\otimes L$.  If
$w_2(\ft)\neq 0$ then the map \eqref{eq:DefnOfIotaOnQuotient} of $\sC_{\fs}$
into $\sC_{\ft}$ restricts to a topological embedding
\begin{equation}
  \label{eq:RedAreSW}
\iota:M_{\fs}
\embed
\sC_{\ft},
\end{equation}
whose image is the subspace of $\sM_{\ft}$ represented by pairs $(A,\Phi)$
which are reducible with respect to the splitting $V=W\oplus W\otimes L$,
with $\Phi\in L^2_k(W^+)$ and $\ft=(\rho,V)$ and $\fs=(\rho,W)$.  If
$w_2(\ft)=0$ then the map \eqref{eq:DefnOfIotaOnQuotient} takes
zero-section points in $M_{\fs}$ to zero-section reducibles in $\sM_{\ft}$,
although this identification of zero-section points need not be injective
when $b_1(X)> 0$.
\end{lem}

\begin{proof}
Given Lemma \ref{lem:TopU1Embedding}, we need only characterize the image
of $\iota$. If $(A,\Phi)$ represents a point in $\sM_{\ft}$ and is reducible
with respect to the splitting $V=W\oplus W\otimes L$, with
$\Phi\in L^2_k(W^+)$, then Lemma
\ref{lem:TopU1Embedding} implies that $(A,\Phi) = \iota(B,\Psi)$, for
some pair $(B,\Psi)\in\tsC_{\fs}$.  Then Lemma
\ref{lem:RestrictionOfPU2MonopoleEquation} implies that $(B,\Psi)$
satisfies the Seiberg-Witten equations
\eqref{eq:U1Monopole} since $(A,\Phi)$ satisfies the $\PU(2)$ monopole
equations \eqref{eq:PT}.
\end{proof}

\subsection{The link of a stratum of reducible monopoles: local structure}
\label{subsec:LinkOfReduciblesLocal}
The construction of the link in $\sM_{\ft}$ of the stratum
$M_{\fs}\embed \sM_{\ft}$ of reducible $\PU(2)$ monopoles occupies this
and the next subsection.  Although for generic perturbations the locus
$M_{\fs}$ of reducible solutions in $\sM_{\ft}$ defined by a
reduction $\ft=\fs\oplus \fs\otimes L$ will be a smooth manifold, the
linearization of the map $\fS$ defined by the $\PU(2)$ monopole equations
\eqref{eq:PT} need not be surjective at a reducible
solution and so the points of $M_{\fs}$ might not be regular
points of $\sM_{\ft}$.  Moreover, unlike the case of the link of $M^w_{\ka}$,
a local model of the link does not suffice as only one of our cohomology
classes extend over $M_{\fs}$ (and that one vanishes in many
cases), so we cannot use geometric representatives to cut down to a generic
point in $M_{\fs}$.

We begin with a definition of link of a stratum in smoothly stratified
space, essentially following Mather \cite{Mather} and Goresky-MacPherson
\cite{GorMacPh}. We need only consider the relatively simple case of a
stratified space with two strata since the lower strata in
\begin{equation}
\label{eq:TopLevelStratification}
\sM_{\ft}
\cong
\sM^{*,0}_{\ft} \cup M^w_{\ka} \cup
\mathop{\cup}\limits_{\fs}M_{\fs}
\end{equation}
do not intersect when {\em $\sM_{\ft}$ contains no reducible, zero-section
  solutions\/}. (The finite union above over $\fs$ is over the subset of all
\spinc structures for which $M_{\fs}$ is non-empty and for which there is a
  splitting $\ft=\fs\oplus\fs'$.)  

\begin{defn}
\label{defn:AbstractLink}
Let $Z$ be a closed subset of a smooth, Riemannian manifold $M$, and
suppose that $Z = Z_0\cup Z_1$, where $Z_0$ and $Z_1$ are locally closed,
smooth submanifolds of $M$ and $Z_1\subset \barZ_0$. (That is, $Z$ is a
smoothly stratified space with two strata in the sense of \cite[Chapter
11]{MMR}.) Let $N_{Z_1}$ be the normal bundle of $Z_1\subset M$ and let
$\sO'\subset N_{Z_1}$ be an open neighborhood of the zero section
$Z_1\subset N_{Z_1}$ such that there is a diffeomorphism $\gamma$,
commuting with the zero section of $N_{Z_1}$ (so $\gamma|_{Z_1} =
\id_{Z_1}$), from $\sO'$ onto an open neighborhood $\gamma(\sO')$ of
$Z_1\subset M$.  Let $\sO\Subset\sO'$ be an open neighborhood of the zero
section $Z_1\subset N_{Z_1}$, where $\barsO =
\sO\cup\rd\sO \subset \sO'$ is a smooth manifold-with-boundary. Then
$L_{Z_1} = Z_0\cap\gamma(\rd\sO)$ is a {\em link of $Z_1$ in $Z_0$\/}.
\end{defn}

In the preceding definition, $Z_0$ will be the intersection of
$\sM^{*,0}_{\ft}$ with a neighborhood of the image of $M_{\fs}$ and $Z_1$
will be identified with $M_{\fs}$.  Although the ambient manifold $M$
in Definition
\ref{defn:AbstractLink} it is not required to be finite-dimensional, we
shall impose this constraint as we need to define a fundamental class for
the sphere bundle of $N_{Z_1}$.  Thus, we will need to define a
finite-dimensional submanifold of $\sC_{\ft}$ containing---as smooth
submanifolds---the image of $M_{\fs}$ and an open neighborhood in
$\sM_{\ft}$ of the image of $M_{\fs}$.  We fulfill these
requirements in the next subsection by using a globalized, stabilized
version of the local Kuranishi model \cite[\S 4.2.4]{DK}, \cite{Kuranishi}
defined by the deformation complex of the $\PU(2)$ monopole equations
\eqref{eq:PT}. For the remainder of this subsection, we describe the
splitting of this deformation complex, at a reducible $\PU(2)$ monopole,
into deformation complexes which are ``tangential'' and ``normal'' to the
stratum $M_{\fs}$.

\subsubsection{Decomposing the $\PU(2)$ elliptic deformation
sequence at reducible pairs}
\label{subsubsec:ReducibleDeformationDecomp}
In this and the next two sub-subsections we describe how the elliptic
deformation complex for $\sM_{\ft}$ at a point $[A,\Phi]$ in the image of
$M_{\fs}\embed \sM_{\ft}$ can be split into normal and tangential
components, where $\fs=(\rho,W)$ and $\ft=(\rho,V)$ with $V=W\oplus
W\otimes L$ as in Lemma \ref{lem:TopU1Embedding}.

Let $(A,\Phi)=\iota(B,\Psi)$ be a pair in $\tsC_{\ft}$, although not
necessarily a solution to the $\PU(2)$ monopole equations. We begin by
considering the deformation sequence \eqref{eq:ConstDetDefComplex}, namely
\begin{equation}
\begin{CD}
C^\8(F_0) @>{d_{A,\Phi}^0}>> C^\8(F_1)
@>{d_{A,\Phi}^1}>> C^\8(F_2),
\end{CD}
\label{eq:ConstDetDefSequence}
\end{equation}
at a point $(A,\Phi) = (B\oplus B\otimes A_L,\Psi\oplus 0)$, where 
$A_L = A_\Lambda\otimes (B^{\det})^*$ and
the vector bundles $F_j$, $j=0,1,2$, are defined by
\begin{equation}
  \label{eq:DefSequenceBundles}
\begin{aligned}
F_0 &= \fg_{\ft}, 
\\
F_1 &= \Lambda^1\otimes \fg_{\ft} \oplus V^+, 
\\
F_2 &= \Lambda^+\otimes\fg_{\ft} \oplus V^-.
\end{aligned}
\end{equation}
We shall also use $L^2_{k+1-j}(F_j)$, the Hilbert spaces of $L^2_{k+1-j}$
sections for $j=0,1,2$, when applying these sequences.  (The sequence is a
complex if and only $\fS(A,\Phi)=0$.)  It will be convenient to define
vector bundle splittings,
\begin{equation}
\label{eq:FjSequenceSplittings}
F_j \cong F_j^t\oplus F_j^n, \quad j=0,1,2,
\end{equation}
using the canonical isomorphism
$\fg_{\ft}\cong i\underline{\RR}\oplus L$ of Lemma \ref{lem:su(E)Split}:
\begin{equation}
\label{eq:HilbertSpaceRedCplxSplitting}
\begin{matrix}
F_0^t = i\Lambda^0\hfill
&\text{and}
&F_0^n = \Lambda^0\otimes_\RR L = L,\hfill 
\\
F_1^t = i\Lambda^1 \oplus W^+
&\text{and}
&F_1^n = \Lambda^1\otimes_\RR L \oplus W^+\otimes L,
\\
F_2^t = i\Lambda^+ \oplus W^-
&\text{and}
&F_2^n = \Lambda^+\otimes_\RR L \oplus W^-\otimes L.
\end{matrix}
\end{equation}
Conversely, the decompositions \eqref{eq:FjSequenceSplittings} yield
inclusions which we write as
\begin{equation}
\label{eq:FjntSequenceInclusions}
\iota:F_j^t\embed F_j 
\quad\text{and}\quad
\iota:F_j^n\embed F_j, \quad j=0,1,2.
\end{equation}
The motivation for the splitting is due to the fact that the component
$C^\8(F_1^t)$ will contain vectors tangent to $M_{\fs}$ while
$C^\8(F_1^n)$ will contain those vectors normal to $M_{\fs}$.
We note in passing that when the above vector bundles are given their
natural fiber inner products, the splittings $F_j = F_j^t\oplus F_j^n$
define isomorphisms of Riemannian vector bundles when the bundles
$F_j^t\oplus F_j^n$ are given their direct sum fiber inner products.

\subsubsection{Decomposing the group actions on the deformation sequence
  bundles} 
\label{subsubsec:DecompGroupActions}
The embedding $\sG_\fs\embed\sG_\ft$ defined in \eqref{defn:GInclusion} and
two $S^1$ actions on $V=W\oplus W\otimes L$ of interest to us, namely
\eqref{eq:ScalarMultAction} and \eqref{eq:VReducibleAction}, induce
corresponding group actions on the bundles $F_j^t$, $F_j^n$ in
\eqref{eq:HilbertSpaceRedCplxSplitting}, arising in the decomposition
\eqref{eq:ReducibleDefComplexSplitting} of the $\PU(2)$-monopole
deformation sequence \eqref{eq:ConstDetDefSequence}.  For ease of later
reference, we record the induced group actions here on $F_j^t$ and $F_j^n$.

Writing $V=W\otimes E$ and $E=\underline{\CC}\oplus L$, we
defined an embedding \eqref{defn:GInclusion} of gauge groups,
$$
\varrho:\sG_\fs\embed\sG_\ft, \quad s\mapsto
\id_W\otimes\begin{pmatrix}s & 0 \\ 0 & s^{-1}\end{pmatrix}.
$$
Then the homomorphism $\Ad:\Aut(E)\to\Aut(\su(E))$, $u\mapsto
u(\,\cdot\,)u^{-1}$, induces an action of $\varrho(s)$ on $\su(E)\cong
i\underline{\RR}\oplus L$ via the isomorphism $\xi\mapsto (\nu,z)$ in
\eqref{eq:su(E)Split},
\begin{align*}
\Ad(\varrho(s))\xi
=
\begin{pmatrix}s & 0 \\ 0 & s^{-1}\end{pmatrix}
\begin{pmatrix}\nu & -\bar z \\ z & \nu\end{pmatrix}
\begin{pmatrix}s^{-1} & 0 \\ 0 & s\end{pmatrix}
=
\begin{pmatrix}\nu & -s^2\bar z \\ s^{-2}z & \nu\end{pmatrix},
\end{align*}
where $\nu\in C^\8(X,i\RR)$ and $z\in C^\8(X,\CC)$. Hence, the induced
action of $\varrho(s)$ on $\su(E)\cong i\underline{\RR}\oplus L$---via the
composition of the embedding $\varrho:\sG_{\fs}\embed\sG_{\ft}$ given by
\eqref{defn:GInclusion}, the map $\Ad:\Aut(V)\to\Aut(\su(V))$, and the
projection $\su(V)\to\su(E)=\fg_{\ft}$ given by \eqref{eq:OrthogDecompsuV}
---is trivial on the factor $i\underline{\RR}$ and acts as scalar
multiplication by $s^{-2}$ on the complex line bundle $L$:
\begin{equation}
  \label{eq:su(E)SplitGaugeGroupAction}
  \Ad\circ\varrho:\sG_\fs \times (i\underline{\RR}\oplus L)
  \to i\underline{\RR}\oplus L,
  \quad (s,(\nu,z))\mapsto (\nu,s^{-2}z).
\end{equation}
Hence, the induced action of $s\in\sG_\fs$ is trivial on the bundles
$F_j^t$ while on $F_j^n$ it acts as $s^{-2}$ on the factors
$\Lambda^j\otimes_\RR L$ and as $s^{-1}$ on the factors
$W^\pm\otimes L$.

The action \eqref{eq:ScalarMultAction} of $S^1$ on $V=W\oplus W\otimes L$
by scalar multiplication induces the trivial action on $\fg_\ft=\su(E)$.
Hence, $S^1$ acts trivially on the factors $i\Lambda^j$ and
$\Lambda^j\otimes_\RR L$ of $F_j^t$, $F_j^n$, whereas it acts
scalar multiplication on the factors $W^\pm$ and $W^\pm\otimes L$.

Finally, we consider the action \eqref{eq:VReducibleAction} of $S^1$ on
$V=W\oplus W\otimes L$ by the identity on the factor $W$ and as scalar
multiplication on $W\otimes L$. The induced action of $S^1$ is trivial on the
factors $F_j^t$ while it acts by scalar multiplication on the factors $F_j^n$.

\subsubsection{Decomposing the differentials in the deformation
sequence at reducible pairs}
\label{subsubsec:DecompDefOps}
We now describe the splitting of the differentials $d_{A,\Phi}^j$
in the sequence, when $(A,\Phi)=\iota(B,\Psi)$ is a reducible pair in
$\tsC_{\ft}$.

Recall from equation \eqref{eq:ConstDetDefComplex} that 
\begin{equation}
d_{A,\Phi}^1(a,\phi) = (D\fS)_{A,\Phi}(a,\phi) =
\left(\begin{matrix} d^+_{\hat A}a -
\tau\rho^{-1}(\Phi\otimes\phi^* + \phi\otimes\Phi^*)_{00} \\
(D_{A}+\rho(\vartheta))\phi + \rho(a)\Phi \end{matrix}\right).
\label{eq:DefnLinearization}
\end{equation}
If $\Phi=\Psi\oplus 0$ and $\phi = \psi\oplus\psi'$, the quadratic
term takes the form
$$
\Phi\otimes\phi^* + \phi\otimes\Phi^*
=
\begin{pmatrix}
\Psi\otimes\psi^* + \psi\otimes\Psi^* & \Psi\otimes\psi^{\prime,*} \\
\psi'\otimes\Psi^* & 0
\end{pmatrix},
$$
and so
\begin{align*}
&\tau\rho^{-1}(\Phi\otimes\phi^* + \phi\otimes\Phi^*)_{00} \\
&\qquad =
\begin{pmatrix}
\half\tau\rho^{-1}(\Psi\otimes\psi^* + \psi\otimes\Psi^*)_0 &
\tau\rho^{-1}(\Psi\otimes\psi^{\prime,*})_0 \\
\tau\rho^{-1}(\psi'\otimes\Psi^*)_0 &
-\half\tau\rho^{-1}(\Psi\otimes\psi^* + \psi\otimes\Psi^*)_0
\end{pmatrix},
\end{align*}
where the diagonal term $\half\tau\rho^{-1}(\Psi\otimes\psi^* +
\psi\otimes\Psi^*)_0$ is in $C^\8(i\Lambda^+)$ while the off-diagonal
term $\tau\rho^{-1}(\psi'\otimes\Psi^*)_0$ is in
$C^\8(\Lambda^+\otimes L)$. Indeed, we have
$$
(\psi'\otimes\Psi^*)_0
\in
C^\8(\fsl(W^+)\otimes_{\CC} L)
=
C^\8(\su(W^+)\otimes_{\RR} L),
$$
since $\fsl(W^+)=\su(W^+)\otimes_\RR\CC$.  The complex-linear Clifford map
restricts to a real-linear isomorphism $\rho:\Lambda^+\to\su(W^+)$, so we
have $\rho^{-1}(\psi'\otimes\Psi^*)_0 \in C^\8(\Lambda^+\otimes L)$.

The Dirac-operator term
$D_{A}+\rho(\vartheta)+\rho(a):C^\8(V^+)\to C^\8(V^-)$ splits as
$$
(D_{A}+\rho(\vartheta))\phi + \rho(a)\Phi
= \begin{pmatrix}
(D_{B}+\rho(\vartheta))\psi + \rho(\al)\Psi \\
(D_{B\otimes A_L}+\rho(\vartheta))\psi' + \rho(\be)\Psi
\end{pmatrix},
$$
noting that 
$$
a = \begin{pmatrix}\al & -\barbe \\ \be & -\al\end{pmatrix},
\quad\text{with}\quad
\al\in C^\8(i\Lambda^1)
\quad\text{and}
\quad\be\in C^\8(\Lambda^1\otimes L).
$$
By Lemma \ref{lem:ReducibleSpinu}, $\hat A=d_{\RR}\oplus A_L$.  The term
$d_{\hat A}^+ a\in C^\8(\Lambda^+\otimes \fg_{\ft})$ then splits as
\begin{equation}
\label{eq:dAplusSplitting}
d_{\hat A}^+a =
\begin{pmatrix}d^+\al & -d_{A_L^*}^+\barbe
\\ d_{A_L}^+\be & -d^+\al\end{pmatrix},
\end{equation}
where the diagonal term $d^+\al$ is in $ C^\8(i\Lambda^+)$, and 
off-diagonal term $d_{A_L}^+\be$ is in $ C^\8(\Lambda^+\otimes L)$.

With respect to these identifications, the linear operator
$d_{\iota(B,\Psi)}^1: C^\8(F_1)\to C^\8(F_2)$ diagonalizes to give
$$
d_{\iota(B,\Psi)}^1
=
\begin{pmatrix}
d_{\iota(B,\Psi)}^{1,t} & 0 \\
0 & d_{\iota(B,\Psi)}^{1,n}
\end{pmatrix}:
 C^\8(F_1^t)\oplus C^\8(F_1^n)\to  C^\8(F_2^t)\oplus C^\8(F_2^n),
$$
where $d_{\iota(B,\Psi)}^{1,t}: C^\8(F_1^t) \to  C^\8(F_2^t)$ and
$d_{\iota(B,\Psi)}^{1,n}: C^\8(F_1^n) \to  C^\8(F_2^n)$ are the
``tangential'' and ``normal'' components of $d_{\iota(B,\Psi)}^1$.  More
explicitly, the operator
$$
d_{\iota(B,\Psi)}^1:
\begin{matrix}
 C^\8(i\Lambda^1) \oplus  C^\8(W^+), \\
\oplus \\
 C^\8(\Lambda^1\otimes L)\oplus C^\8(W^+\otimes L)
\end{matrix}
\too
\begin{matrix}
 C^\8(i\Lambda^+) \oplus  C^\8(W^-), \\
\oplus \\
 C^\8(\Lambda^+\otimes L)\oplus C^\8(W^-\otimes L)
\end{matrix}
$$
diagonalizes and the tangential component of $d_{\iota(B,\Psi)}^1$ is
given by
$$
d_{\iota(B,\Psi)}^{1,t}:
 C^\8(i\Lambda^1) \oplus  C^\8(W^+)
\too
 C^\8(i\Lambda^+) \oplus  C^\8(W^-),
$$
where we define 
\begin{equation}
d_{\iota(B,\Psi)}^{1,t}(\al,\psi)
=
\begin{pmatrix}
d^+\al - \half\tau\rho^{-1}(\Psi\otimes\psi^* + \psi\otimes\Psi^*)_0 \\
(D_{B}+\rho(\vartheta))\psi + \rho(\al)\Psi
\end{pmatrix},
\label{eq:ReducibleTangentDerivative}
\end{equation}
Note that $d_{\iota(B,\Psi)}^{1,t}$ matches the Seiberg-Witten
differential $d_{B,\Psi}^1$ in equation \eqref{eq:d1SW}, aside from a
factor of $2$ in the $d^+$ component: the scaling factor has no
significance since we are only interested in the kernels and cokernels of
these operators. The normal component of $d_{\iota(B,\Psi)}^1$ is given by
$$
d_{\iota(B,\Psi)}^{1,n}:
 C^\8(\Lambda^1\otimes L)\oplus C^\8(W^+\otimes L)
\too
 C^\8(\Lambda^+\otimes L)\oplus C^\8(W^-\otimes L),
$$
where we define
\begin{equation}
d_{\iota(B,\Psi)}^{1,n}(\be,\psi')
=
\begin{pmatrix}
d^+_{A_L}\beta
   -\tau\rho^{-1}(\psi'\otimes\Psi^*)_0 \\
(D_{B\otimes A_L}+\rho(\vartheta))\psi'+\rho(\beta)\Psi
\end{pmatrix}.
\label{eq:ReducibleNormalDerivative}
\end{equation}
This completes our description of the diagonalization of the differential
$d_{\iota(B,\Psi)}^1$.

It remains to discuss the rather simpler diagonalization of the
differential $d_{\iota(B,\Psi)}^0$ with respect to the splitting of
$ C^\8(F_0)$ and $ C^\8(F_1)$ into normal and tangential
components. {}From Lemma
\ref{lem:su(E)Split} we see that elements $\zeta = \zeta_0$ of
$C^\8(\fg_{\ft})$, recalling that
$L^2_{k+1}(\fg_{\ft})= T_{\id}\sG_{\ft}$, may be uniquely written as
\begin{equation}
\zeta
=
\begin{pmatrix}
f            &-\bar\kappa \\
\kappa &-f
\end{pmatrix}
,
\label{eq:GaugeGroupLieAlgDecomp}
\end{equation}
where $f\in C^\8(X,i\RR)$ is an imaginary-valued
function and $\kappa \in C^\8(L)$. Thus,
$$
\zeta  = f \oplus\kappa\in C^\8(F_0^t)\oplus C^\8(F_0^n).
$$
Recall from equation \eqref{eq:d0PT} that the
differential $d^0_{A,\Phi}$ is defined by
$$
d^0_{A,\Phi}\zeta = (d_{\hat A}\zeta,-\zeta\Phi)
\in  C^\8(\Lambda^1\otimes\fg_{\ft})\oplus  C^\8(V^+),
\qquad \zeta\in  C^\8(\Lambda^1\otimes\fg_{\ft}).
$$
Note that
\begin{align*}
\zeta\Phi
&=
\begin{pmatrix}
f              &-\bar\kappa \\
\kappa &-f
\end{pmatrix}
\cdot
\begin{pmatrix}
\Psi \\
0
\end{pmatrix}
=
\begin{pmatrix}
f\Psi \\
\kappa\Psi
\end{pmatrix}, \\
d_{\hat A}\zeta
&=
\begin{pmatrix}
df                                   &-d_{A_L^*}\bar\kappa \\
d_{A_L}\kappa  &-df
\end{pmatrix}.
\end{align*}
So, at the pair $(A,\Phi) =
\iota(B,\Psi) = (B\oplus B\otimes A_L,\Psi\oplus 0)$ we have 
$$
d^0_{\iota(B,\Psi)}\zeta
=
\left(
\begin{pmatrix}
df                                   &-d_{A_L^*}\bar\kappa \\
d_{A_L}\kappa  &-df
\end{pmatrix},
\begin{pmatrix}
-f\Psi \\
-\kappa\Psi
\end{pmatrix}
\right).
$$
With respect to these identifications, the linear operator
$d_{\iota(B,\Psi)}^0: C^\8(F_0)\to C^\8(F_1)$
diagonalizes to give
$$
d_{\iota(B,\Psi)}^0
=
\begin{pmatrix}
d_{\iota(B,\Psi)}^{0,t} & 0 \\
0 & d_{\iota(B,\Psi)}^{0,n}
\end{pmatrix}:
 C^\8(F_0^t)\oplus C^\8(F_0^n)\too  C^\8(F_1^t)\oplus C^\8(F_1^n),
$$
where $d_{\iota(B,\Psi)}^{0,t}: C^\8(F_0^t) \to  C^\8(F_1^t)$ and
$d_{\iota(B,\Psi)}^{0,n}: C^\8(F_0^n) \to  C^\8(F_1^n)$
are the ``tangential'' and
``normal'' components of $d_{\iota(B,\Psi)}^0$.
More explicitly, the operator
$$
d_{\iota(B,\Psi)}^0:
\begin{matrix}
 C^\8(i\Lambda^0))\\
\oplus \\
 C^\8(L)
\end{matrix}
\too
\begin{matrix}
 C^\8(i\Lambda^1) \oplus  C^\8(W^+) \\
\oplus \\
 C^\8(\Lambda^1\otimes L)\oplus C^\8(W^+\otimes L)
\end{matrix}
$$
diagonalizes and the tangential component of $d_{\iota(B,\Psi)}^0$ is
given by
\begin{equation}
d_{\iota(B,\Psi)}^{0,t} :
 C^\8(i\Lambda^0)
\too
 C^\8(i\Lambda^1)\oplus C^\8(W^+),
\label{eq:ReducibleTangentGauge}
\end{equation}
where we define
$$
d_{\iota(B,\Psi)}^{0,t}f
= (df,-f\Psi).
$$
Note that $d_{\iota(B,\Psi)}^{0,t}$ matches the gauge group
differential $d_{B,\Psi}^0$ in equation \eqref{eq:d0SW}. The normal
component of $d_{\iota(B,\Psi)}^0$ is given by 
\begin{equation}
d_{\iota(B,\Psi)}^{0,n}:
 C^\8(L)
\too
 C^\8(\Lambda^1\otimes L)\oplus C^\8(W^+\otimes L)
\label{eq:ReducibleNormalGauge},
\end{equation}
where we define
$$
d_{\iota(B,\Psi)}^{0,n}\kappa
=(d_{A_L}\kappa,-\kappa\Psi).
$$
This completes our description of the diagonalization of the differential
$d_{\iota(B,\Psi)}^0$.

Recall from \cite[Equation (2.38)]{FL1} that the $L^2$ adjoint of
$d_{A,\Phi}^0$ is given by
\begin{equation}
\label{eq:DefineSliceOperator}
d_{A,\Phi}^{0,*}(a,\phi) = d_{\hat A}^*a - (\cdot\Phi)^*\phi,
\qquad (a,\phi) \in  C^\8(\Lambda^1\otimes\fg_{\ft})\oplus C^\8(V^+).
\end{equation}
The $L^2$ adjoint $d_{\iota(B,\Psi)}^{0,*}: C^\8(F_1)\to  C^\8(F_0)$
then diagonalizes to give
$$
d_{\iota(B,\Psi)}^{0,*}
=
\begin{pmatrix}
d_{\iota(B,\Psi)}^{0,t,*} & 0 \\
0 & d_{\iota(B,\Psi)}^{0,n,*}
\end{pmatrix}:
 C^\8(F_1^t)\oplus C^\8(F_1^n)\to
 C^\8(F_0^t)\oplus  C^\8(F_0^n),
$$
where $d_{\iota(B,\Psi)}^{0,t,*}: C^\8(F_1^t) \to  C^\8(F_0^t)$
is given explicitly by
$$
d_{\iota(B,\Psi)}^{0,t,*} :
 C^\8(i\Lambda^1) \oplus  C^\8(W^+)
\to
 C^\8(i\Lambda^0),
$$
with
\begin{equation}
d_{\iota(B,\Psi)}^{0,t,*}(\al,\psi)
= d^*\al - \Imag\langle\Psi,\psi\rangle,
\label{eq:ReducibleTangentGaugeAdjoint}
\end{equation}
and $d_{\iota(B,\Psi)}^{0,n,*}: C^\8(F_1^n) \to  C^\8(F_0^n)$ is
explicitly given by
$$
d_{\iota(B,\Psi)}^{0,n,*} :
 C^\8(\Lambda^1\otimes L)\oplus C^\8(W^+\otimes L)
\to
 C^\8(L),
$$
with
\begin{equation}
d_{\iota(B,\Psi)}^{0,n,*}(\beta,\psi')
= d^*_{A_L}\beta + \psi'\otimes\Psi^*.
\label{eq:ReducibleNormalGaugeAdjoint}
\end{equation}
Of course, these $L^2$ adjoints are defined in the usual way by
\begin{align*}
(d^*\al - \Imag\langle\Psi,\psi\rangle,f)_{L^2}
&=
(\al,df)_{L^2} - (\psi,f\Psi)_{L^2},
\\
(d^*_{A_L}\beta + \psi'\otimes\Psi^*,\kappa)_{L^2}
&=
(\beta,d_{A_L}\kappa)_{L^2}
    + (\psi',\kappa\Psi)_{L^2},
\end{align*}
for all $f\in C^\8(i\Lambda^0)$ and $\kappa\in C^\8(L)$,
where
\begin{align*}
(f\Psi,\psi)_{L^2} &= (f,\overline{\langle\psi,\Psi\rangle})_{L^2}
= (f,\langle\Psi,\psi\rangle)_{L^2}, \\
(\kappa\Psi,\psi')_{L^2}
& = (\kappa,\psi'\otimes\Psi^*)_{L^2}.
\end{align*}
This completes our description of the diagonalization of the $L^2$ adjoint
$d_{\iota(B,\Psi)}^{0,*}$.

The above discussion implies that the deformation sequence
\eqref{eq:ConstDetDefSequence} splits when $(A,\Phi)=\iota(B,\Psi)$:
\begin{equation}
\label{eq:ReducibleDefComplexSplitting}
\begin{gathered}
\begin{CD}
 C^\8(F_0^t) @> d^{0,t}_{\iota(B,\Psi)} >>
 C^\8(F_1^t) @> d^{1,t}_{\iota(B,\Psi)} >>
 C^\8(F_2^t),
\end{CD} \\
\begin{CD}
 C^\8(F_0^n) @> d^{0,n}_{\iota(B,\Psi)} >>
 C^\8(F_1^n) @> d^{1,n}_{\iota(B,\Psi)} >>
 C^\8(F_2^n).
\end{CD}
\end{gathered}
\end{equation}
These sequences form complexes if and only if $\fS(A,\Phi)=0$, in which
case they give the {\em tangential\/} and {\em normal deformation
complexes\/} for the $\PU(2)$ monopole equations at a reducible solution,
so it is convenient to consider the rolled-up sequences
\begin{gather*}
\begin{CD}
 C^\8(F_1^t) @> \sD^t_{\iota(B,\Psi)} >>
 C^\8(F_0^t)\oplus  C^\8(F_2^t),
\end{CD}
\\
\begin{CD}
 C^\8(F_1^n) @> \sD^n_{\iota(B,\Psi)} >>
 C^\8(F_0^n)\oplus  C^\8(F_2^n)
\end{CD},
\end{gather*}
comprising the tangential and normal components of the rolled-up
deformation sequence
$$
\sD_{\iota(B,\Psi)}:
 C^\8(F_1) \to  C^\8(F_0)\oplus  C^\8(F_2).
$$
for the $\PU(2)$ monopole equations.
The tangential and normal components of the deformation operator
$\sD_{\iota(B,\Psi)} = d^{0,*}_{\iota(B,\Psi)} +
d^1_{\iota(B,\Psi)}$ are given by
$$
\sD_{\iota(B,\Psi)}^t = d^{0,t,*}_{\iota(B,\Psi)} +
d^{1,t}_{\iota(B,\Psi)}
\quad\text{and}\quad
\sD_{\iota(B,\Psi)}^n = d^{0,n,*}_{\iota(B,\Psi)} +
d^{1,n}_{\iota(B,\Psi)},
$$
so that we have a decomposition
\begin{equation}
\sD_{\iota(B,\Psi)}
= \sD^t_{\iota(B,\Psi)}\oplus\sD^n_{\iota(B,\Psi)}.
\label{eq:DefOperatorReducSplit}
\end{equation}
More explicitly, we see from
\eqref{eq:ReducibleTangentDerivative} and
\eqref{eq:ReducibleTangentGaugeAdjoint}
that the tangential component
\begin{equation}
\sD_{\iota(B,\Psi)}^t :
 C^\8(i\Lambda^1)
\oplus C^\8(W^+)
\too
 C^\8(i\Lambda^0)
\oplus C^\8(i\Lambda^+)
\oplus C^\8(W^-)
\label{eq:TangentReducDefOperator}
\end{equation}
of $\sD_{\iota(B,\Psi)}$ is given by
$$
\sD^t_{\iota(B,\Psi)}(\al,\psi)
=
\begin{pmatrix}
d^*\al - \Imag\langle\Psi,\psi\rangle \\
d^+\al
- \half\tau\rho^{-1}(\Psi\otimes\psi^* + \psi\otimes\Psi^*)_0 \\
(D_{B}+\rho(\vartheta))\psi + \rho(\al)\Psi
\end{pmatrix},
$$
while from \eqref{eq:ReducibleNormalDerivative} and
\eqref{eq:ReducibleNormalGaugeAdjoint}
the normal component
\begin{equation}
\sD_{\iota(B,\Psi)}^n :
 C^\8(\Lambda^1\otimes L)
\oplus
 C^\8(W^+\otimes L)
\too
 C^\8(L)
\oplus
 C^\8(\Lambda^+\otimes L)
\oplus
 C^\8(W^-\otimes L)
\label{eq:NormalReducDefOperator}
\end{equation}
of $\sD_{\iota(B,\Psi)}$ is given by
$$
\sD^n_{\iota(B,\Psi)}(\beta,\psi')
=
\begin{pmatrix}
d^*_{A_L}\beta + \psi'\otimes\Psi^* \\
d^+_{A_L}\be
- \tau\rho^{-1}(\psi'\otimes\Psi^*)_0 \\
(D_{B\otimes A_L}+\rho(\vartheta))\psi'
+ \rho(\beta)\Psi
\end{pmatrix}.
$$
Note the operators $d^{0,n}_{\iota(B,\Psi)}$ and
$d^{1,n}_{\iota(B,\Psi)}$, and thus $\sD^n_{\iota(B,\Psi)}$, are
complex linear when the vector bundles
$F_j^n$, $j=0,1,2$, are given the natural
complex structures induced by the $S^1$ action on $L$.

If $\fS(A,\Phi)=0$, the two sequences in
\eqref{eq:ReducibleDefComplexSplitting} define cohomology groups
$H^{\bullet,t}_{\iota(B,\Psi)}$ and $H^{\bullet,n}_{\iota(B,\Psi)}$,
respectively, which we may compare with the cohomology groups
$H^{\bullet}_{\iota(B,\Psi)}$ and $H^{\bullet}_{B,\Psi}$ of the $\PU(2)$
monopole deformation complex \eqref{eq:ConstDetDefComplex} at the pair
$\iota(B,\Psi)\in\tilde\sM_{\ft}$ and the cohomology groups
$H^{\bullet}_{B,\Psi}$ of the Seiberg-Witten deformation complex
\eqref{eq:SWDeformationComplex} at the pair $(B,\Psi)\in\tilde M_{\fs}$:

\begin{lem}
\label{lem:IdentityOfReducibleCohomology}
Continue the above notation and require that $\iota(B,\Psi)$ be a
(reducible) $\PU(2)$ monopole. For the elliptic deformation complex
\eqref{eq:ConstDetDefComplex} defined
by the $\PU(2)$ monopole equations \eqref{eq:PT} and the automorphism group
$\sG_{\ft}$, we have the following canonical isomorphisms of cohomology groups,
$$
H^{\bullet,t}_{\iota(B,\Psi)} \cong H^{\bullet}_{B,\Psi},
$$
and hence,
$$
H^{\bullet}_{\iota(B,\Psi)} \cong H^{\bullet}_{B,\Psi}
\oplus H^{\bullet,n}_{\iota(B,\Psi)}.
$$
Moreover,  if $\Psi\not\equiv 0$ then
$H^{0}_{B,\Psi}=0$, while if
$H^{2}_{B,\Psi}=0$ then $H^2_{\iota(B,\Psi)} \cong
H^{2,n}_{\iota(B,\Psi)}$.
\end{lem}

\begin{proof}
The isomorphisms identifying the cohomology of the tangential deformation
complex with that of the Seiberg-Witten complex follow
immediately from a comparison of these two complexes.

If $\Psi\not\equiv 0$ then $H^{0}_{B,\Psi}=0$
and \eqref{eq:ReducibleNormalGauge} implies
$H^{0,n}_{\iota(B,\Psi)}=0$ and so
$H^0_{\iota(B,\Psi)} = 0$. Moreover, if
$H^{2}_{B,\Psi}=0$, then $H^2_{\iota(B,\Psi)} =
H^{2,n}_{\iota(B,\Psi)}$.
\end{proof}

Transversality for the Seiberg-Witten equations (Proposition
\ref{prop:SmoothFamilyOfReducibles}) implies $H^{2}_{B,\Psi}=0$.

\subsection{The link of a stratum of reducible monopoles: global structure}
\label{subsec:LinkOfReduciblesGlobal}
Our task in this subsection is to construct an ambient finite-dimensional,
smooth submanifold $\sM_\ft(\Xi,\fs)\subset\sC_{\ft}$ containing $M_{\fs}$
as a smooth submanifold, as required by the definition of the link in
Definition \ref{defn:AbstractLink}.  Recall that $\fS$ is the section of
the $S^1$-equivariant, infinite-rank ``obstruction bundle''
\begin{equation}
\fV = \tsC_{\ft}\times_{\sG_{\ft}} L^2_{k-1}(F_2) \to \sC_{\ft}
\label{eq:S1EquivInfRankObstBundle}
\end{equation}
defined by the $\PU(2)$ monopole equations \eqref{eq:PT} and this
section need not vanish transversely along $M_{\fs}$.

To motivate the construction of the ambient manifold $\sM_\ft(\Xi,\fs)$
given in this subsection, suppose temporarily that the cokernel of
$\sD_{A,\Phi}$ has constant rank as $[A,\Phi]$ varies in the image of
$M_{\fs}\embed \sM_{\ft}$ (that is, no spectral flow occurs). Then we
obtain a finite-rank, $S^1$-equivariant vector bundle $\Coker\bsD$ over
$M_{\fs}$, with fibers $\Coker\sD_{A,\Phi}$.  Let $2\nu$ be the least
positive eigenvalue of the Laplacian $\Delta_{A,\Phi}
=\sD_{A,\Phi}\sD_{A,\Phi}^*$ as $[A,\Phi]$ varies along the image of the
compact manifold $M_{\fs}$ and let $\Pi_{A,\Phi;\nu}$ denote the $L^2$
orthogonal projection from $L^2_{k-1}(F_2)$ onto the subspace spanned by the
eigenvectors of $\De_{A,\Phi}$ with eigenvalue less than $\nu$.  The vector
bundle $\Coker\bsD$ over $\iota(M_{\fs})\subset\sC_\ft$ then extends to a
vector bundle $\Xi_\nu= \Ker \Pi_\nu^\perp\bDelta =
\Coker\Pi_\nu^\perp\bsD$, where $\Pi_\nu^\perp=\id-\Pi_\nu$, of the same
rank over an open neighborhood of $\iota(M_{\fs})\subset \sC_{\ft}$.
Arguing as in Lemma \ref{lem:ReducibleSubmanifold}, one can see that the
space $M_{\fs}$ would be a smooth submanifold of $\sC_{\ft}$.  Then both
$M_{\fs}$ and an open neighborhood in $\sM_{\ft}^{*,0}$ of $\iota(M_{\fs})$
would be smooth submanifolds of the $S^1$-invariant ``thickened'' 
moduli space,
$$
\sM_{\ft}(\Xi_\nu,\fs) = (\Pi_\nu^\perp\fS)^{-1}(0) \subset\sC_{\ft},
$$
a finite-dimensional, smooth, $S^1$-invariant manifold  which serves
as the ambient, finite-dimensional, smooth manifold ``$M$'' of Definition
\ref{defn:AbstractLink}. Ambient manifolds of this form,
defined by spectral projections as above, have been used
by Donaldson \cite{DonPoly}, Friedman-Morgan \cite{FrM}, and Taubes
\cite{TauIndef}, \cite{TauStable} to describe
neighborhoods of the stratum $\{[\Theta]\}\times\Sym^\kappa(X)$ containing
the trivial connection $\Theta$ in the Uhlenbeck compactification of the
moduli space $M_\kappa^w$ of anti-self-dual connections when $w=0$.

In practice, one cannot guarantee that $\Coker\bsD$ will either vanish or
even have constant rank due to spectral flow, so we must resort to a more
general construction of an $S^1$-equivariant vector bundle $\Xi$ over an
open neighborhood of the image of $M_{\fs}$ in $\sC_{\ft}$ which ``spans''
$\Coker\bsD$ along $M_{\fs}$. The method we employ is an extension of one
used by Atiyah and Singer to construct the index bundle or determinant-line
bundle of a family of elliptic operators \cite[pp.  122--127]{AS4},
\cite[\S 5.1.3 \& \S 5.2.1]{DK}. The Atiyah-Singer method has also been
exploited by Furuta \cite{Furuta}, T-J. Li \& Liu \cite{LiLiu}, J. Li \&
Tian \cite{LiTian}, and Ruan \cite{RuanGafa}, \cite{RuanSW}, \cite{RuanGW}
to construct certain ``global Kuranishi models'' (or ``virtual'' or
``thickened'' moduli spaces) parameterizing spaces of Seiberg-Witten
monopoles or pseudo-holomorphic curves; related ideas are contained in
\cite{Brussee}, \cite{SiebertGW}.  The principal difference between our
construction and those of Li-Liu, Li-Tian, or Ruan is that the vector
bundle replacing $\Coker\bsD$ is defined over an open neighborhood in the
original configuration space $\sC_{\ft}$ rather than on an artificial,
augmented configuration space, such as $\CC^r\times \sC_{\ft}$. In this
sense, our construction is closer to that of \cite{DonPoly}, \cite{FrM},
\cite{TauIndef}, \cite{TauStable} and corresponds to the alternative
construction of index bundles discussed in \cite[pp. 153--166]{AtiyahK},
\cite[\S 1.7]{BoossBleecker}; we find this second stabilization technique
more convenient when constructing links of lower-level strata of reducibles
via gluing in \cite{FL3}, \cite{FL4}, \cite{FLLevelOne}.

\subsubsection{Smooth embeddings}
We recall that our abstract definition of a link of the lower stratum of a
two-stratum space (see Definition \ref{defn:AbstractLink}) requires an
ambient smooth manifold containing the strata. As a first step in the
construction of a {\em finite-dimensional\/} ambient manifold, we use the
decompositions of \S \ref{subsec:LinkOfReduciblesLocal} to show that the
topological embedding $M_{\fs}\hookrightarrow \sC^0_{\ft}$ of Lemma
\ref{lem:RedAreSW} is smooth, so $M_{\fs}$ is a smooth submanifold of
$\sC^0_{\ft}$.

\begin{lem}
\label{lem:ReducibleSubmanifold}
Let $\ft$ be a \spinu structure on a closed, oriented, smooth
four-manifold $X$, with reduction $\ft=\fs\oplus \fs\otimes L$.
Let $\iota:\sC^0_{\fs}\hookrightarrow\sC^0_{\ft}$ be the map in Lemma
\ref{lem:RedAreSW}. Then the following hold:
\begin{itemize}
\item
The map $\iota:\sC^0_{\fs}\to\sC^0_{\ft}$ is a smooth immersion, and
\item The map $\iota:M^{0}_{\fs}\hookrightarrow\sC^0_{\ft}$ is a smooth
  embedding, so $M^{0}_{\fs}$ is a submanifold of $\sC^0_{\ft}$.
\end{itemize}
\end{lem}

\begin{proof}
Let $[B,\Psi]$ be a point in $\sC^0_{\fs}$: by the slice
theorem for
$\sC^0_{\fs}$, the restriction of the projection map $\tsC^0_{\fs}\to
\sC^0_{\fs}= \tsC^0_{\fs}/\sG_{\fs}$ to a small enough open
neighborhood of $(B,\Psi)$ in the slice $(B,\Psi)+\Ker
d_{B,\Psi}^{*,0}$ gives a smooth parameterization of an open
neighborhood of $[B,\Psi]$ in $\sC^0_{\fs}$. Similarly, by the slice
theorem \cite[Proposition 2.8]{FL1} for $\sC^0_{\ft}$, the restriction of the
projection map $\tsC^0_{\ft}\to \sC^0_{\ft}=\tsC^0_{\ft}/\sG_{\ft}$ to a
small enough open neighborhood of $\iota(B,\Psi)$ in the slice
$\iota(B,\Psi)+\Ker d_{\iota(B,\Psi)}^{*,0}$ gives a smooth
parameterization of an open neighborhood of $[\iota(B,\Psi)]$ in
$\sC^0_{\ft}$.
Comparing the operators $d^{0,*}_{B,\Psi}$ in
equation \eqref{eq:d0SWL2adjoint}
and $d^{0,t,*}_{\iota(B,\Psi)}$ in equation
\eqref{eq:ReducibleTangentGaugeAdjoint},
we see that the differential of the smooth embedding
$$
\iota:\tsC^0_{\fs} \to \tsC^0_{\ft},
\qquad
(B,\Psi) \mapsto
(B\oplus B\otimes A_L,\Psi\oplus 0),
$$
of Lemma \ref{lem:TopU1Embedding} restricts to an isomorphism
$$
D\iota:\Ker d_{B,\Psi}^{*,0}
\cong
\Ker d_{\iota(B,\Psi)}^{*,0,t}
\subset
\Ker d^{0,t,*}_{A,\Phi}\oplus\Ker d^{0,n,*}_{A,\Phi}
=
\Ker d_{B,\Psi}^{*,0}.
$$
Hence, the differential $D\iota:T_{B,\Psi}\sC^0_{\fs} \to
T_{\iota(B,\Psi)}\sC^0_{\ft}$ is injective and the induced maps
$\iota:\sC^0_{\fs}\to \sC^0_{\ft}$ and $\iota:M^{0}_{\fs}\to
\sC^0_{\ft}$ are smooth immersions.  Since the map
$\iota:M^{0}_{\fs}\to \sC^0_{\ft}$ is a topological embedding according
to Lemma \ref{lem:RedAreSW}, and is a smooth immersion, the map
$\iota:M^{0}_{\fs}\to \sC^0_{\ft}$ is a smooth embedding.
\end{proof}

\subsubsection{Construction of the stabilized cokernel bundle}
Again, we assume that there are no zero-section reducible monopoles
in $M_{\fs}$.
For any representative $(A,\Phi)$ of a point in the image of
$\iota:M_{\fs}\embed \sM_{\ft}$, we have $\Coker (D\fS)_{A,\Phi} = \Coker
d_{A,\Phi}^{1,n} \subset
L^2_{k-1}(F_2^n)$, since $(D\fS)_{A,\Phi} = d_{A,\Phi}^1$ (by definition) and
$\Coker d_{A,\Phi}^{1,t} = 0$. The splitting of the last term of the
elliptic deformation complex \eqref{eq:ConstDetDefSequence} at a reducible
$\PU(2)$ monopole into tangential and normal deformation components
\eqref{eq:ReducibleDefComplexSplitting} yields a
splitting of the corresponding complexes of Hilbert bundles,
\begin{equation}
\label{eq:S1EquivInfRankHilbertObstBundle}
\fV
=
\tsC_{\ft}^0\times_{\sG_{\ft}}L^2_{k-1}(F_2)
\cong
\tsC_{\ft}^0\times_{\sG_{\ft}}L^2_{k-1}(F_2^t\oplus F_2^n),
\end{equation}
when restricted to $\iota(M_{\fs})$.
This splitting can be seen by using a reduction of the structure group
$\sG_{\ft}$ to $\sG_{\fs}$:
\begin{equation}
\label{eq:DefineV}
\begin{aligned}
\fV|_{\iota(M_{\fs})}
&\cong
\iota(\tM_{\fs})\times_{\sG_{\fs}}L^2_{k-1}(F_2)
\cong
\fV^t\oplus\fV^n, 
\\
\text{where}\quad
\fV^t &= \iota(\tM_{\fs})\times_{\sG_{\fs}}L^2_{k-1}(F_2^t),
\\
\fV^n&= \iota(\tM_{\fs})\times_{\sG_{\fs}}L^2_{k-1}(F_2^n)
\end{aligned}
\end{equation}
The action \eqref{eq:VReducibleAction} of $S^1$ on $V=W\oplus W\otimes L$
induces the trivial action of $S^1$ on the bundle $\fV^t$ and the
standard action by complex multiplication on $\fV^n$; see the action
\eqref{eq:SO(3)ReducibleAction} of $S^1$ on $\fg_{\ft} \cong
i\underline{\RR}\oplus L$ induced by the action of $S^1$ on $V$.

Note that although the splitting \eqref{eq:DefineV} is not well-defined
away from the stratum $\iota(M_{\fs})$, since the full gauge group
$\sG_{\ft}$ does not preserve the splitting of the fiber, $L^2_{k-1}(F_2)
\cong L^2_{k-1}(F_2^t)\oplus L^2_{k-1}(F_2^n)$, the circle action
on $L^2_{k-1}(F_2)$ induced by the circle action
\eqref{eq:VReducibleAction} on $V$ does preserve this splitting, with the
circle acting trivially on $L^2_{k-1}(F_2^t)$ and by complex multiplication
on $L^2_{k-1}(F_2^n)$. Hence, the vector bundle
\eqref{eq:S1EquivInfRankHilbertObstBundle} is $S^1$ equivariant with
respect to the circle action on the fiber $L^2_{k-1}(F_2^t\oplus F_2^n)$
induced by the circle action \eqref{eq:VReducibleAction} on $V$.

We use the preceding observation below to define a finite-rank, real
subbundle $\Xi\subset\fV$ with an almost complex structure when
restricted to $\iota(M_{\fs})$. (If the requirement that $\Xi\to
\iota(M_{\fs})$ has an almost complex structure were dropped then we would
not require the splitting into tangential and normal components in order to
construct $\Xi$.)  We shall need the following well-known consequence of
Kuiper's result that the unitary group of a Hilbert space is contractible
\cite[Theorem 3]{Kuiper} (see also \cite[\S I.7, p. 67]{BoossBleecker} or
\cite[p. 208]{LM}):

\begin{thm}
\cite[p. 29]{Kuiper}
\label{Kuiper}
Let $M$ be a compact topological space or a space which has the homotopy
type of a CW complex.  Then every vector bundle over $M$, with fiber an
infinite-dimensional, real or complex separable Hilbert space $\fH$ and
structure group $\GL(\fH)$, is trivial.
\end{thm}

We have the following version of Kuiper's result in the smooth category:

\begin{cor}
\label{KuiperSmooth}
Let $M$ be a smooth manifold which is compact or has the homotopy type of a
CW complex.  Then every $C^\8$ vector bundle over $M$, with fiber an
infinite-dimensional, real or complex separable Hilbert space $\fH$ and
structure group $\GL(\fH)$, has a global $C^\8$ trivialization.
\end{cor}

\begin{proof}
If $\fW$ is a vector bundle over $M$ with fiber a Hilbert space $\fH$, then
Theorem \ref{Kuiper} yields a $C^0$ trivialization $\btau:\fW\to
M\times\fH$, that is, a $C^0$ section $\btau$ of the $C^\8$ Hilbert bundle
$\Hom(\fW,M\times\fH)$ over $M$ which gives a linear isomorphism on each
fiber. Now suppose that $\btau_\8$ is a $C^\8$ section of the bundle
$\Hom(\fW,M\times\fH)$. If $\btau_\8$ is chosen so that
$\|\btau(p)-\btau_\8(p)\|$ is sufficiently small for each $p\in M$, then
$\btau_\8$ is an isomorphism on each fiber and gives the desired global,
$C^\8$ trivialization.
\end{proof}

Note that if we are given a vector bundle $\fW$ over a space $M$ with fiber
a separable Hilbert space $\fH$ and structure group $G\subset \U(\fH)$,
then we have an isomorphism of vector
bundles $\fW \cong \Fr(\fW)\times_{\U(\fH)}\fH$, where $\Fr(\fW)$ is the
principal $\U(\fH)$ bundle of unitary frames for $\fW$.  Recall that a {\em
quasi-subbundle\/} $V'\to M$ of a smooth vector bundle $\pi:V\to M$ over a
manifold $M$ is a closed subspace $V'\subset V$ such that
$\pi|_{V'}$ is still surjective and each fiber of $\pi|_{V'}$ is a
linear subspace \cite[Definition 1.2]{Karoubi}.

\begin{thm}
\label{thm:DefnOfStabilizeBundle}
Assume that the moduli space $M_{\fs}$ contains no zero-section pairs. Then
there is an open neighborhood $\sU$ of the subspace $\iota(M_{\fs})$ in
$\sC_{\ft}$, which does not contain any zero-section pairs or other
reducibles, and a finite-rank, smooth, trivial,
vector subbundle $\Xi\rightarrow\sU$ of $\fV \to \sU$, 
which is $S^1$ equivariant with respect to the circle action
\eqref{eq:VReducibleAction}, such that the following hold:
\begin{enumerate}
\item The restriction of $\Xi$ to $\iota(M_{\fs})$ is a
complex vector bundle.
\item The smooth bundle map $\Pi_{\Xi^\perp}:\fV\to\Xi^\perp$ defined by
the fiberwise $L^2$-orthogonal projection $\Pi_{\Xi^\perp}$
onto the subbundle $\Xi^\perp\to\sU$ restricts to a surjective fiber map
$\Pi_{\Xi^\perp}:\Ran(D\fS)_{A,\Phi}\to\Xi_{A,\Phi}^\perp$ for any point
$[A,\Phi]\in\sU$.
\item
If $(A,\Phi)$ is an $L^2_\ell$ representative of a point in $\sU$, for some
integer $k\le\ell\le\8$, then the fiber $\Xi_{A,\Phi}$ is contained in
$L^2_\ell(F_2) = L^2_\ell(\La^+\otimes\fg_{\ft})\oplus L^2_\ell(V^-)$.
\end{enumerate}
\end{thm}

\begin{proof}
By Corollary \ref{KuiperSmooth}, there is a smooth
trivialization of $\fV^n$ given by a smooth isomorphism of
complex Hilbert bundles,
$$
\btau: \fV^n\cong \iota(M_{\fs})\times L^2(F^n_2),
$$
as the space $\sG_{\fs}$ (of $L^2_{k+1}$ gauge transformations)
is a subgroup of the unitary group of the Hilbert space
$L^2(F^n_2)$. We have an equality of $L^2$-orthogonal complements,
$$
(\Ran D\fS_{A,\Phi})^\perp 
= (\Ran d^{1,n}_{A,\Phi})^\perp \subset L^2(F^n_2),
$$
for every $(A,\Phi)\in \iota(\tM_{\fs})$, by Lemma
\ref{lem:IdentityOfReducibleCohomology} and the fact that $\Ran
d^{1,t}_{A,\Phi} = L^2(F^t_2)$ by Proposition
\ref{prop:SmoothFamilyOfReducibles}, for generic parameters $\tau$.  Thus,
the family of operators $\{d^{1,n}_{A,\Phi}:(A,\Phi)\in \iota(\tM_{\fs})\}$
defines quasi-subbundles $\Ran\bd^{1,n}$ and $(\Ran\bd^{1,n})^\perp$ of
$\fV^n\to \iota(M_{\fs})$ with fibers
$$
\Ran \bd^{1,n}|_{[A,\Phi]}
=
\{(A',\Phi',\Ran d^{1,n}_{A',\Phi'}):(A',\Phi')\in [A,\Phi]\}/\sG_{\ft}.
$$
For each point $[A,\Phi]\in \iota(M_{\fs})$, we define
$$
\bV_{[A,\Phi]}^\perp
=
\btau(\Ran\bd^{1,n}|_{[A,\Phi]})
\subset L^2(F^n_2).
$$
Observe that there is an open neighborhood $U_{[A,\Phi]}\subset
\iota(M_{\fs})$ with the property that for all $[A',\Phi']\in U_{[A,\Phi]}$,
the map
$$
\btau(\Ran\bd^{1,n}|_{[A',\Phi']})\to
\bV_{[A,\Phi]}^\perp
$$
defined by $L^2$-orthogonal projection onto $\bV_{[A,\Phi]}^\perp$ is
surjective.  Since $\iota(M_{\fs})$ is compact, it has a finite subcover
$U_{[A_\alpha,\Phi_\alpha]}$ of such neighborhoods.  If $\bV=\oplus_\alpha
\bV_{[A_\alpha,\Phi_\alpha]}$, then $\bV$ is a finite-dimensional, complex
subspace of $L^2(F^n_2)$.  If we define $\Xi''\to \iota(M_{\fs})$ by setting
$$
\Xi''=\btau^{-1}\left(\iota(M_{\fs})\times \bV\right),
$$
then $\Xi$ will be a complex, finite-dimensional, trivial
subbundle of $\fV^n|_{\iota(M_{\fs})}$ such that fiberwise $L^2$-orthogonal
projection onto the subbundle $\Xi^\perp\subset \fV^n|_{\iota(M_{\fs})}$ is
surjective when restricted to the quasi-subbundle $\Ran\bd^{1,n}$.

We now extend the bundle $\Xi''\to \iota(M_{\fs})$ to a subbundle of $\fV$
over an open neighborhood $\sU$ of $\iota(M_{\fs})$ in $\sC^0_{\ft}$, which
does not contain any other reducible or zero-section pairs, and which is
$S^1$ equivariant with respect to the circle action
\eqref{eq:VReducibleAction}.  The space $M_{\fs}$ is a smooth submanifold
of the Riemannian manifold $\sC^0_{\ft}$ by Lemma
\ref{lem:ReducibleSubmanifold} and so it has an $S^1$-equivariant normal
bundle $\pi:\fN\to \iota(M_{\fs})$ and an $S^1$-invariant tubular
neighborhood given by an $S^1$-equivariant diffeomorphism
$$
\bg:\sO\subset\fN\to \sU\subset\sC^0_{\ft},
$$
from an open, $S^1$-invariant neighborhood $\sO$ of the zero
section of $\fN$ onto an open, $S^1$-invariant neighborhood of
$\iota(M_{\fs})$ in $\sC^0_{\ft}$ (see, for example,
\cite[p. 306]{Bredon}). The bundle projection $\pi$ and the
diffeomorphism $\bg$ define an $S^1$-equivariant, $C^\8$
retraction, 
$$
\br
=
\pi\circ\bg^{-1}:\sU\subset\sC^0_{\ft}\to \iota(M_{\fs}).
$$
The complex vector bundle $\Xi''\to \iota(M_{\fs})$ then extends to a
vector bundle $\br^*\Xi''\to\sU$, a subbundle of
$\br^*(\fV|_{\iota(M_{\fs})})\to\sU$, which are both $S^1$ equivariant with
respect to the circle action \eqref{eq:VReducibleAction}.  Since the map
$\br$ is an $S^1$-equivariant, $C^\8$ retraction, there is an
$S^1$-equivariant, $C^\8$ isomorphism \cite[Theorem 4.1.5]{Hirsch},
\cite{Husemoller},
$$
\boldf:\fV|_\sU \to \br^*(\fV|_{\iota(M_{\fs})}).
$$
We obtain a $C^\8$ subbundle of the vector bundle $\fV|_\sU$ by
setting
$$
\Xi' = \boldf^{-1}(\br^*\Xi'') \to \sU,
$$
both $S^1$ equivariant with respect to the action \eqref{eq:VReducibleAction}.
Because $\Xi''\to \iota(M_{\fs})$ is isomorphic to
$\iota(M_{\fs})\times\CC^{r_{\Xi}}$ as a complex vector bundle, for some
$r_\Xi\in\NN$, we obtain an isomorphism of $S^1$-equivariant vector bundles,
$$
\Xi'\cong \sU\times\CC^{r_{\Xi}},
$$
since the maps $\boldf$ and $\br$ are $S^1$ equivariant; the circle
acts non-trivially on $\sU$, except along the stratum $\iota(M_{\fs})$, and
acts by complex multiplication on $\CC^{r_{\Xi}}$.

By construction, the fiberwise $L^2$-orthogonal projection
$\fV|_{A,\Phi}\to (\Xi_{A,\Phi}')^\perp$ restricts to a surjective map
$\Ran d^1_{A,\Phi}\to(\Xi_{A,\Phi}')^\perp$ for any pair $(A,\Phi)$
representing a point in $\sU$, after shrinking $\sU$ if necessary.

Given an $L^2_k$ pair $(A,\Phi)$ representing a point in
$\sU\subset\sC^0_{\ft}$, our construction yields a subspace
$$
\Xi_{A,\Phi}' \subset L^2(F_2)
=
L^2(\Lambda^+\otimes\fg_{\ft}) \oplus L^2(V^-).
$$
If $\ell\ge k$ is an integer, it does not necessarily follow that
$\Xi_{A,\Phi}'$ is contained in the subspace of $L^2_{\ell-1}$ pairs when
$(A,\Phi)$ is an $L^2_\ell$ pair. However, for any $t>0$ the heat operator
$$
\exp(-t(1+d_{A,\Phi}^1d_{A,\Phi}^{1,*})):
L^2(\Lambda^+\otimes\fg_{\ft}) \oplus L^2(V^-)
\to
L^2_{\ell-1}(\Lambda^+\otimes\fg_{\ft}) \oplus L^2_{\ell-1}(V^-)
$$
is a bounded, $\sG_{\ft}$-equivariant, $S^1$-equivariant, linear map and
we can define
$$
\Xi_{A,\Phi}
=
\exp(-t(1+d_{A,\Phi}^1d_{A,\Phi}^{1,*}))\Xi_{A,\Phi}'.
$$
For small enough $t = t(\sU)$, the approximation properties of the heat
kernel (see \cite[Lemma A.1]{FL1} for a similar application),
ensure that
\begin{itemize}
\item
$\Xi\to\sU$ is a trivial, $S^1$-equivariant  vector bundle,
with the same rank as $\Xi'$, and
\item
$L^2$-orthogonal projection $\fV|_{A,\Phi}\to\Xi_{A,\Phi}^\perp$
restricts to a surjective map $\Ran d^1_{A,\Phi}\to\Xi_{A,\Phi}^\perp$ for any
pair $(A,\Phi)$ representing a point in $\sU$.
\end{itemize}
We let $\Pi_\Xi:\fV\to\Xi$ be the $C^\8$ bundle map defined by
fiberwise $L^2$-orthogonal projection. This completes the proof.
\end{proof}

\subsubsection{Definition of the thickened moduli space and the link}
With Theorem \ref{thm:DefnOfStabilizeBundle} in place, we can finally
construct the required ambient, finite-dimensional, smooth manifold of
Definition \ref{defn:AbstractLink} and the link of the singular stratum of
reducibles.

\begin{defn}
\label{defn:ThickenedModuliSpace}
Assume $M_{\fs}$ contains no zero-section pairs.  Let $\Xi$ be a
finite-rank, smooth, $S^1$-equivariant vector bundle over an open
neighborhood of $\iota(M_{\fs})$ in $\sC_{\ft}$ such that, as
in Theorem \ref{thm:DefnOfStabilizeBundle}, $L^2$-orthogonal projection
gives a surjective map of quasi vector bundles $\Ran D\fS\to\Xi^\perp$ over
$M_{\fs}$. We say that $\Xi$ is a {\em stabilizing bundle\/} for $\Ran D\fS$
and call 
\begin{equation}
\label{eq:DefineThickenedModuli}
\sM_{\ft}(\Xi,\fs) = (\Pi_{\Xi^\perp}\fS)^{-1}(0) \subset\sC_{\ft}
\end{equation}
the {\em thickened moduli space\/} defined by $\Xi$.
If $(A,\Phi)\in\tsC_{\ft}$ represents $[A,\Phi]\in\sC_{\ft}$,
we write $\Xi|_{(A,\Phi)}$ for the subspace of $L^2_{k-1}(F_2)$
representing the fiber of $\Xi$ over $[A,\Phi]$.  We then define 
\begin{equation}
N_{\ft}(\Xi,\fs) =\Ker(\Pi_{\Xi^\perp}\bsD^n),
\label{eq:NormalBundle}
\end{equation}
a complex, finite-rank, smooth vector bundle over $M_{\fs}$ with fibers
$$
N_{\ft}(\Xi,\fs)|_{[A,\Phi]}
\cong
\Ker\left(d_{A,\Phi}^{0,n,*} + \Pi_{\Xi^\perp} d_{A,\Phi}^{1,n}\right),
$$
noting that $\sD^n_{A,\Phi} = d_{A,\Phi}^{0,n,*} + d_{A,\Phi}^{1,n}$ for
$[A,\Phi]\in \iota(M_{\fs})$.
\end{defn}

\begin{thm}
\label{thm:ThickenedModuliSpace}
Suppose that the \spinu structure $\ft$ admits a reduction
$\ft=\fs\oplus \fs\otimes L$.
Assume $M_{\fs}$ contains no zero-section pairs.  Then the
following hold:
\begin{enumerate}
\item There is an $S^1$-invariant, open neighborhood $\sU$ of 
  $\iota(M_{\fs})$ in $\sC_{\ft}$ such that the zero locus $\sU\cap
  \sM_{\ft}(\Xi,\fs)$ is regular and so a manifold of dimension $\dim
  \sM_{\ft}+\rank_\RR\Xi$.
\item The space $M_{\fs}$ is a smooth, $S^1$-invariant submanifold of
  $\sM_{\ft}(\Xi,\fs)$.
\item The bundle $N_{\ft}(\Xi,\fs)$ is a normal
  bundle for the submanifold $\iota:M_{\fs}\embed \sM_{\ft}(\Xi,\fs)$ and
  the tubular map is equivariant with respect to the circle action on
  $N_{\ft}(\Xi,\fs)$ given by the trivial action on the base $M_{\fs}$ and
  complex multiplication on the fibers, and the circle action on
  $\sM_{\ft}(\Xi,\fs)$ induced from the $S^1$ action
  \eqref{eq:VReducibleAction}. 
\item The restriction of the section $\fS$ to $\sM_{\ft}(\Xi,\fs)$
  takes values in $\Xi$ and vanishes transversely on
  $\sM_{\ft}(\Xi,\fs) - \iota(M_{\fs})$.
\end{enumerate}
\end{thm}

\begin{proof}
The space $\sM_{\ft}(\Xi,\fs)$ is the zero locus of the section $\fT =
\Pi_{\Xi^\perp}\fS$ of a vector subbundle $\Xi^\perp \to \sU \subset
\sC_{\ft}$ of $\fV|_\sU$ constructed in
Theorem \ref{thm:DefnOfStabilizeBundle},
for some open neighborhood $\sU$ of $\iota(M_{\fs})$ in $\sC_{\ft}$.
For any point $[A,\Phi] \in
\iota(M_{\fs})\subset\fS^{-1}(0)$, we have
$$
(D\fT)_{A,\Phi} = \Pi_{\Xi^\perp} (D\fS)_{A,\Phi}.
$$
(The differential of $\Pi_\Xi$ does not appear here since $\fS[A,\Phi] =
0$.) According to Definition \ref{defn:ThickenedModuliSpace}, the
$L^2$-orthogonal projection $\Pi_{\Xi^\perp}$ gives a surjective map
$$
\Pi_{\Xi^\perp}: \Ran(D\fS)_{A,\Phi} \to \Xi_{A,\Phi}^\perp
$$
and thus $(D\fT)_{A,\Phi}$ is surjective at all points $[A,\Phi]$ in the
image of $M_{\fs}$.  Surjectivity is an open condition, so we may
assume that $D\fT$ is surjective on the open neighborhood $\sU$ of
$\iota(M_{\fs})\subset\sC_{\ft}$, after shrinking $\sU$ if necessary.
The zero locus of $\fT$ in this open set is regular and thus a smooth
submanifold of $\sC_{\ft}$, which gives Assertion (1).

By Lemma \ref{lem:ReducibleSubmanifold} the space $M_{\fs}$ is a
smooth submanifold of $\sC_{\ft}$ and as its image is contained in
$\sM_{\ft}(\Xi,\fs)$, it is also a smooth submanifold of $\sM_{\ft}(\Xi,\fs)$
and Assertion (2) follows, as the zero locus of an $S^1$-equivariant
section is $S^1$ invariant. 

We observe that $N_{\ft}(\Xi,\fs)$ is the normal bundle of
$M_{\fs}$ in $\sM_{\ft}(\Xi,\fs)$, since
\begin{equation}
\label{eq:TangThickenedModSpaceSplit}
\begin{aligned}
T_{A,\Phi}\sM_{\ft}(\Xi,\fs)
&=
\Ker d_{A,\Phi}^{0,*}\cap\Ker\Pi_{\Xi^\perp}(D\fS)_{A,\Phi}
\\
&=
\left(\Ker d_{A,\Phi}^{0,t,*}\oplus \Ker d_{A,\Phi}^{0,n,*}\right)
\cap\left(\Ker d_{A,\Phi}^{1,t}\oplus
\Ker\Pi_{\Xi^\perp} d_{A,\Phi}^{1,n}\right)
\\
&=
\Ker\sD^t_{A,\Phi}\oplus\Ker\Pi_{\Xi^\perp}\sD^n_{A,\Phi}
\\
&=
T_{A,\Phi}\iota(M_{\fs})\oplus N_{\ft}(\Xi,\fs)|_{A,\Phi}.
\end{aligned}
\end{equation}
In the second equality above, we make use of the fact that the fibers of
the vector bundle $\Xi\to \iota(M_{\fs})$ are contained in
$L^2_{k-1}(F_2^n)$ and so (by definition) $\Pi_\Xi = 0$ on
$L^2_{k-1}(F_0^t)\oplus L^2_{k-1}(F_2^t)\supset\Ran
\sD^t_{A,\Phi}$. Also, the splitting
\eqref{eq:TangThickenedModSpaceSplit} of tangent spaces corresponds to the
splitting of Hilbert spaces $L^2_k(F_1)\cong L^2_k(F_1^t)\oplus
L^2_k(F_1^n)$, defined by the subspaces
\eqref{eq:HilbertSpaceRedCplxSplitting}. Hence, the isomorphism 
\eqref{eq:TangThickenedModSpaceSplit} is equivariant with respect to the
circle action \eqref{eq:VReducibleAction} on the subspace
$T_{A,\Phi}\sM_{\ft}(\Xi,\fs)$ of $L^2_k(F_1)$,
the trivial action on the subspace $T_{A,\Phi}\iota(M_{\fs})$ of
$L^2_k(F_1^t)$, and complex multiplication on the subspace
$N_{\ft}(\Xi,\fs)|_{A,\Phi}$ of $L^2_k(F_1^n)$.
This proves Assertion (3).

Because $\sM_{\ft}(\Xi,\fs)$ is given by
$\sU\cap(\Pi_{\Xi^\perp}\fS)^{-1}(0)=\sU\cap\fS^{-1}(\Xi)$, we see that
$\fS$ takes values in $\Xi$ on $\sM_{\ft}(\Xi,\fs)$ and so
$\Pi_{\Xi^\perp}\fS = \fS$ on $\sM_{\ft}(\Xi,\fs)$.  According to Theorem
\ref{thm:Transversality}, our transversality result for the section $\fS$
of $\fV|_\sU = \Xi^\perp\oplus\Xi$ defined by the $\PU(2)$ monopole
equations \eqref{eq:PT}, the section $\fS$ vanishes transversely on $\sU -
\iota(M_{\fs})$ with zero locus $(\sU\cap \sM_{\ft})- \iota(M_{\fs})$.
This implies that for each $[A,\Phi]$ in $(\sU -
\iota(M_{\fs}))\cap\fS^{-1}(0)$, the differential
$$
(D\fS)_{A,\Phi}:T_{[A,\Phi]}\sC^{*,0}_{\ft}\to \fV|_{[A,\Phi]}
$$
is surjective. But we have the identifications
\begin{align*}
T_{[A,\Phi]}\sM_{\ft}(\Xi,\fs)
&=
\Ker (D(\Pi_{\Xi^\perp}\fS))_{[A,\Phi]}
\\
&=
\Ker (\Pi_{\Xi^\perp}D\fS)_{[A,\Phi]}
\quad\text{(since $\fS(A,\Phi)=0$)}
\\
&=
(D\fS)_{[A,\Phi]}^{-1}(\Xi|_{[A,\Phi]}),
\end{align*}
and so one has a surjective differential
$$
(D\fS)_{[A,\Phi]}: T_{[A,\Phi]}\sM_{\ft}(\Xi,\fs) \to \Xi|_{[A,\Phi]},
$$
for $[A,\Phi] \in (\sU-\iota(M_\fs))\cap\fS^{-1}(0)$.  Thus, the sections
$\Pi_{\Xi^\perp}\fS$ and $\fS$---which are equal on
$\sM_{\ft}(\Xi,\fs)$---vanish transversely when restricted to
$\sM_{\ft}(\Xi,\fs)$, proving Assertion (4).
\end{proof}

The equivariant tubular neighborhood theorem (see \cite[p. 306]{Bredon} for
the finite-dimensional, $G$-equivariant case and \cite{Lang} for the
infinite-dimensional case) provides an embedding
\begin{equation}
\bga:\sO\hookrightarrow\sC_{\ft},
\label{eq:NormalBundleNbhdEmbedding}
\end{equation}
mapping an open, $S^1$-invariant 
neighborhood $\sO$ of the zero-section $M_{\fs}\subset
N_{\ft}(\Xi,\fs)$ onto an open neighborhood of the submanifold
$\iota:M_{\fs}\embed \sM_{\ft}(\Xi,\fs)$ which covers the embedding
$\iota$; the map \eqref{eq:NormalBundleNbhdEmbedding} is $S^1$-equivariant
with respect to scalar multiplication on the fibers of $N_{\ft}(\Xi,\fs)$
and the circle action induced by \eqref{eq:VReducibleAction} on
$\sC_{\ft}$. 

The map $\bga$ then descends to a homeomorphism, and a diffeomorphism on
smooth strata, from the zero locus $\bvarphi^{-1}(0)/S^1$ in
$N_{\ft}(\Xi,\fs)/S^1$ onto an open neighborhood of $\iota(M_{\fs})$
in the actual moduli space, $\sM_{\ft}$, where we define
\begin{equation}
\bvarphi = \Pi_\Xi\fS\circ\bga,
\label{eq:ObstructionSection}
\end{equation}
to be an $S^1$-equivariant {\em obstruction section\/}
over $\sO\subset N_{\ft}(\Xi,\fs)$ of the $S^1$ 
equivariant {\em obstruction bundle\/}
\begin{equation}
\bga^*\Xi \to N_{\ft}(\Xi,\fs).
\label{eq:ObstructionBundle}
\end{equation}
This descends to a vector bundle 
$$
(\bga^*\Xi)/S^1 \to N_{\ft}(\Xi,\fs)/S^1
$$
on the complement of the
zero section, $M_{\fs}\subset N_{\ft}(\Xi,\fs)/S^1$, whose
Euler class may be computed from
$$
(\bga^*\Xi)/S^1 \cong
(\pi_N^*\Xi)/S^1 \to N_{\ft}(\Xi,\fs)/S^1,
$$
where $\pi_N:N_{\ft}(\Xi,\fs)\rightarrow M_{\fs}$ is the
projection.

We can now construct the link of $M_{\fs}$ in $\sM_{\ft}$,
following Definition \ref{defn:AbstractLink}.

\begin{defn}
\label{defn:LinkOfReducible}
Assume $M_{\fs}$ contains no zero-section pairs.
Let $N^\eps_{\ft}(\Xi,\fs)\subset N_{\ft}(\Xi,\fs)$ be the sphere bundle of
fiber vectors of length $\eps$ and set 
$$
\PP N_{\ft}(\Xi,\fs) = N^\eps_{\ft}(\Xi,\fs)/S^1.
$$
The {\em link of the stratum $\iota(M_{\fs})\subset \sM_{\ft}$\/}
of reducible $\PU(2)$ monopoles is defined by
\begin{equation}
\label{eq:DefineReducibleLink}
\bL_{\ft,\fs}
=
\bga\left(\bvarphi^{-1}(0)\cap \PP N_{\ft}(\Xi,\fs)\right)
\subset
\sM^{*,0}_{\ft}/S^1.
\end{equation}
The orientation for $\bL_{\ft,\fs}$ is defined by the orientation on
$M_{\fs}$ and in turn from the homology orientation $\Omega$ (see
\cite[\S 6.6]{MorganSWNotes}) and the complex structure on the fibers of
$N_{\ft}(\Xi,\fs)$ given by the $S^1$ action. We treat the quotient $\PP
N_{\ft}(\Xi,\fs)$ as the complex projectivization of a complex vector
bundle; we use the complex orientation on the obstruction bundle
$\bga^*\Xi$ in \eqref{eq:ObstructionBundle}.
\end{defn}

The following lemma shows that the link $\bL_{\ft,\fs}$ can be
represented homologically.

\begin{lem}
\label{lem:TransverseObstruction}
Assume $M_{\fs}$ contains no zero-section pairs.  Then for
generic $\eps$, the section $\bvarphi$ of the obstruction bundle
$\bga^*\Xi$ in \eqref{eq:ObstructionBundle} vanishes transversely on
$N^\eps_{\ft}(\Xi,\fs)$.
\end{lem}

\begin{proof}
{}From the final statement of Theorem \ref{thm:ThickenedModuliSpace}, we see
that $\bvarphi$ vanishes transversely on $\sM_{\ft}(\Xi,\fs) -
\iota(M_{\fs})$ and so $\bvarphi$ cuts out the zero locus, $\sM_{\ft} -
\iota(M_{\fs})$, as a regular submanifold of $\sM_{\ft}(\Xi,\fs) -
\iota(M_{\fs})$. Then for generic values of $\eps$, the zero locus will
intersect $N^\eps_{\ft}(\Xi,\fs)$ transversely.
\end{proof}

Lemma \ref{lem:TransverseObstruction} implies that 
\begin{equation}
\label{eq:DefineHomologyOfReducibleLink}
[\bL_{\ft,\fs}]
=
e\left((\bga^*\Xi)/S^1\right)\cap [\PP N_{\ft}(\Xi,\fs)]
\end{equation}
is the homology class of the link in Definition \ref{defn:LinkOfReducible}.

\subsubsection{Group actions and lifts of the normal bundle embedding to
the pre-configuration space}
\label{subsubsec:GroupActionsNormalBundleLifts}
In \cite{FL2b} we shall need a lift of the $S^1$-equivariant diffeomorphism
$\bgamma$ from a neighborhood $\sO\subset N_{\ft}(\Xi,\fs)$ of the
zero-section $M_{\fs}$ onto on open neighborhood of $\iota(M_{\fs})$ in the
thickened moduli space $\sM_{\ft}(\Xi,\fs)\subset\sC_{\ft}$,
\begin{equation}
\label{eq:PreConfigNormalBundleEmbedding}
\tilde\bgamma:\tilde\sO\subset\tN_{\ft}(\Xi,\fs) 
\to \tilde\sM_{\ft}(\Xi,\fs),
\end{equation}
where $\tilde\sO\subset\tN_{\ft}(\Xi,\fs)$ is an $S^1$ and $\sG_{\fs}$
invariant open neighborhood of $\tM_{\fs}$ and $\tilde\sM_{\ft}(\Xi,\fs)
= \pi^{-1}(\sM_{\ft}(\Xi,\fs))$. It is convenient to describe the
construction here.  As usual, we need only consider the case where
$M_{\fs}$ contains no zero-section pairs.

To see what should be the ``correct'' equivariance properties of the lift
\eqref{eq:PreConfigNormalBundleEmbedding}, we first
consider the obvious extension $\iota:N_{\ft}(\Xi,\fs)\to\sC_{\ft}$ of the
embedding $\iota:M_{\fs}\embed\sC_{\ft}$, since $\bgamma$ approximates this
extension on a small open neighborhood of the zero section, $M_{\fs}$.  Let
$\tM_{\fs}=\pi^{-1}(M_{\fs})\subset\tsC_{\fs}$ be the preimage of
$M_{\fs}\subset\sC_{\fs}$ under the projection $\pi:\tsC_{\fs}\to
\tsC_{\fs}/\sG_{\fs}$, and let
$$
\tN_{\ft}(\Xi,\fs) = \pi^*N_{\ft}(\Xi,\fs) \to \tM_{\fs}
$$
be the $\sG_{\fs}$-equivariant pullback bundle, so
$$
\tN_{\ft}(\Xi,\fs)
\subset
\tM_{\fs} \times F_1^n,
$$
where $F_1^n\subset F_1\cong F_1^t\oplus F_1^n$ is the complex Hilbert
subspace in \eqref{eq:HilbertSpaceRedCplxSplitting}. The bundle
$\tN_{\ft}(\Xi,\fs)$ is complex, with trivial circle action on $\tM_{\fs}$
and action by complex multiplication on $F_1^n$.

The $S^1$-equivariant map $\iota:N_{\fs}(\Xi,\fs)\to\sC_{\ft}$,
$[B,\Psi,\eta] \mapsto [\iota(B,\Psi)+\eta]$, is covered by
\begin{equation}
\label{eq:LiftSWNormalBundleInclusion}
\iota:\tN_{\ft}(\Xi,\fs) \to \tsC_{\ft}, 
\quad
(B,\Psi,\eta) \mapsto \iota(B,\Psi)+\eta.
\end{equation}
This map is equivariant with respect to the embedding
$\varrho:\sG_{\fs}\embed\sG_{\ft}$ in definition \eqref{defn:GInclusion},
for the following domain and range $\sG_{\fs}$ actions:
\begin{itemize}
\item
The action of $s\in\sG_{\fs}$ on $\tN_{\ft}(\Xi,\fs)$ given by the usual
gauge group action on $(B,\Psi)\in\tM_{\fs}$ and the action of $\sG_{\fs}$ on
$\eta=(\beta,\psi')\in F_1^n$ induced by the isomorphism $F_1\cong
F_1^t\oplus F_1^n$ and the embedding $\varrho:\sG_{\fs}\embed\sG_{\ft}$, as
described in \S \ref{subsubsec:DecompGroupActions}, so
$$
(s,(B,\Psi,\beta,\psi'))\mapsto (s(B,\Psi),s^{-2}\beta,s^{-1}\psi').
$$
\item
The action of $\varrho(s)\in\sG_{\ft}$ on $\tsC_{\ft}$ and $F_1$ in 
definition \eqref{defn:GInclusion}, so
$$
(s,\iota(B,\Psi)+\eta)
\mapsto \varrho(s)(\iota(B,\Psi) + \eta).
$$
The fact that $\varrho(s)\iota(B,\Psi) = \iota(s(B,\Psi))$ was noted in
Lemma \ref{lem:TopU1Embedding}, while we see that $\iota(s^{-2}\beta,
s^{-1}\psi') = \varrho(s)\iota(\beta,\psi')$ by the remarks following
equation \eqref{eq:su(E)SplitGaugeGroupAction} in
\S \ref{subsubsec:DecompGroupActions}.
\end{itemize}

The map \eqref{eq:LiftSWNormalBundleInclusion} is also $S^1$ equivariant
for the following domain and range circle actions:
\begin{itemize}
\item
The trivial $S^1$ action on $\tM_{\fs}$ and action by complex
multiplication on $F_1^n$,
$$
(e^{i\theta},(B,\Psi,\eta)) \mapsto (B,\Psi,e^{i\theta}\eta).
$$
\item
The $S^1$ action \eqref{eq:VReducibleAction} on $\tsC_{\ft}$, induced by
the trivial action on the factor $W$ of $V=W\oplus W\otimes L$ and complex
multiplication on $W\otimes L$, so on the image
\eqref{eq:LiftSWNormalBundleInclusion} one has
$$
(e^{i\theta},\iota(B,\Psi)+\eta) \mapsto \iota(B,\Psi)
+e^{i\theta}\eta,
$$
recalling that the points $\iota(B,\Psi)$ are fixed by the action
\eqref{eq:VReducibleAction}. 
\end{itemize}

We now turn to the construction of the map
\eqref{eq:PreConfigNormalBundleEmbedding}, which we shall require to have
the same $\sG_{\fs}$ and $S^1$ equivariance properties as the map
\eqref{eq:LiftSWNormalBundleInclusion}.   
If $[B,\Psi]\in M_{\fs}$, then Proposition 2.8 in
\cite{FL1} yields an open neighborhood $\tU_{\iota(B,\Psi)}$ of
$\iota(B,\Psi)$ in the slice $\iota(B,\Psi)+\Ker d_{\iota(B,\Psi)}^{0,*}$,
such that projection onto $\pi(\tU_{\iota(B,\Psi)}) = U_{\iota(B,\Psi)}$
gives a local parameterization for an open neighborhood of
$[\iota(B,\Psi)]$ in $\sM_{\ft}(\Xi,\fs)$; the pair $\iota(B,\Psi)$ has
trivial stabilizer in $\sG_{\ft}$ since $\Psi\not\equiv 0$.  Because
$M_{\fs}$ is compact we can assume (by shrinking $\sO$ if necessary) that
$$
\bgamma[B,\Psi,\eta] \in U_{[\iota(B,\Psi)]},
$$ 
for all $[B,\Psi,\eta]\in \sO\subset N_{\ft}(\Xi,\fs)$. Hence, 
$$
\bgamma[B,\Psi,\eta] = [\iota(B,\Psi) + \bgamma_0(B,\Psi,\eta)]
$$ 
for a slice element $\bgamma_0(B,\Psi,\eta)\in \Ker
d_{\iota(B,\Psi)}^{0,*}$ uniquely determined by $\iota(B,\Psi,\eta)$, and
we can therefore set
\begin{equation}
\label{eq:DefnLiftedNormalEmbedding}
\tilde\bgamma(B,\Psi,\eta) 
= 
\iota(B,\Psi) + \bgamma_0(B,\Psi,\eta),
\quad
(B,\Psi,\eta)\in \tilde\sO\subset \tN_{\ft}(\Xi,\fs),
\end{equation}
where $\tilde\sO$ is the preimage of $\sO=\tilde\sO/\sG_{\fs}$.
Note that $\bgamma_0(B,\Psi,\eta) \approx \eta$. If $s\in\sG_{\fs}$, then 
the preceding equation yields
\begin{equation}
\label{eq:GtLiftedNormalEmbedding}
\varrho(s)\tilde\bgamma(B,\Psi,\eta) 
= 
\iota(s(B,\Psi)) + \varrho(s)\bgamma_0(B,\Psi,\eta),
\end{equation}
and we again have $\varrho(s)\bgamma_0(B,\Psi,\eta) \in \Ker
d_{\iota(s(B,\Psi))}^{0,*}$, by the $\sG_{\ft}$-equivariance of the slice
condition and the fact that $\varrho(s)\iota(B,\Psi)=\iota(s(B,\Psi))$. On
the other hand,
\begin{equation}
\label{eq:GsLiftedNormalEmbedding}
\tilde\bgamma(s(B,\Psi,\eta)) 
= 
\iota(s(B,\Psi)) + \bgamma_0(s(B,\Psi,\eta)),
\end{equation}
while $\bgamma[s(B,\Psi,\eta)] = \bgamma[B,\Psi,\eta] \in
\sO=\tilde\sO/\sG_{\fs}$. Hence, comparing equations
\eqref{eq:GtLiftedNormalEmbedding} and \eqref{eq:GsLiftedNormalEmbedding},
we see that we must have
$$
\tilde\bgamma(s(B,\Psi,\eta)) 
=
\varrho(s)\tilde\bgamma(B,\Psi,\eta).
$$
Therefore, the map $\tilde\bgamma$ has the same $\sG_{\fs}$-equivariance
properties as the map $\iota:\tN_{\ft}(\Xi,\fs)\to\tsC_{\ft}$ in
\eqref{eq:LiftSWNormalBundleInclusion}; a very similar argument shows that
it has the same $S^1$-equivariance properties.

\subsection{The link of a stratum of reducible monopoles: Chern character
of the normal bundle}
\label{subsec:TopOfNormal}
To compute intersection pairings on the space $\PP N_{\ft}(\Xi,\fs)$, we
need to know the Chern classes of $N_{\ft}(\Xi,\fs)$. The Chern character
of the vector bundle $N_{\ft}(\Xi,\fs)\to M_{\fs}$ is computed
by observing that, as elements of $K(M_{\fs})$, we have
$[N_{\ft}(\Xi,\fs)] = \ind\bsD^n + [\Xi]$, where $\bsD^n$ is the normal
component \eqref{eq:NormalReducDefOperator} of the family of deformation
operators $\bsD$ parameterized by $M_{\fs}$.  For convenience in
this section, we shall often omit explicit mention of the embedding map
$\iota:M_{\fs}\embed \sM_{\ft}(\Xi,\fs)$ and write the bundles
$\iota^*N_{\ft}(\Xi,\fs)$ and $\iota^*(\ind\bsD^n)$ over $M_{\fs}$
simply as $N_{\ft}(\Xi,\fs)$ and $\ind\bsD^n$, respectively. We then use
the Atiyah-Singer index theorem for families to express
$\ch(N_{\ft}(\Xi,\fs))$ in terms of the cohomology classes on
$H^\bullet(\sC^0_{\fs};\RR)$ described in \S \ref{subsubsec:AbelianCohom}.

{}From \eqref{eq:NormalReducDefOperator} we obtain a family of elliptic
differential operators
%
%
$$
\begin{CD}
C^\8(\bE)     @>{\bsD^n}>>             C^\8(\bF)   \\
@VVV                         @VVV         \\
M_{\fs}\times X   @>{\id}>> M_{\fs}\times X
\end{CD}
$$
where $C^\8(\,\cdot\,)$ denotes the space of smooth sections of the
families of finite-rank vector bundles $\bE$, $\bF$ over
$M_{\fs}\times X$ defined by
\begin{equation}
\label{eq:UniversalBundle}
\begin{aligned}
\bE &= \tM_{\fs}\times_{\sG_{\fs}}F_1^n
 \\
&= \tM_{\fs}\times_{\sG_{\fs}}\left(
(\Lambda^1\otimes L)\oplus W^+\otimes L\right),
 \\
\bF &= \tM_{\fs}\times_{\sG_{\fs}}(F_0^n\oplus F_2^n)
\\
&= \tM_{\fs}\times_{\sG_{\fs}}
\left((\Lambda^0\oplus\Lambda^+)\otimes L\oplus
W^-\otimes L\right).
\end{aligned}
\end{equation}
Recall from the paragraph following \eqref{eq:su(E)SplitGaugeGroupAction}
that an element $s\in\sG_{\fs}$ acts on the bundles listed in
\eqref{eq:UniversalBundle} as multiplication by $s^{-2}$ on the factors
$\La^j\otimes L$ and as multiplication by $s^{-1}$ on the factors
$W^\pm\otimes L$.

Recall that $ N_{\ft}(\Xi,\fs) = \Ker\Pi_{\Xi^\perp}\bsD^n$ by
definition \eqref{eq:NormalBundle}. On the other hand, by
equation \eqref{eq:NormalReducDefOperator}, we see that
$$
\Coker\Pi_{\Xi^\perp}\sD^n_{\iota(B,\Psi)}
=
\Coker d_{\iota(B,\Psi)}^{0,n,*}\oplus
\Coker\Pi_{\Xi^\perp} d_{\iota(B,\Psi)}^{1,n}.
$$
For any $[B,\Psi]$, we have $\Coker d_{\iota(B,\Psi)}^{0,*} = 0$,
as the stabilizer of $\iota(B,\Psi)$ is {\em trivial in $\sG_{\ft}$\/},
so $\Coker d_{\iota(B,\Psi)}^{0,n,*} = 0$,
while Lemma \ref{lem:IdentityOfReducibleCohomology}
and Proposition \ref{prop:SmoothFamilyOfReducibles} implies that
$\Coker d_{\iota(B,\Psi)}^{1,t} = 0$. Hence,
$$
\Coker\Pi_{\Xi^\perp} d_{\iota(B,\Psi)}^{1,t} = 0
$$ 
and
\begin{align*}
\Coker\Pi_{\Xi^\perp}\sD^n_{\iota(B,\Psi)}
&=
\Coker\Pi_{\Xi^\perp} d_{\iota(B,\Psi)}^{1,n}\oplus
\Coker\Pi_{\Xi^\perp} d_{\iota(B,\Psi)}^{1,t}
\\
&=
\Coker\Pi_{\Xi^\perp} d_{\iota(B,\Psi)}^1
=
\Coker\Pi_{\Xi^\perp}(D\fS)_{\iota(B,\Psi)}.
\end{align*}
Because $L^2$-orthogonal projection from
$\Ran(D\fS)_{\iota(B,\Psi)}$ surjects onto $\Xi_{\iota(B,\Psi)}^\perp$
for all points $[B,\Psi]\in
M_{\fs}$, we obtain the identity
$$
\Coker\Pi_{\Xi^\perp}\bsD^n = \Coker\Pi_{\Xi^\perp} \cong \Xi.
$$
The subbundle $\Xi$ of $\fV \to M_{\fs}$ is trivial by the
construction of Theorem \ref{thm:DefnOfStabilizeBundle}. From the
stabilization construction of \cite[\S I.7.B]{BoossBleecker} and the preceding
remarks we see that, as elements of the $K$-theory group $K(M_{\fs})$,
\begin{equation}
\label{eq:NormalBundleIsNormalStab}
\begin{aligned}
\ind\bsD^n
&= [\Ker\Pi_{\Xi^\perp}\sD^n]-[\Coker\Pi_{\Xi^\perp}\sD^n]
\\
&= [N_{\ft}(\Xi,\fs)]-[\Xi]
\\
&= [N_{\ft}(\Xi,\fs)]-[M_{\fs}\times\CC^{r_\Xi}],
\end{aligned}
\end{equation}
where $r_\Xi = \rank_\CC\Xi$.
To compute $\ch(\ind\bsD^n)$ we shall apply the Atiyah-Singer index theorem for
families of Dirac operators.

\begin{prop}
\label{prop:AtiyahSinger}
\cite{AS4}, \cite[Theorem 5.1.16]{DK}
Let $X$ be a closed, oriented, smooth four-manifold with \spinc structure
$(\rho,W)$. Suppose $\sE\rightarrow T\times X$ is a locally
trivial family of complex
vector bundles over $X$, parameterized by
a compact space $T$,
with a connection $A_t$ on the bundle
$\sE_t = \sE|_{\{t\}\times X}$ for all $t\in T$.
Then the Chern character of the index bundle of the family of
Dirac operators parameterized by $T$,
%
%
$$
\begin{CD}
 C^\8(\sE\otimes W^+)   @>{\bD}>>         C^\8(\sE\otimes W^-)    \\
@VVV                         @VVV         \\
T\times X   @>{\id}>> T\times X
\end{CD}
$$
given by $D_{A_t}: C^\8(\sE_t\otimes W^+)\to
 C^\8(\sE_t\otimes W^-)$, $t\in T$, is
\begin{align*}
\ch(\ind (\bD,\sE,W))
&=
\ch(\sE)\ch(\ind(D,W))
\\
&=
\ch(\sE)e^{\half c_1(W^+)}\hat\sA(X)/[X],
\end{align*}
where $D: C^\8(W^+)\to C^\8(W^-)$ is the Dirac operator defined by
the given \spinc structure.
\end{prop}

The index of the family of operators in Proposition
\ref{prop:AtiyahSinger} defines a group homomorphism \cite[p. 184]{DK},
$$
\ind(\bD,\cdot\ ): K(T\times X) \to K(T),
$$
by taking the element $[\sE]$ of $K(T\times X)$ to $\ind(\bD,\sE,W) =
[\Ker\bD]-[\Coker\bD]$ in $K(T)$, a virtual vector bundle over the
parameter space $T$.  We now compute $\ch(N_{\ft}(\Xi,\fs))$ in the
following steps:
\begin{itemize}
\item Identify $\bsD^n$ with the sum of a pair of families of Dirac
  operators, $\bD'$ and $\bD''$ .
\item Compute $\ch(N_{\ft}(\Xi,\fs))$, using Proposition
  \ref{prop:AtiyahSinger} to compute the Chern characters of the index
  bundles of these families of Dirac operators.
\end{itemize}
The first step is accomplished in part by introducing the following
families of operators,
%
%
\begin{equation}
\label{eq:DeltaAndDiracL2Families}
\begin{CD}
 C^\8(\bE') @>{\bdelta}>>  C^\8(\bF') \\
@VVV                         @VVV         \\
M_{\fs}\times X   @>{\id}>> M_{\fs}\times X
\end{CD}
\qquad\text{and}\qquad
\begin{CD}
 C^\8(\bE'')   @>{\bD''}>>    C^\8(\bF'')    \\
@VVV                         @VVV         \\
M_{\fs}\times X   @>{\id}>> M_{\fs}\times X
\end{CD}
\end{equation}
where $\bE = \bE'\oplus \bE''$ and $\bF = \bF'\oplus \bF''$,
and we have defined
\begin{equation}
\label{eq:LambdaWSplittingUnivBundle}
\begin{matrix}
\bE' = \tM_{\fs}\times_{\sG_{\fs}}
                    (\Lambda^1\otimes L)\hfill
&\text{and}
&\bE'' = \tM_{\fs}\times_{\sG_{\fs}} (W^+\otimes L),\hfill
\\
\bF' = \tM_{\fs}\times_{\sG_{\fs}}
                 ((\Lambda^0\oplus\Lambda^+)\otimes L)\hfill
&\text{and}
&\bF''  = \tM_{\fs}\times_{\sG_{\fs}} (W^-\otimes L).\hfill
\end{matrix}
\end{equation}
Explicitly, if $(B,\Psi)$ is a pair representing a point $[B,\Psi]$ in
$M_{\fs}$, so the spin connection in the pair
$\iota(B,\Psi)\in\tilde\sM_\ft$ is given by $B\oplus B\otimes A_L$ and
$A_L=A_\Lambda\otimes (B^{\det})^*$ is the induced unitary connection on
the subbundle $L$ of $\fg_{\ft}\cong i\ubarRR\oplus L$, the operators on
the fibers defined by these families are then given by
\begin{equation}
  \label{eq:LambdaWSplitUnivBundleOps}
  \begin{aligned}
\delta_{A_L} \equiv
d^*_{A_L}+d^+_{A_L}
&:
 C^\8(\Lambda^1\otimes L)
\to
 C^\8((\Lambda^0\oplus \Lambda^+)\otimes L),
\\
D''_{B\otimes A_L}
&:
 C^\8(W^+\otimes L)
\to
 C^\8(W^-\otimes L),
\end{aligned}
\end{equation}
where $D''_{B\otimes A_L}$ is the Dirac operator.  We use the preceding
families to rewrite $\ind \bsD^n$ in terms of index bundles whose Chern
characters are more readily computable:

\begin{lem}
\label{lem:IdentifyNormalWithDirac}
Continue the above notation.  Then, as elements of $K(M_{\fs})$,
$$
\ind \bsD^n = \ind \bdelta + \ind \bD''.
$$
\end{lem}

\begin{proof}
As an element of $K(M_{\fs})$, the index bundle $\ind\bsD^n$
depends only on the homotopy class of the leading symbol \cite[Theorem
III.8.6]{LM}.  Thus, the index bundle of the family $(B,\Psi)\mapsto
\sD^n_{B,\Psi}$ is equivalent to that of the family
$(B,\Psi)\mapsto \sD^n_{B,0}$, where $\sD^n_{B,0} =
\delta_{A_L} + D''_{B\otimes A_L}$.
\end{proof}

We now identify the family of operators $\bdelta$ in
\eqref{eq:DeltaAndDiracL2Families} with a family of Dirac operators.
(Identifications of this type were used in \cite{AHS} to compute the index
of the elliptic deformation complex for the anti-self-dual equation.)  The
vector bundle $\End(W)\cong W\otimes W^*$ is a Clifford module by Clifford
multiplication on the factor $W$ in the tensor product and so we obtain a
Dirac operator,
\begin{equation}
  \label{eq:End(W)DiracOperator}
D: C^\8(W^\pm\otimes W^*)\to  C^\8(W^\mp\otimes W^*),
\end{equation}
on $\End(W)$ \cite[p. 122-123]{LM}.  Under the identification
$\Lambda^\bullet\otimes_\RR\CC$ with $\End(W)$ given by Clifford
multiplication, the operator $d^*+d$ on $\Lambda^\bullet$ is identified
with the Dirac operator \eqref{eq:End(W)DiracOperator} (see \cite[Theorem
II.5.12]{LM}).  Restricting the domain of Clifford multiplication and
tensoring with the line bundle $L$ gives isomorphisms
\begin{equation}
\label{eq:CliffMultIsom}
\begin{aligned}
\Lambda^1\otimes_\RR L
&\cong \Hom(W^+,W^-)\otimes_\CC L,
\\
(\Lambda^0\oplus\Lambda^+)\otimes_\RR L
&\cong \End(W^+)\otimes_\CC L,
\end{aligned}
\end{equation}
and therefore isomorphisms
$$
\bE' \cong \bE_-'
\quad\text{and}\quad
\bF' \cong \bE_+',
$$
where the bundles $\bE_{\pm}'$ are defined by
\begin{equation}
  \label{eq:bEpmLambdaBundleDefns}
\begin{aligned}
\bE_-'
&=
\tM_{\fs}\times_{\sG_{\fs}}
(\Hom(W^+,W^-)\otimes L),
\\
\bE_+'
&=
\tM_{\fs}\times_{\sG_{\fs}}
(\End(W^+)\otimes L).
\end{aligned}
\end{equation}
The operator $d^*_{A_L}+d^+_{A_L}$ is the restriction of
$d^*_{A_L}+d_{A_L}$ to $\Lambda^1\otimes_\RR L$ composed with the
projection from $(\Lambda^0\oplus\Lambda^2)\otimes_\RR L$ to
$(\Lambda^0\oplus\Lambda^+)\otimes_\RR L$.  Thus, the isomorphisms
\eqref{eq:CliffMultIsom} identify the family of operators $\bdelta$ in
\eqref{eq:DeltaAndDiracL2Families} with the family of Dirac operators:
%
%
\begin{equation}
\label{eq:DiracEndWL2L1Family}
\begin{CD}
 C^\8(\bE_-') @>{\bD'}>>
 C^\8(\bE_+')    \\
@VVV                         @VVV         \\
M_{\fs}\times X   @>{\id}>> M_{\fs}\times X
\end{CD}
\end{equation}
Explicitly, if $(B,\Psi)$ is a pair representing a point $[B,\Psi]$
in $M_{\fs}$, the Dirac operator on the fiber given by this family,
\begin{equation}
\label{eq:DefnDeltaDirac}
D'_{B\otimes A_L}:
 C^\8(\Hom(W^+,W^-)\otimes L)
\to
 C^\8(\End(W^+)\otimes L),
\end{equation}
is defined by the \spinc connection on $\End(W)$ induced by the \spinc
connection $B$ on $W$ and the unitary connection $A_L=A_\Lambda\otimes
(B^{\det})^*$ on $L$.  The preceding identification yields:

\begin{lem}
\label{lem:IdentifyDeltaWithDirac}
Continue the above notation. Then, as elements of $K(M_{\fs})$,
$$
\ind \bdelta
=
\ind \bD'.
$$
\end{lem}

We now begin the second step, which is to compute the Chern characters of
the index bundles of the families of Dirac operators in
\eqref{eq:DeltaAndDiracL2Families} and \eqref{eq:DiracEndWL2L1Family}.
The following technical lemma helps to identify
the universal bundles $\bE_\pm'$, $\bE''$, and $\bF''$.

\begin{lem}
\label{lem:DiagonalQuotient}
If $Q_i\rightarrow M$, $i=1,2$, are $S^1$ bundles over a manifold $M$ and
$L_i=Q_i\times_{S^1}\CC$, $i=1,2$ are the associated complex line bundles,
and $k\in\ZZ$, then the following hold:
\begin{enumerate}
\item If $V\to M$ is a complex vector bundle, and
  $e^{i\theta}\in S^1$ acts on the fiber product $Q_1\times_M V$ by
  $e^{i\theta}\cdot(q_1,v)=(e^{i\theta}q_1,e^{ik\theta}v)$, then
  $(Q_1\times_M V)/S^1\cong L_1^{-k}\otimes V$.
\item If $e^{i\theta}\in S^1$ acts on the fiber product $Q_1\times_MQ_2$ by
  $e^{i\theta}\cdot(q_1,q_2)=(e^{i\theta}q_1,e^{ik\theta}q_2)$, then the
  first Chern class of the $S^1$-bundle $(Q_1\times_M Q_2)/S^1\rightarrow M$
  is $c_1(Q_2)-kc_1(Q_1)$, where the action of $S^1$ on $(Q_1\times_M
  Q_2)/S^1$ is induced by the $S^1$ action on $Q_2$ of weight one.
\end{enumerate}
\end{lem}

\begin{proof}
  The associated line bundles $L_i=Q_i\times_{S^1}\CC$ are given by the
  quotients of $Q_i\times\CC$ by the relation $(q_i,z)\sim
  (e^{-i\theta}q_i,e^{i\theta}z)$, for $(q_i,z)\in Q_i\times\CC$ and
  $e^{i\theta}\in S^1$; under the same relation, $Q_1\times_{S^1}S^1= Q_1$.
  A tensor product of a complex line bundle $L$ with a
  complex vector bundle $V$ is given by the quotient of the fiber product
  $L\times_MV$ by the relation $(z,v)\sim (w^{-1}z,wv)$, for $z\in L,v\in
  V$, and $w\in\CC^*$. Hence,
\begin{align*}
L_1^{-k}\otimes V
&=
\{([q_1,z],v)\in L_1\times_M V: ([q_1,z],v)\sim ([q_1,wz],w^kv), w\in\CC^*\}
\\
&=
\{([q_1,e^{i\mu}],v)\in Q_1\times_M V: ([q_1,e^{i\mu}],v)\sim 
([q_1,e^{i\theta}e^{i\mu}],e^{ik\theta}v), e^{i\theta}\in S^1\}
\\
&=
\{(p_1,v)\in Q_1\times_M V: (p_1,v)\sim (e^{i\theta}p_1,e^{ik\theta}v),
e^{i\theta}\in S^1\}
\quad\text{(where $p_1=[q_1,e^{i\mu}]$)}
\\
&= (Q_1\times_M V)/S^1,
\end{align*}
where in the second line above we can assume without loss that $z\neq 0$.
This proves Assertion (1) and Assertion (2) follows trivially from this.
\end{proof}

Since it will not cause confusion, we will write $\LL_{\fs}$ for the
restriction of the universal bundle of equation \eqref{eq:SWUnivLineBundle}
to $M_{\fs}\times X$,
$$
\LL_{\fs}= \tM_{\fs}\times_{\sG_{\fs}}\underline{\CC},
$$
and define
\begin{equation}
\label{eq:SlantSWUnivLineBundle}
\LL_{\fs,x}=\LL_{\fs}|_{M_{\fs}\times\{x\}}.
\end{equation}
By construction, $c_1(\LL_{\fs,x})=c_1(\LL_\fs)/x=\mu_{\fs}(x)$ where
$x\in H_0(X;\ZZ)$ is a generator.  Let $\pi_M:M_{\fs}\times X\to
M_{\fs}$ and $\pi_X:M_{\fs}\times X\to X$ be the
projections.

\begin{lem}
\label{lem:IndexBundleIsomorphism}
Assume $M_{\fs}$ contains no zero-section pairs.  If $\br^*\bDelta\to
M_{\fs}\times X$ is the restriction of the pullback by the retraction $\br$
of Lemma \ref{lem:Retraction} of the line bundle \eqref{eq:JacobianBundle}
then, using $\Hom(W^+,W^-)\cong W^-\otimes W^{+,*}$ and $\End(W^+)\cong
W^+\otimes W^{+,*}$, there are isomorphisms,
\begin{align*}
\begin{array}{ll}
\bE'_-
\cong \sE'\otimes \pi_X^*W^-,
&\quad
\bE'_+
\cong \sE'\otimes \pi_X^*W^+,
\\
\bE''
\cong \sE''\otimes \pi_X^*W^+,
&\quad
\bF''
\cong \sE''\otimes \pi_X^*W^-,
\end{array}
\end{align*}
where we have defined
\begin{equation}
\label{eq:DefnsE'sE''}
\begin{aligned}
\sE'
&=
(\pi_M^*\LL_{\fs,x}\otimes\br^*\bDelta)^{\otimes 2}
\otimes
\pi_X^*(W^{+,*}\otimes L),\\
\sE''
&=
\pi_M^*\LL_{\fs,x}\otimes\br^*\bDelta\otimes \pi_X^*L.
\end{aligned}
\end{equation}
\end{lem}

\begin{proof}
There are isomorphisms,
\begin{equation}
\label{eq:IndexIsomorphism1}
\bE'_{\pm}\cong
(\tM_{\fs}\times_{\sG_{\fs}} L)
\otimes
\Hom(W^+,W^\pm),
\end{equation}
where, from equation \eqref{eq:su(E)SplitGaugeGroupAction}, the action of
$s\in\sG_{\fs}$ on $\tM_{\fs}\times L$ is given by $(B,\Psi,z)\mapsto
(s(B,\Psi),s^{-2}z)$.  Let $\tM_{\fs}\times_{\sG_{\fs}}(X\times S^1)$ be
the unit sphere bundle of $\LL_{\fs}$.  If $e^{i\mu}\in S^1$ acts on
$$
\left(\tM_{\fs}\times_{\sG_{\fs}}(X\times S^1)\right)
\times_{M_{\fs}\times X}\pi_X^*L
$$
by $[(B,\Psi),x,e^{i\theta},z]\mapsto
[(B,\Psi),x,e^{i\mu}e^{i\theta},e^{-i2\mu}z]$---where
$(B,\Psi)\in\tM_{\fs}$, $x\in X$, $e^{i\theta}\in S^1$ and $z\in
\pi_X^*L$---then we obtain an isomorphism of complex line bundles
\begin{equation}
\label{eq:IndexIsomorphism2}
\begin{aligned}
\left(\left(\tM_{\fs}\times_{\sG_{\fs}}(X\times S^1)\right)
\times_{M_{\fs}\times X}\pi_X^*L\right)/S^1
& \to
\tM_{\fs}\times_{\sG_{\fs}} L,
\\
[(B,\Psi),x,e^{i\theta},z]
&\mapsto
[(B,\Psi),e^{2i\theta}z].
\end{aligned}
\end{equation}
Lemma \ref{lem:DiagonalQuotient} then implies that
\begin{equation}
\label{eq:IndexIsomorphism3}
\tM_{\fs}\times_{\sG_{\fs}} L
\cong
\LL_{\fs}^{\otimes 2}\otimes \pi_X^*L
\end{equation}
and equation \eqref{eq:IndexIsomorphism3} gives
\begin{equation}
\bE'_{\pm}
\cong
\LL_{\fs}^{\otimes 2}\otimes \pi_X^*(L \otimes W^{+,*}\otimes W^\pm).
\end{equation}
The desired expression for $\bE'_\pm$ then follows from the isomorphism
$\LL_{\fs}\cong \LL_{\fs,x}\otimes\br^*\bDelta$ given
by Lemma \ref{lem:MuMapDescription}.

In the bundles
$$
\bE''
\cong
\tM_{\fs}\times_{\sG_{\fs}} (W^+\otimes L)
\quad\text{and}\quad
\bF''\cong
\tM_{\fs}\times_{\sG_{\fs}} (W^-\otimes L),
$$
an element $s\in\sG_{\fs}$ acts on $\tM_{\fs}\times W^\pm\otimes L$ by
$((B,\Psi),\Psi')\mapsto (s(B,\Psi),s^{-1}\Psi')$, where $\Psi'\in
 C^\8(W^\pm\otimes L)$, as noted in the remark following equation
\eqref{eq:su(E)SplitGaugeGroupAction}. Lemma \ref{lem:DiagonalQuotient} and
the argument yielding equation \eqref{eq:IndexIsomorphism3} then imply that
there are isomorphisms of complex vector bundles
\begin{equation}
\label{eq:IndexIsomorphism4}
\bE''
\cong
\LL_{\fs}\otimes \pi_X^*(W^+\otimes L)
\quad\text{and}\quad
\bF''
\cong
\LL_{\fs}\otimes \pi_X^*(W^-\otimes L).
\end{equation}
The isomorphisms in the conclusion of the lemma now follow from equation
\eqref{eq:IndexIsomorphism4} and the isomorphism $\LL_{\fs}\cong
\pi_M^*\LL_{\fs,x}\otimes \br^*\bDelta$
implied by equation \eqref{eq:ChernClassOfUniversalSW}.
\end{proof}

Given Lemma \ref{lem:IndexBundleIsomorphism}, the index bundles of the
families of Dirac operators in \eqref{eq:DeltaAndDiracL2Families} and
\eqref{eq:DiracEndWL2L1Family} now take the shape:
%
%
$$
\begin{CD}
 C^\8(\sE'\otimes W^-)@>{\bD'}>>
 C^\8(\sE'\otimes W^+) \\
@VVV                         @VVV         \\
M_{\fs}\times X   @>{\id}>> M_{\fs}\times X
\end{CD}
\quad\text{and}\quad
\begin{CD}
 C^\8(\sE''\otimes W^+) @>{\bD''}>>  C^\8(\sE''\otimes W^-)    \\
@VVV                         @VVV         \\
M_{\fs}\times X   @>{\id}>> M_{\fs}\times X
\end{CD}
$$
We abuse notation slightly by continuing to denote the families of Dirac
operators on these isomorphic bundles by $\bD'$ and $\bD''$.

The final step in the computation of the Chern character of
$N_{\ft}(\Xi,\fs)$ is to compute the Chern characters of these families.

Given the decomposition \eqref{eq:NormalReducDefOperator} of
$\sD_{\iota(B,\Psi)}^n$, it will be convenient to define
\begin{equation}
\label{eq:SplittingNormalIndex}
\ind_\CC\sD_{\iota(B,\Psi)}^n = n_s = n_s' + n_s'',
\quad\text{where}\quad
\begin{aligned}
n_s'
&=
\ind_\CC (d^*_{A_L}+d^+_{A_L}),
\\
n_s''
&=
\ind_\CC D''_{B\otimes A_L}.
\end{aligned}
\end{equation}
Viewing $E=\underline{\CC}\oplus L$ and $\su(E)=i\underline{\RR}\oplus L$,
we can compute the complex index of $d^*_{A_L}+d^+_{A_L}$ from the real
index of $d^*_{\hat A}+d^+_{\hat A}$ on $C^\8(\Lambda^1\otimes\su(E))$ (for
example, see \cite[Equation (4.2.22)]{DK}) and the fact (Lemma
\ref{lem:ReducibleSpinu}) that $d^*_{\hat A}+d^+_{\hat A}$ decomposes as
the direct sum of $d^*+d^+$ on $C^\8(i\Lambda^1)$ and $d^*_{A_L}+d^+_{A_L}$
on $C^\8(\Lambda^1\otimes L)$; the complex index of the \spinc Dirac
operator is given in \cite[p.  47]{MorganSWNotes}. Thus,
\begin{equation}
\label{eq:DefOfNormalIndices}
\begin{aligned}
n_s'(\ft,\fs)
&=-\left(c_1(\fs)-c_1(\ft)\right)^2 - \textstyle{\frac{1}{2}}(\chi+\sigma),
\\
n_s''(\ft,\fs)
&=
\textstyle{\frac{1}{8}}(\left(2c_1(\ft)-c_1(\fs)\right)^2-\sigma).
\end{aligned}
\end{equation}
We can now state and prove one of the main results of our article:

\begin{thm}
\label{thm:ChernCharacterOfNormal}
Let $\ft$ be a \spinu structure over a closed, oriented, Riemannian, smooth
four-manifold $X$.  Assume $\ft$ admits a splitting $\ft=\fs\oplus
\fs\otimes L$, where $L$ is a complex line bundle.  Suppose that there are
no zero-section pairs in $M_{\fs}$.  Let $\mu_{\fs}$ be the Seiberg-Witten
$\mu$-map, as in definition \eqref{eq:SWMuMap}.  Let $x\in H_0(X;\ZZ)$ be
the positive generator and for $\{\gamma_i\}$ a basis for
$H_1(X;\ZZ)/\Tor$, let $\{\gamma_i^*\}$ be the dual basis for $H^1(X;\ZZ)$
introduced in Definition \ref{defn:RelatedBasis}. Let $r_\Xi =
\rank_\CC\Xi$ denote the rank of the stabilizing bundle.  Then the Chern
character of the normal bundle $N_{\ft}(\Xi,\fs)$ of the stratum
$M_{\fs}\embed \sM_{\ft}(\Xi,\fs)$ is
\begin{equation}
\label{eq:ChernCharacterOfNormal}
\begin{aligned}
\ch( N_{\ft}(\Xi,\fs))
&=
r_\Xi + n_s'' e^{\mu_{\fs}(x)} + n_s' e^{2\mu_{\fs}(x)} 
\\
&\quad -
8\sum_{i<j}
\langle \gamma_i^*\gamma^*_j(c_1(\ft)-c_1(\fs)),[X]\rangle
e^{2\mu_{\fs}(x)}\mu_{\fs}(\gamma_i\gamma_j)
\\
&\quad +
\textstyle{\frac{1}{2}}\sum_{i<j}
\langle \gamma_i^*\gamma_j^*(2c_1(\ft)-c_1(\fs)),[X]\rangle
e^{\mu_{\fs}(x)}\mu_{\fs}(\gamma_i\gamma_j)
\\
&\quad +
(e^{\mu_{\fs}(x)}-32e^{2\mu_{\fs}(x)})\sum_{i<j<k<\ell}
\langle \gamma_i^*\gamma_j^*\gamma_k^*\gamma_\ell^*,[X]\rangle
\mu_{\fs}(\gamma_i\gamma_j\gamma_k\gamma_\ell).
\end{aligned}
\end{equation}
where $\mu_{\fs}(\gamma_i\gamma_j)
=\mu_{\fs}(\gamma_i)\mu_{\fs}(\gamma_j)$ and similarly for
$\mu_{\fs}(\gamma_i\gamma_j\gamma_k\gamma_\ell)$.
\end{thm}

\begin{proof}
For convenience in the proof, we write $\mu=\mu_{\fs}(x)$ and
$\gamma_i^{J,*}=\mu_{\fs}(\gamma_i)$, as in Lemma \ref{lem:MuMapDescription}.
The K-theory identification \eqref{eq:NormalBundleIsNormalStab} of
$N_{\ft}(\Xi,\fs)$, Lemmas \ref{lem:IdentifyNormalWithDirac} and
\ref{lem:IdentifyDeltaWithDirac}, and the homomorphism property
of the Chern character (see, for example, \cite[Proposition III.11.16]{LM})
imply that
\begin{equation}
\label{eq:ChernCharNormalBundle}
\begin{aligned}
\ch(N_{\ft}(\Xi,\fs))
&= \ch(\ind \bsD^n)+\ch(\Xi)
\\
&=
\ch\bigl(\ind \bD'\bigr)+\ch(\ind \bD'')+r_\Xi.
\end{aligned}
\end{equation}
Recall that $\ind(D^*,\sE,W)$ denotes the index bundle of the family of
operators obtained, as in Proposition \ref{prop:AtiyahSinger}, by twisting
the Dirac operator $D^*: C^\8(W^-)\to  C^\8(W^+)$ by a family of
connections on the bundle $\sE$. Note that $\ch(\ind(D^*,W)) =
-\ch(\ind(D,W))$, where $D^*: C^\8(W^-)\to C^\8(W^+)$.  We now use
Lemma \ref{lem:IndexBundleIsomorphism} and Proposition
\ref{prop:AtiyahSinger} to partly compute the Chern characters of the index
bundle $\ind \bD'$.  Bundles over $X$
will be considered to be bundles over $M_{\fs}\times X$; the
pullback $\pi_X^*$ will be omitted.
\begin{equation}
\label{eq:ChernCharacterDelta}
\begin{aligned}
{}&\ch\bigr(\ind \bD'\bigl)
= \ch(\ind(D^*,\sE',W))
= \ch(\sE')\ch(\ind(D^*,W))
\\
&\quad = -\ch\left((\pi_M^*\LL_{\fs,x}\otimes\br^*\bDelta)^{\otimes 2}
\otimes W^{+,*}\otimes L\right)\ch(\ind(D,W))
\\
&\quad =
-\pi_M^*\ch(\LL_{\fs,x})^2\br^*\ch(\bDelta)^2
\ch(W^{+,*})e^{c_1(L)} e^{\half c_1(W^+)}\hat\sA(X)/[X].
\end{aligned}
\end{equation}
Similarly, we use Lemma \ref{lem:IndexBundleIsomorphism} and
Proposition \ref{prop:AtiyahSinger} to partly compute the Chern character
of the index bundle $\ind \bD''$:
\begin{equation}
\label{eq:ChernCharacterDirac}
\begin{aligned}
\ch(\ind \bD'')
&= \ch(\ind(D,\sE'',W))
= \ch(\sE'')\ch(\ind(D,W))
\\
&= \ch(\pi_M^*\LL_{\fs,x}\otimes\br^*\bDelta\otimes L)
e^{\half c_1(W^+)}\hat\sA(X)/[X]
\\
&= \pi_M^*\ch(\LL_{\fs,x})\br^*\ch(\bDelta)e^{c_1(L)}
e^{\half
c_1(W^+)}\hat\sA(X)/[X].
\end{aligned}
\end{equation}
In the first lines of \eqref{eq:ChernCharacterDelta} and
\eqref{eq:ChernCharacterDirac} we simply rewrite the
index bundles using the notation of Proposition \ref{prop:AtiyahSinger}.

For the calculation of $\ch(\ind \bD'')$ in
\eqref{eq:ChernCharacterDirac}, we compute $\ch(\bDelta)$ using the
expression for $\br^*c_1(\bDelta)$ in Lemmas \ref{lem:CanonicalLineBundle}
and \ref{lem:MuMapDescription}:
$$
\br^*\ch(\bDelta)
=\prod_{i=1}^{b_1(X)}e^{\gamma_i^{J,*}\times\gamma_i^*}
=\prod_{i=1}^{b_1(X)}\left(1+\gamma_i^{J,*}\times\gamma_i^*\right).
$$
Because $c_1(\LL_{\fs,x})=\mu$ as noted
before Lemma \ref{lem:IndexBundleIsomorphism}, we see
$\pi_M^*\ch(\LL^{*}_{\fs,x})=\pi_M^*e^{\mu}$. We shall
write $\mu$ for $\pi_M^*\mu$.
Applying this to \eqref{eq:ChernCharacterDirac} and
noting that
$c_1(L)+\half c_1(W^+)=c_1(\ft)-\half c_1(\fs)$, we obtain
\begin{align*}
&\ch(\ind \bD'') \\
&\quad =
e^{\mu}
\left(\prod_{i=1}^{b_1(X)}e^{\gamma_i^{J,*}\times\gamma_i^*}\right)
e^{c_1(\ft)-\textstyle{\frac{1}{2}} c_1(\fs)}
\left(1-\textstyle{\frac{1}{24}}p_1(X)\right)/[X] \\
&\quad =
e^{\mu}
\left(1+\sum_{i<j}(\gamma_i^{J,*}\gamma_j^{J,*})
\times(\gamma_i^*\gamma_j^*)
+\sum_{i<j<k<\ell}
(\gamma_i^{J,*}\gamma_j^{J,*}\gamma_k^{J,*}\gamma_\ell^{J,*})\times
(\gamma_i^*\gamma_j^*\gamma_k^*\gamma_\ell^*)\right) \\
&\qquad
\times\left( 1+\textstyle{\frac{1}{2}}(2c_1(\ft)-c_1(\fs))
+\textstyle{\frac{1}{8}}(2c_1(\ft)-c_1(\fs))^2\right)
\left(1-\textstyle{\frac{1}{24}}p_1(X)\right)/[X],
\end{align*}
and therefore, using $\si=\frac{1}{3}\langle p_1(X),[X]\rangle$
(see \cite[Theorem I.3.1]{BPV})
\begin{equation}
\label{eq:ChernCharIndDiracL2}
\begin{aligned}
\ch(\ind \bD'')
&=
e^{\mu}\Bigl( \textstyle{\frac{1}{8}} (2c_1(\ft)-c_1(\fs))^2
               -\textstyle{\frac{1}{8}} \sigma
\\
&\quad+\textstyle{\frac{1}{2}}\sum_{i<j}(\gamma_i^{J,*}\gamma_j^{J,*})
\langle\gamma_i^*\gamma_j^*(2c_1(\ft)-c_1(\fs)),[X]\rangle\\
&\quad + \sum_{i<j<k<\ell}
(\gamma_i^{J,*}\gamma_j^{J,*}\gamma_k^{J,*}\gamma_\ell^{J,*})
\langle(\gamma_i^*\gamma_j^*\gamma_k^*\gamma_\ell^*),[X]\rangle\Bigr),
\end{aligned}
\end{equation}
completing the calculation of $\ch(\ind \bD'')$.

We now complete the calculation of $\ch(\ind \bdelta)$
in \eqref{eq:ChernCharacterDelta}.  To compute $\ch(W^{+,*})$ we
observe that $c_1(W^{+,*})=-c_1(W^+)$ while $c_2(W^{+,*}) = c_2(W^+)$ and
so, using this and the isomorphism $\su(W^+)\cong \Lambda^+$, we have
\begin{align*}
\ch(W^{+,*})
&=
2 + c_1(W^{+,*}) 
+ \textstyle{\frac{1}{2}}\left(c_1(W^{+,*})^2 - 2c_2(W^{+,*}\right)
\\
&=
2 - c_1(W^+) + \textstyle{\frac{1}{2}}c_1(W^+)^2
+ \textstyle{\frac{1}{4}}\left(p_1(\su(W^+))-c_1(W^+)^2\right)
\\
&=
2 - c_1(W^+) + \textstyle{\frac{1}{4}}\left(p_1(\Lambda^+)+c_1(W^+)^2\right).
\end{align*}
Applying this to the terms involving $W^+$ in
\eqref{eq:ChernCharacterDelta} and writing $c_1(W^+) = c_1(\fs)$ yields
\begin{align*}
&-\ch(W^{+,*}) e^{\half c_1(W^+)}\hat\sA(X)/[X]
\\
&\quad =
-\left( 2-c_1(\fs)+\textstyle{\frac{1}{4}}(c_1(L)^2+p_1(\Lambda^+))\right)
\left(1+\textstyle{\frac{1}{2}} c_1(\fs) +\textstyle{\frac{1}{8}}
  c_1(L)^2\right)
\left(1-\textstyle{\frac{1}{24}}p_1(X)\right)/[X]
\\
&\quad =
-(2+\textstyle{\frac{1}{4}} p_1(\Lambda^+))
\left(1-\textstyle{\frac{1}{24}}p_1(X)\right)/[X].
\end{align*}
Hence, using the preceding expression in \eqref{eq:ChernCharacterDelta}
and $c_1(L)=c_1(\ft)-c_1(\fs)$, we see that
\begin{align*}
&\ch(\ind \bdelta)
\\
&=
-\pi_M^*\ch(\LL_{\fs,x})^2\br^*\ch(\bDelta)^2\ch(L)
(2+\textstyle{\frac{1}{4}} p_1(\Lambda^+))
\left(1-\textstyle{\frac{1}{24}}p_1(X)\right)/[X]
\\
&=
-e^{2\mu}
\left(\prod_{i=1}^{b_1(X)}\left(1+2\gamma_i^{J,*}\times\gamma_i^*\right)\right)
\left(1+(c_1(\ft)-c_1(\fs))+\textstyle{\frac{1}{2}}(c_1(\ft)-c_1(\fs))^2\right)\\
&\qquad\qquad
\times\left(2+\textstyle{\frac{1}{4}} p_1(\Lambda^+)\right)
\left(1-\textstyle{\frac{1}{24}}p_1(X)\right)/[X].
\end{align*}
Simplifying the preceding expression for
$\ch(\ind \bdelta)$ and recalling that
$p_1(\Lambda^+)=2\chi+3\sigma$ (from \cite[Satz 1.5]{HirzebruchHopf}), yields
\begin{align*}
&\ch(\ind \bdelta) \\
&=
-e^{2\mu}
\Bigl(1+4\sum_{i<j}(\gamma_i^{J,*}\gamma_j^{J,*})
\times(\gamma_i^*\gamma_j^*)
+16\sum_{i<j<k<\ell}
(\gamma_i^{J,*}\gamma_j^{J,*}\gamma_k^{J,*}\gamma_\ell^{J,*})
\times (\gamma_i^*\gamma_j^*\gamma_k^*\gamma_\ell^*)\Bigr) \\
&\qquad\qquad
\times\left(1+(c_1(\ft)-c_1(\fs))
+\textstyle{\frac{1}{2}}(c_1(\ft)-c_1(\fs))^2\right)
        \left(2+\textstyle{\frac{1}{4}} p_1(\Lambda^+)\right)
        \left(1-\textstyle{\frac{1}{24}}p_1(X)\right)/[X],
\end{align*}
and thus,
\begin{equation}
\label{eq:ChernCharIndDelta}
\begin{aligned}
\ch(\ind \bdelta)
&=
e^{2\mu}\Bigl(-(c_1(\ft)-c_1(\fs))^2-\textstyle{\frac{1}{2}}\sigma
-\textstyle{\frac{1}{2}}\chi
\\
&\quad
- 8\sum_{i<j}(\gamma_i^{J,*}\gamma_j^{J,*})
\langle(\gamma_i^*\gamma_j^*)(c_1(\ft)-c_1(\fs)),[X]\rangle
\\
&\quad
-32\sum_{i<j<k<\ell}
(\gamma_i^{J,*}\gamma_j^{J,*}\gamma_k^{J,*}\gamma_\ell^{J,*})
\langle \gamma_i^*\gamma_j^*\gamma_k^*\gamma_\ell^*,[X]\rangle\Bigr).
\end{aligned}
\end{equation}
The desired expression for $\ch(N_{\ft}(\Xi,\fs))$ follows from
\eqref{eq:DefOfNormalIndices}, \eqref{eq:ChernCharNormalBundle},
\eqref{eq:ChernCharIndDiracL2}, and \eqref{eq:ChernCharIndDelta}.
\end{proof}

Finally, we calculate the total Chern class $c(N_{\ft}(\Xi,\fs))$, as an
element of rational cohomology, under a simplifying assumption.

\begin{cor}
\label{cor:SimpleNormalChernClass}
Continue the hypotheses of Theorem \ref{thm:ChernCharacterOfNormal} and
assume that $\alpha\smile\alpha'=0$, for every $\alpha,\alpha'\in
H^1(X;\ZZ)$.  Then, as elements of $H^\bullet(M_{\fs};\RR)$,
$$
c(N_{\ft}(\Xi,\fs))
=
\left(1+2\mu_{\fs}(x)\right)^{n_s'}\left(1+\mu_{\fs}(x)\right)^{n_s''}.
$$
\end{cor}

\begin{proof}
For convenience, we write $\mu=\mu_{\fs}(x)$. The Chern character
determines the Chern polynomial as an element of
rational cohomology (see the formula in
\cite[Problem 16-A]{MilnorStasheff} or \cite[pp. 156--157]{ACGH}), so
$\ch(N_{\ft}(\Xi,\fs))$ determines
$$
c_t(N_{\ft}(\Xi,\fs)) = \sum_{i=0}^{r_\Xi} c_i(N_{\ft}(\Xi,\fs))t^i
=
\prod_{j=1}^{r_\Xi}(1+\alpha_i t),
$$
where the $\alpha_i$ are the Chern roots of $N_{\ft}(\Xi,\fs)$.
Theorem \ref{thm:ChernCharacterOfNormal} implies that, with the
constraint on $H^1(X;\ZZ)$,
$$
\ch( N_{\ft}(\Xi,\fs))
=
\sum_{j=1}^{r_\Xi} \exp(\alpha_j)
=
r_\Xi+n_s'e^{2\mu}
+
n_s''e^{\mu}.
$$
Suppose $n_s'\geq 0$ and $n_s''\geq 0$. If
$\prod_{j=1}^{r_\Xi}(1+\alpha_i t)
=
(1+2\mu t)^{n_s'}(1+\mu t)^{n_s''}$,
then $\alpha_j = +2\mu$ for $1\le j\le n_s'$ and $\alpha_j = \mu$
for $n_s'<j\le n_s'+n_s''$, while $\alpha_j = 0$ for
$n_s'+n_s'' < j\le r_\Xi$. Then one can easily see that
$\ch(N_{\ft}(\Xi,\fs))$ is equal to the Chern character associated to the
Chern polynomial $(1+2\mu t)^{n_s'}(1+\mu t)^{n_s''}$.
The case where either $n_s'$ or $n_s''$ is
negative follows from the observation that if $G\in K(M_{\fs})$,
then $\ch(-G)=\ch(G)$ while $c(-G)=c(G)^{-1}$.
\end{proof}


\appendix
\section{Abundant four-manifolds}
\label{sec:Abundance}
Our goal in this section is to prove the 

\begin{thm}
\label{thm:Abundance}
Every compact, complex algebraic, simply connected surface with $b_2^+ \geq 3$
is abundant.
\end{thm}

In \cite[p. 175]{FKLM} we asserted without further explanation that simply
connected, minimal, complex algebraic surfaces of general type were
abundant. On the other hand, some of the fake K3-surfaces of
\cite{GompfMrowka} fail to be abundant.  If log transforms are performed on
tori in three distinct nuclei then the intersection form on $B^{\perp}$ is
a degenerate form with three-dimensional radical and having an
$-E_{8}\oplus -E_{8}$ summand \cite[p. 175]{FKLM}. We apply Theorem
\ref{thm:Abundance} to produce classes $\Lambda\in H^2(X;\ZZ)$ which are
orthogonal to the SW-basic classes, with square equal to prescribed even
integers. As far as we can tell, one can always find such classes $\Lambda$
for compact four-manifolds with $b_2^+\geq 3$, even if non-abundant: it is
interesting problem to determine if indeed this is true.

We are extremely grateful to Andr\'as Stipsicz for describing a proof of
Theorem \ref{thm:Abundance} in the case where $X$ is a minimal surface of
general type with odd intersection form: see Lemmas \ref{lem:CrudeKexists}
and \ref{lem:UnstableRangeKexists}, as well as the ideas for the case of
odd minimal surfaces of general type in \S
\ref{subsec:CptCplxSimpConnSurface}---these are all due to him. The
argument for Lemma \ref{lem:CrudeKexists} relies on the ``odd four-square
theorem'' (Lemma \ref{lem:FourOddSquare}), for whose proof we are indebted
to A. Agboola \cite{AgboolaJune2000}.

For the remainder of this section, unless further restrictions are
mentioned, we suppose $X$ is a compact, connected, smooth four-manifold
with an orientation for which $b_2^+(X)>0$. 

\subsection{Four-manifolds whose basic classes are multiples of a fixed class}
We consider the cases of even and odd intersection forms separately,
beginning with the even case:

\begin{lem}
\label{lem:OneBasicClass}
Suppose $(L,Q)$ is an indefinite, integral, unimodular lattice and that
$\kappa\in L$. If $Q$ is even and $|\sigma(Q)|\leq \rank(Q)-4$ then 
$\kappa^\perp$, the $Q$-orthogonal complement of $\kappa$ in $L$, contains
a hyperbolic sublattice. 
\end{lem}

\begin{proof}
  We may assume without loss that $\kappa$ is primitive because a
  sublattice $H$ is orthogonal to $\kappa=d\kappa'$ if and only if it is
  orthogonal to $\kappa'$, where $d\in\ZZ$ and $\kappa'$ is primitive.
  
  {}From the classification of indefinite, integral, unimodular forms (for
  example, see Theorem 1.2.21 in \cite{GompfStipsicz}) we have
  \begin{equation}
    \label{eq:classification}
(L,Q)
\cong 
\textstyle{\frac{1}{8}}\sigma(Q)E_8 
\oplus 
\textstyle{\frac{1}{2}}(\rank(Q)-|\sigma(Q)|)H.
  \end{equation}
  By hypothesis, $\rank(Q)-|\sigma(Q)| \geq 4$, so if $R=H\oplus H$ then
  $(R,Q|_R)$ is a sublattice of $(L,Q)$. Let $\ka_R$ denote the component
  of $\ka$ in $R$. Because $R$ is even, we have $Q(\ka_R,\ka_R)=2h$ for some
  $h\in\ZZ$.  Let $e_1,e_2$ and $f_1,f_2$ be bases for the two hyperbolic
  sublattices of $R$, so $Q(e_1,e_2)=Q(f_1,f_2)=1$ and all other pairings
  of these four vectors vanish. Define $v=e_1+ h e_2$, so $v^2=2h$ and $v$
  is primitive.  The hyperbolic sublattice $F=\ZZ f_1+\ZZ f_2 \subset R$ is
  orthogonal to $v$.  According to Theorem 1 in \cite{WallUnimodQuadForms},
  the orthogonal group of $(R,Q|_R)$ acts transitively on primitive vectors
  of a given square.  Hence, there is an automorphism $A$ of $R$ such that
  $Q(Ax,Ay) = Q(x,y)$ for all $x,y\in R$ and $Av=\ka_R$.  Because $F$ is
  contained in $v^\perp$, then $A(F)\subset \ka_R^\perp$.  Hence, the
  hyperbolic sublattice $A(F)$ is contained in $\ka^\perp$.
\end{proof}

\begin{cor}
\label{cor:OneBasicClassEven}
Let $X$ be a simply-connected spin four-manifold with $b_2^+(X)\geq 3$.  If
the SW-basic classes of $X$ are multiples of a class $K\in H^2(X;\ZZ)$,
then $X$ is abundant.
\end{cor}

\begin{proof}
  We consider the intersection form $Q$ to be a form on $H^2(X;\ZZ)$.  By
  Rochlin's theorem (see Theorem 1.2.29 in \cite{GompfStipsicz}), we know
  that $\sigma(X) \equiv 0\pmod{16}$ and so the classification
  \eqref{eq:classification} of forms yields $(H^2(X;\ZZ),Q)\cong 2k
  E_8\oplus l H$ with $k,l\in\ZZ$. A result of Furuta \cite[Theorem
  1.2.31]{GompfStipsicz} then implies that $l\geq 2|k| + 1$, where
  $k=\frac{1}{16}\sigma(X)$ and $l=\frac{1}{2}(b_2(X)-\sigma(X))$.  If
  $\sigma(X)=0$ and $b_2^+(X)\geq 3$, then $b_2(X)\geq 4$ and so $0\leq
  \rank(Q)-4$; if $\sigma(X)\neq 0$, then $k\neq 0$ and Furuta's theorem
  implies $l\geq 2$ and so again $|\sigma(Q)|\leq \rank(Q)-4$. Therefore,
  $K^\perp$ contains a hyperbolic sublattice by Lemma
  \ref{lem:OneBasicClass}.
\end{proof}

We turn to the more complicated case, where the intersection form is odd.

\begin{lem}
\label{lem:CrudeKexists}
Suppose $(L,Q)$ is an indefinite, integral, unimodular lattice.  If $Q$ is
odd with $b_2^+\geq 5$ and $b_2^-\geq 3$, and $\kappa\in L$ is characteristic,
then $\kappa^\perp$ contains a hyperbolic sublattice.
\end{lem}

\begin{proof}
  We may assume without loss that $\kappa$ is primitive for, if not, write
  $\kappa =d\kappa'$ where $d\in\ZZ$ and $\kappa'\in L$ is primitive. Then
  $Q(\kappa,x)=dQ(\kappa',x)\equiv Q(x,x)\pmod{2}$ for all $x\in L$ since
  $\kappa$ is characteristic. If $d$ were even we would have $Q(x,x)\equiv
  0\pmod{2}$ for all $x$, contradicting our hypothesis that $Q$ is an odd
  form. Hence, $d$ is odd and $dQ(\kappa',x)\equiv Q(\kappa',x)\equiv
  Q(x,x)\pmod{2}$ for all $x$ and thus $\kappa'$ is also
  characteristic. Then a sublattice $H$ is orthogonal to $\kappa$ if and
  only if it is orthogonal to $\kappa'$.
  
  The classification of indefinite, integral, unimodular forms shows that,
  because $Q$ is odd, we can find a basis $\{e_i,f_j\}$ for $L$ for which
  $e_i^2=1$, $f_i^2=-1$, and
$$
(L,Q)
\cong
\left(\mathop{\oplus}\limits_{i=1}^{b_2^+}\ZZ e_i\right) 
\oplus 
\left(\mathop{\oplus}\limits_{j=1}^{b_2^-}\ZZ f_j\right).
$$
For odd integers $a_1,a_2,a_3,a_4,c_1$ yet to be determined, choose
$$
\ell 
= 
a_1e_1+a_2e_2+a_3e_3+a_4e_4+\sum_{i=5}^{b_2^+}e_i
+ c_1f_1+\sum_{j=2}^{b_2^-}f_j \in L.
$$
The element $\ell$ is primitive since at least one basis coefficient is
equal to one. By hypothesis, $b_2^+\geq 5$ and $b_2^-\geq 3$, so we can define
$$
H = \ZZ(e_5 + f_2) + \ZZ(e_5 + f_3)
$$
and observe that $(H, Q|_H)$ is hyperbolic and orthogonal to $\ell$. Then
\begin{align*}
Q(\ell,\ell)
&= 
a_1^2+a_2^2+a_3^2+a_4^2 + (b_2^+ - 4) - c_1^2 - (b_2^- - 1)
\\
&=
a_1^2+a_2^2+a_3^2+a_4^2 + \sigma(Q) - c_1^2 - 3,
\end{align*}
and so
\begin{equation}
  \label{eq:SquareEllformula}
Q(\ell,\ell)- \sigma(Q) + c_1^2 + 3
=
a_1^2+a_2^2+a_3^2+a_4^2.
\end{equation}
Since the coefficients of $\ell$ are odd, we have $Q(\ell,e_i)\equiv
Q(\ell,f_j) \equiv 1 \pmod{2}$ for all $i,j$; thus $\ell$ is
characteristic and we have $\ell^2\equiv \sigma(Q)\pmod{8}$.
As $c_1$ is odd, so we can write $c_1 = 2u+1$ for some $u\in \ZZ$. The
left-hand side of the preceding equation therefore yields
$$
Q(\ell,\ell)- \sigma(Q) + c_1^2 + 3
\equiv 
(2u+1)^2+3
\equiv
4 + 4u(u+1)
\equiv 
4 \pmod{8}.
$$
Now any positive integer which is congruent to $4\pmod{8}$ can be
written as the sum of four odd squares (see Lemma \ref{lem:FourOddSquare}).
So we select any odd integer $c_1$ for which $Q(\kappa,\kappa) - \sigma(Q)
+ c_1^2 + 3 > 0$ and then choose odd integers $a_1,a_2,a_3,a_4$ so that
\begin{equation}
  \label{eq:SquareKappaformula}
a_1^2+a_2^2+a_3^2+a_4^2
=
Q(\kappa,\kappa) - \sigma(Q) + c_1^2 + 3.
\end{equation}
Therefore equations \eqref{eq:SquareEllformula} and
\eqref{eq:SquareKappaformula} give 
$$
Q(\ell,\ell)= Q(\kappa,\kappa).
$$
A result of Wall, (see \cite[Proposition 1.2.18]{GompfStipsicz}) implies
that the orthogonal group of $(L,Q)$ acts transitively on the primitive,
characteristic elements with a given square.  Hence, we can find an
orthogonal automorphism $A$ of $(L,Q)$ with $A\ell=\ka$.  Then $A(H)$ is a
hyperbolic sublattice of $L$ which is orthogonal to $\ka$.
\end{proof}

Lagrange's theorem tells us that every positive integer $k$ is the sum of four
integral squares \cite[Theorem 20.5]{HardyWright}. While the following
refinement of this result must surely be well known, we cannot find
references to it in standard texts on elementary number theory, so we include
a proof here which was generously supplied to us by Adebisi  Agboola
\cite{AgboolaJune2000}. Andr\'as Stipsicz has pointed out to us that 
Lemma \ref{lem:FourOddSquare} was used by Dieter Kotschick in
\cite{KotschickAlmostComplex} to show that every finitely presentable group
is the fundamental group of a closed, almost complex four-manifold.

\begin{lem}
\label{lem:FourOddSquare}
A positive integer is the sum of four odd squares if and only if it is
congruent to $4$ modulo $8$.
\end{lem}

\begin{proof}
  Let $k$ be a positive integer.  The set of squares modulo $8$ is $0,1,4$.
  Hence, if $k$ is the sum of four odd squares, then we must have
  $k\equiv 4\pmod{8}$.
  
  Conversely, suppose that $k\equiv 4\pmod{8}$. By Lagrange's theorem we
  can write
$$
    k = w^2 + x^2 + y^2 + z^2.
$$
Since $k$ is even, we must have one of the following possibilities (up to 
rearranging the terms on the right-hand side):
\alphenumi
\begin{enumerate}
\item $w, x$ are even, and $y, z$ are odd.
\item $w, x, y, z$ are all even.
\item $w, x, y, z$ are all odd (as desired).
\end{enumerate}
Now since $k\equiv 4\pmod{8}$, it is not hard to check that case (a) cannot
happen. In fact, it is easy to see that if case (a) is true, then
$k\equiv 2\text{ or }6\pmod{8}$.

To see that case (c) can occur, write $k = 8m+4$, where $m$ is a
non-negative integer, so
$$
          k = 1^2 + (8m+3).
$$
Legendre's theorem tells us that a positive integer is the sum of three
integral squares if and only if it cannot be written in the form
$4^a(8b+7)$, for some $a,b$ (see \cite[\S 20.10]{HardyWright} for the
statement and \cite{Ankeny} for a proof).  Now $1$ is odd, and since $8m+3$
cannot be of the form $4^a(8b+7)$, we can express $8m+3$ as a sum of three
squares. So we can write
$$
          k = 1^2 + x^2 + y^2 + z^2.
$$
Since $1$ is odd, and we have ruled out case (a) above, it follows that
$x, y, z$ are all odd.
\end{proof}

It remains to consider an ``unstable range'', where $b_2^+(Q)$ and
$b_2^-(Q)$ are small. Define $\chi(Q)=2+b_2(Q)$ and note that
$2\chi(Q)+3\sigma(Q) = 4 + 5b_2^+(Q) - b_2^-(Q)$.  

\begin{lem}
\label{lem:UnstableRangeKexists}
Suppose $(L,Q)$ is an integral, indefinite, unimodular lattice and that
$\kappa \in L$. Assume $Q$ is odd and $\kappa$ is characteristic. If one of
the following hold, then $\kappa^\perp$ contains a hyperbolic sublattice:
\alphenumi
\begin{enumerate}
\item
$b_2^+ = 3$ and $b_2^-\geq 5$,
\item
$b_2^+ = 3$ and $2\leq b_2^-\leq 4$, and $\kappa$ has square
$2\chi(Q)+3\sigma(Q)$. 
\end{enumerate}
\end{lem}

\begin{proof}
Interchanging the role of $b_2^+$ and $b_2^-$ in the proof of 
Lemma \ref{lem:CrudeKexists} takes care of case (a). Thus we need only
consider case (b).  

Continuing the notation of the proof of Lemma \ref{lem:CrudeKexists}, we
choose
$$
\ell = 3e_1 + 3e_2 + e_3 + \sum_{j=1}^{b_2^-}f_j \in L.
$$
Plainly, $\ell$ is characteristic, primitive, and is orthogonal to the
hyperbolic sublattice $H=\ZZ(e_3+f_1) + \ZZ(e_3+f_2)$ of $(L,Q)$, while 
(as $2\chi(Q)+3\sigma(Q) = 19-b_2^-$)
$$
Q(\ell,\ell)
=
19-b_2^-
=
Q(\kappa,\kappa).
$$
Since $Q$ is odd, we cannot have $b_2^-=3$ (which would give
$Q(\kappa,\kappa) = 16$) so $b_2^-=2$ or $b_2^-=4$, which gives
$Q(\kappa,\kappa) = 17$ or $15$, respectively, neither of which is
divisible by $d^2$, $d\in\ZZ$, unless $d=\pm 1$. Thus $\kappa$ is primitive.

Just as in the proof of Lemma \ref{lem:CrudeKexists}, Wall's theorem
implies that we can find an orthogonal automorphism $A$ of $(L,Q)$ with
$A\ell=\ka$, since $\kappa, \ell$ are both primitive, characteristic, and have
equal square. Then $A(H)$ is a hyperbolic sublattice of $\kappa^\perp$. This
takes care of case (b) and completes the proof of the lemma.
\end{proof}

Lemmas \ref{lem:CrudeKexists} and \ref{lem:UnstableRangeKexists} thus
yield:

\begin{cor}
\label{cor:OneBasicClassOdd}
Let $X$ be a simply-connected four-manifold having odd intersection form
$Q_X$, with $b_2^+(X)\geq 5$ and $b_2^-(X)\geq 3$ or $b_2^+(X)=3$ and
$b_2^-(X)\geq 2$.  Suppose that the SW-basic classes of $X$ are integer
multiples of a class $K\in H^2(X;\ZZ)$, where $K$ is characteristic.  Then
$X$ is abundant.
\end{cor}

Combining Corollaries
\ref{cor:OneBasicClassEven} and \ref{cor:OneBasicClassOdd} yields:

\begin{prop}
\label{prop:OneBasicClass}
Let $X$ be a simply-connected four-manifold.  Suppose the SW-basic classes
of $X$ are multiples of a class $K\in H^2(X;\ZZ)$, where $K$ is
characteristic.  If any one of the following hold, then $X$ is abundant:
\alphenumi
\begin{enumerate}
\item
$Q_X$ is even with $b_2^+(X)\geq 3$,
\item
$Q_X$ is odd with $b_2^+(X)\geq 5$ and $b_2^-(X)\geq 3$,
\item
$Q_X$ is odd with $b_2^+(X)=3$ and $b_2^-(X)\geq 5$, 
\item
$Q_X$ is odd with $b_2^+(X)=3$ and $2\leq b_2^-(X)\leq 4$, and $K^2 =
2\chi(X)+3\sigma(X)$. 
\end{enumerate}
\end{prop}

\subsection{Compact, complex algebraic, simply connected surfaces}
\label{subsec:CptCplxSimpConnSurface}
We now combine the results of the preceding subsection to prove the
principal result of this appendix:

\begin{proof}[Proof of Theorem \ref{thm:Abundance}]
  Let $X$ be a compact, complex algebraic, simply connected surface with
  $b_2^+(X) \geq 3$. Suppose $X$ is minimal.  The Enriques-Kodaira
  classification then implies that $X$ is one of the following (see
  \cite[\S 3.4]{GompfStipsicz}):
\alphenumi
\begin{enumerate}
\item
A $K3$ surface,
\item
An elliptic surface,
\item
A surface of general type.
\end{enumerate}
The cases where $X$ is diffeomorphic to $\CC\PP^2$,
$\CC\PP^2\#\overline{\CC\PP}^2$, or $\CC\PP^1\times\CC\PP^1$ (when $X$ has
Kodaira dimension $-\8$, see \cite[p.  88 \& Theorem
3.4.13]{GompfStipsicz}) are eliminated by our requirement that
$b_2^+(X)\geq 3$.

If $X$ is elliptic, then it is diffeomorphic to $E(n)_{p,q}$, for some
$n,p,q \in\NN$, $p\leq q$, $\gcd(p,q)=1$ (see Theorems 3.4.12 and 3.4.13
in\cite{GompfStipsicz}); if $X$ is a $K3$ surface, then it is diffeomorphic
to the surface $E(2)$ (see Theorem 3.4.9 in \cite{GompfStipsicz}), which is
included in this family as $E(n)_{1,1}=E(n)$ (see
\cite[p. 83]{GompfStipsicz} for the construction of this family).

Let $f$ be the homology class of a regular fiber of $E(n)_{p,q}$ and
observe that $f_{pq}=\frac{1}{pq}f$ is a primitive, integral homology
class.  According to \cite{FSRationalBlowdown} (see \cite[Theorem
3.3.6]{GompfStipsicz}), the SW-basic classes of $E(n)_{p,q}$ are multiples
of the Poincar\'e dual $\PD[f_{pq}]$. Thus, $X$ is abundant by Proposition
\ref{prop:OneBasicClass} (cases (a), (b), and (c))
and the observation that $b_2^+(E(n)_{p,q})=2n-1$
and $b_2^-(E(n)_{p,q})=10n-1$ (see \cite[Lemma 3.3.4]{GompfStipsicz}), so
$b_2^+(X)\geq 3$ and $b_2^-(X)\geq 19$ for $n\ge 2$.

If $X$ is a minimal algebraic surface of general type, then its SW-basic
classes are $\pm K_X$ by \cite{FriedmanDonSW}, where $K_X$ is the canonical
class.  If $Q_X$ is even, then $X$ is abundant by Proposition
\ref{prop:OneBasicClass}.  If $Q_X$ is odd, the Bogomolov-Miyaoka-Yau
inequality, $c_1^2(X)\le 3c_2(X)$ (see \cite[Theorem VII.1.1 (iii)]{BPV}),
implies that
$$
3\sigma \leq \chi,
$$
since $c_1^2(X)=K_X^2=2\chi+3\sigma$ and $c_2(X)=\chi$. Thus
\begin{equation}
  \label{eq:BMYbpm}
  b_2^+(X)\le 2b_2^-(X)+1.
\end{equation}
Equality in \eqref{eq:BMYbpm}, or $c_1^2(X) = 3c_2(X)$, holds only if
the universal covering space of $X$ is the closed unit ball in $\CC^2$ by
\cite[Theorem 7.2.24]{GompfStipsicz} (or see \cite[Corollary I.15.5]{BPV} and
the discussion in \cite[p. 230]{BPV} or
\cite[Theorem 4]{Miyaoka}).  Since $X$ is simply-connected and closed by
hypothesis, we must have
$$
b_2^+(X) < 2 b_2^-(X)+1.
$$
If $b_2^+(X)=3$, then the preceding inequality yields $b_2^-(X)\ge 2$
while if $b_2^+(X)\ge 5$, it yields $b_2^-(X)\ge 3$.  Hence, for $Q_X$ odd,
$X$ is again abundant by Proposition \ref{prop:OneBasicClass}.

This proves the theorem for minimal surfaces.  If $X$ is an abundant,
smooth four-manifold, there is a hyperbolic sublattice $H \subset
H^2(X;\ZZ)$ such that $K$ is orthogonal to $H$ if $K$ is an SW-basic class.
The blow-up $X\#\overline{\CC\PP}^2$ is also abundant, since we may view
$H$ as a hyperbolic sublattice of $H^2(X\#\overline{\CC\PP}^2;\ZZ)$, all
SW-basic classes of $X\#\overline{\CC\PP}^2$ have the form $K\pm\PD[e]$
(see \cite{FSTurkish}), and such classes are again orthogonal to $H$.
Hence, if $X$ is abundant, all its blow-ups are abundant too.
This completes the proof of the theorem.
\end{proof}


\ifx\undefined\bysame
\newcommand{\bysame}{\leavevmode\hbox to3em{\hrulefill}\,}
\fi


\begin{thebibliography}{10}

\bibitem{AgboolaJune2000}
A.~Agboola, {\em Private communication}, June 2000.

\bibitem{Ankeny}
N.~C. Ankeny, {\em Sums of three squares}, Proc. Amer. Math. Soc. {\bf 8}
  (1957), 316--319.

\bibitem{ACGH}
E.~Arbarello, M.~Cornalba, P.~A. Griffiths, and J.~Harris, {\em Geometry of
  algebraic curves}, Springer, New York, 1985.

\bibitem{AtiyahK}
M.~F. Atiyah, {\em K-theory}, Addison-Wesley, New York, 1969.

\bibitem{AHS}
M.~F. Atiyah, N.~J. Hitchin, and I.~M. Singer, {\em Self-duality in
  four-dimensional {R}iemannian geometry}, Proc. Royal Soc. London {\bf A 362}
  (1978), 425--461.

\bibitem{AS4}
M.~F. Atiyah and I.~M. Singer, {\em Index of elliptic operators. {IV}}, Ann.
  Math. {\bf 93} (1971), 119--138.

\bibitem{BPV}
W.~Barth, C.~Peters, and A.~Van de~Ven, {\em Compact complex surfaces},
  Springer, New York, 1984.

\bibitem{BerlineGetzlerVergne}
N.~Berline, E.~Getzler, and M.~Vergne, {\em Heat kernels and {D}irac
  operators}, Springer-Verlag, Berlin, 1992.

\bibitem{BoossBleecker}
B.~Booss and D.~D. Bleecker, {\em Topology and analysis: The {A}tiyah-{S}inger
  index formula and gauge-theoretic physics}, Springer, New York, 1985.

\bibitem{Bredon}
G.~Bredon, {\em Introduction to compact transformation groups}, Academic Press,
  New York, 1972.

\bibitem{Brussee}
R.~Brussee, {\em The canonical class and the {$C^\infty$}-properties of
  {K\"a}hler surfaces}, New York J. Math. (electronic) {\bf 2} (1996),
  103--146, alg-geom/9503004.

\bibitem{DonApplic}
S.~K. Donaldson, {\em An application of gauge theory to four-dimensional
  topology}, J. Differential Geom. {\bf 18} (1983), 279--315.

\bibitem{DonPoly}
\bysame, {\em Polynomial invariants for smooth four-manifolds}, Topology {\bf
  29} (1990), 257--315.

\bibitem{DK}
S.~K. Donaldson and P.~B. Kronheimer, {\em The geometry of four-manifolds},
  Oxford Univ. Press, Oxford, 1990.

\bibitem{FeehanGenericMetric}
P.~M.~N. Feehan, {\em Generic metrics, irreducible rank-one {PU(2)} monopoles,
  and transversality}, Comm. Anal. Geom. {\bf 8} (2000), to appear,
  math.DG/9809001.

\bibitem{FKLM}
P.~M.~N. Feehan, P.~B. Kronheimer, T.~G. Leness, and T.~S. Mrowka, {\em {PU(2)}
  monopoles and a conjecture of {M}ari{\~n}o, {M}oore, and {P}eradze}, Math.
  Res. Lett. {\bf 6} (1999), 169--182, math.DG/9812125.

\bibitem{FL2}
P.~M.~N. Feehan and T.~G. Leness, {\em {PU(2)} monopoles. {II}: {H}ighest-level
  singularities and relations between four-manifold invariants}, dg-ga/9712005
  (v1).

\bibitem{FL2b}
\bysame, {\em {PU(2)} monopoles. {II}: {T}op-level {S}eiberg-{W}itten moduli
  spaces and {W}itten's conjecture in low degrees}, J. Reine Angew. Math., to
  appear, dg-ga/9712005.

\bibitem{FL3}
\bysame, {\em {PU(2)} monopoles. {III}: {E}xistence of gluing and obstruction
  maps}, submitted to a print journal; math.DG/9907107.

\bibitem{FL4}
\bysame, {\em {PU(2)} monopoles. {IV}: {S}urjectivity of gluing maps}, in
  preparation.

\bibitem{FLLevelOne}
\bysame, {\em {PU(2)} monopoles, level-one {S}eiberg-{W}itten moduli spaces,
  and {W}itten's conjecture in low degrees}, preprint, submitted to a print
  journal.

\bibitem{FL5}
\bysame, {\em {PU(2)} monopoles. {V}}, in preparation.

\bibitem{FLGeorgia}
\bysame, {\em {PU(2)} monopoles and relations between four-manifold
  invariants}, Topology Appl. {\bf 88} (1998), 111--145, dg-ga/9709022.

\bibitem{FL1}
\bysame, {\em {PU(2)} monopoles. {I}: {R}egularity, {U}hlenbeck compactness,
  and transversality}, J. Differential Geom. {\bf 49} (1998), 265--410,
  dg-ga/9710032.

\bibitem{FS1}
R.~Fintushel and R.~Stern, {\em {SO(3)}-connections and the topology of
  4-manifolds}, J. Differential Geom. {\bf 20} (1984), 523--539.

\bibitem{FSTurkish}
\bysame, {\em Immersed spheres in 4-manifolds and the immersed {T}hom
  conjecture}, Turkish J. Math. {\bf 19} (1995), 145--157.

\bibitem{FSBlowUp}
\bysame, {\em Blow-up formulas for {D}onaldson invariants}, Ann. of Math. {\bf
  143} (1996), 529--546, alg-geom/9405002.

\bibitem{FSRationalBlowdown}
\bysame, {\em Rational blowdowns of smooth $4$-manifolds}, J. Differential
  Geom. {\bf 46} (1997), 181--235, alg-geom/9505018.

\bibitem{FU}
D.~Freed and K.~K. Uhlenbeck, {\em Instantons and four-manifolds}, 2nd ed.,
  Springer, New York, 1991.

\bibitem{FriedmanDonSW}
R.~Friedman, {\em Donaldson and {S}eiberg-{W}itten invariants of algebraic
  surfaces}, Algebraic geometry (Santa Cruz, 1995), Amer. Math. Soc.,
  Providence, RI, 1997, alg-geom/9605006, pp.~85--100.

\bibitem{FrM}
R.~Friedman and J.~W. Morgan, {\em Smooth four-manifolds and complex surfaces},
  Springer, Berlin, 1994.

\bibitem{Furuta}
M.~Furuta, {\em The {$10/8$} conjecture}, unpublished.

\bibitem{GompfMrowka}
R.~E. Gompf and T.~S.Mrowka, {\em Irreducible four-manifolds need not be
  complex}, Ann. of Math. (2) {\bf 138} (1993), 61--111.

\bibitem{GompfStipsicz}
R.~E. Gompf and A.~I. Stipsicz, {\em $4$-manifolds and {K}irby calculus},
  American Mathematical Society, Providence, RI, 1999.

\bibitem{GorMacPh}
M.~Goresky and R.~MacPherson, {\em Stratified {M}orse theory}, Springer, New
  York, 1980.

\bibitem{HardyWright}
G.~H. Hardy and E.~M. Wright, {\em An introduction to the theory of numbers},
  fifth ed., The Clarendon Press--Oxford University Press, New York, 1979.

\bibitem{Hirsch}
M.~W. Hirsch, {\em Differential topology}, Springer, New York, 1976.

\bibitem{HirzebruchHopf}
F.~Hirzebruch and H.~Hopf, {\em Felder von {F}l{\"a}chenelementen in
  4-dimensionalen {M}annigfaltigkeiten}, Math. Ann. {\bf 136} (1958), 156--172.

\bibitem{Husemoller}
D.~Husemoller, {\em Fiber bundles}, second ed., Springer, New York, 1966.

\bibitem{Karoubi}
M.~Karoubi, {\em {$K$}-theory}, Springer, Berlin, 1978.

\bibitem{Koschorke}
U.~Koschorke, {\em Infinite-dimensional {K}-theory and characteristic classes
  of {F}redholm bundle maps}, Global Analysis (F.~E. Browder, ed.), Proc. Symp.
  Pure Math., vol. 15-I, Amer. Math. Soc., Providence, RI, 1970, pp.~95--133.

\bibitem{KotschickAlmostComplex}
D.~Kotschick, {\em All fundamental groups are almost complex}, Bull. London
  Math. Soc. {\bf 24} (1992), 377--378.

\bibitem{KMThom}
P.~B. Kronheimer and T.~S. Mrowka, {\em The genus of embedded surfaces in the
  projective plane}, Math. Res. Lett. {\bf 1} (1994), 797--808.

\bibitem{KMStructure}
\bysame, {\em Embedded surfaces and the structure of {D}onaldson's polynomial
  invariants}, J. Differential Geom. {\bf 43} (1995), 573--734.

\bibitem{KMContact}
\bysame, {\em Monopoles and contact structures}, Invent. Math. {\bf 130}
  (1997), 209--255.

\bibitem{Kuiper}
N.~H. Kuiper, {\em The homotopy type of the unitary group of {H}ilbert space},
  Topology {\bf 3} (1965), 19--30.

\bibitem{Kuranishi}
M.~Kuranishi, {\em New proof for the existence of locally complete families of
  complex structures}, Proc. Conf. on Complex Analysis (A.~Aeppli, E.~Calabi,
  and H.~R{\"o}hrl, eds.), Springer, New York, 1965, pp.~142--154.

\bibitem{Lang}
S.~Lang, {\em Differential and {R}iemannian manifolds}, third ed., Springer,
  New York, 1995.

\bibitem{LM}
H.~B. Lawson and M-L. Michelsohn, {\em Spin geometry}, Princeton Univ. Press,
  Princeton, NJ, 1988.

\bibitem{LiTian}
J.~Li and G.~Tian, {\em Virtual moduli cycles and {G}romov-{W}itten invariants
  of general symplectic manifolds}, Topics in symplectic $4$-manifolds (Irvine,
  CA, 1996), Internat. Press, Cambridge, MA, 1998, alg-geom/9608032,
  pp.~47--83.

\bibitem{LiLiu}
T-J. Li and A-K. Liu, {\em General wall-crossing formula}, Math. Res. Lett.
  {\bf 2} (1995), 797--810.

\bibitem{Mather}
J.~N. Mather, {\em Stratifications and mappings}, Dynamical systems (Proc.
  Sympos., Univ. Bahia, Salvador, 1971) (M.~M. Peixoto, ed.), Academic Press,
  New York, 1973, pp.~195--232.

\bibitem{MilnorStasheff}
J.~W. Milnor and J.~D. Stasheff, {\em Characteristic classes}, Princeton Univ.
  Press, Princeton, NJ, 1974.

\bibitem{Miyaoka}
Y.~Miyaoka, {\em Algebraic surfaces with positive indices}, Classification of
  algebraic and analytic manifolds (Katata, 1982) (K.~Ueno, ed.), Birkh\"auser,
  Boston, MA, 1983, pp.~281--301.

\bibitem{MooreWitten}
G.~Moore and E.~Witten, {\em Integration over the $u$-plane in {D}onaldson
  theory}, Adv. Theor. Math. Phys. {\bf 1} (1997), 298--387, hep-th/9709193.

\bibitem{MorganSWNotes}
J.~W. Morgan, {\em The {S}eiberg-{W}itten equations and applications to the
  topology of smooth four-manifolds}, Princeton Univ. Press, Princeton, NJ,
  1996.

\bibitem{MorganMrowkaPoly}
J.~W. Morgan and T.~S. Mrowka, {\em A note on {D}onaldson's polynomial
  invariants}, Internat. Math. Res. Notes {\bf 10} (1992), 223--230.

\bibitem{MMR}
J.~W. Morgan, T.~S. Mrowka, and D.~Ruberman, {\em {$L^2$} moduli spaces and a
  vanishing theorem for {D}onaldson polynomial invariants}, Internat. Press,
  Cambridge, MA, 1994.

\bibitem{MrowkaPrincetonMorseTalk}
T.~S. Mrowka, {\em Marston {M}orse memorial lectures}, Institute for Advanced
  Study, Princeton, NJ, April, 1999.

\bibitem{MrowkaOzsvathYu}
T.~S. Mrowka, P.~S. Ozsv{\'a}th, and B.~Yu, {\em Seiberg-{W}itten monopoles on
  {S}eifert fibered spaces}, Comm. Anal. Geom. {\bf 5} (1997), 685--791.

\bibitem{OTQuaternion}
C.~Okonek and A.~Teleman, {\em Quaternionic monopoles}, Comm. Math. Phys. {\bf
  180} (1996), 363--388, alg-geom/9505029.

\bibitem{OTWall}
\bysame, {\em Seiberg-{W}itten invariants for manifolds with {$b^+=1$} and the
  universal wall crossing formula}, Internat. J. Math. {\bf 7} (1996),
  811--832, alg-geom/9603003.

\bibitem{PTLectures}
V.~Y. Pidstrigatch, Lectures at the {N}ewton {I}nstitute in {D}ecember 1994,
  {O}berwolfach in {M}ay 1996, and the {N}ewton {I}nstitute in {N}ovember 1996.

\bibitem{PTLocal}
V.~Y. Pidstrigatch and A.~N. Tyurin, {\em Localisation of {D}onaldson
  invariants along the {S}eiberg-{W}itten classes}, dg-ga/9507004.

\bibitem{RuanGafa}
Y.~Ruan, {\em Symplectic topology and extremal rays}, Geom. Funct. Anal. {\bf
  3} (1993), 395--430.

\bibitem{RuanSW}
\bysame, {\em Virtual neighborhoods and the monopole equations}, Topics in
  symplectic $4$-manifolds (Irvine, CA, 1996), Internat. Press, Cambridge, MA,
  1998, alg-geom/9611021., pp.~101--116.

\bibitem{RuanGW}
\bysame, {\em Virtual neighborhoods and pseudo-holomorphic curves}, Turkish J.
  Math. {\bf 23} (1999), 161--231, alg-geom/9611021.

\bibitem{SalamonSWBook}
D.~Salamon, {\em Spin geometry and {S}eiberg-{W}itten invariants},
  Birkh{\"a}user, Boston, to appear.

\bibitem{SiebertGW}
B.~Siebert, {\em Symplectic {G}romov-{W}itten invariants}, New trends in
  algebraic geometry (Warwick, 1996), Cambridge Univ. Press, Cambridge, 1999,
  alg-geom/9608005, pp.~375--424.

\bibitem{Spanier}
E.~H. Spanier, {\em Algebraic topology}, Springer, New York, 1966.

\bibitem{TauIndef}
C.~H. Taubes, {\em Self-dual connections on 4-manifolds with indefinite
  intersection matrix}, J. Differential Geom. {\bf 19} (1984), 517--560.

\bibitem{TauStable}
\bysame, {\em The stable topology of self-dual moduli spaces}, J. Differential
  Geom. {\bf 29} (1989), 162--230.

\bibitem{TelemanGenericMetric}
A.~Teleman, {\em Moduli spaces of {PU(2)}-monopoles}, Asian J. Math. {\bf 4}
  (2000), 391--435, math.DG/9906163.

\bibitem{WallUnimodQuadForms}
C.~T.~C. Wall, {\em On the orthogonal groups of unimodular quadratic forms},
  Math. Ann. {\bf 147} (1962), 328--338.

\bibitem{Witten}
E.~Witten, {\em Monopoles and four-manifolds}, Math. Res. Lett. {\bf 1} (1994),
  769--796, hep-th/9411102.

\end{thebibliography}
\end{document}